\documentclass[11pt]{article}

\usepackage{enumitem}

\usepackage[backref=false,colorlinks=true, linkcolor=blue, citecolor=blue, urlcolor=blue, pdfborder={0 0 0}]{hyperref}

\usepackage[left=1in,top=1in,right=1in,bottom=1in,letterpaper]{geometry}

\usepackage{amssymb,amsmath, amsthm}
\usepackage{mathtools}

\usepackage{cancel}

\allowdisplaybreaks

\usepackage{indentfirst}

\usepackage{inputenc}
\usepackage{graphicx}
\usepackage{epstopdf}

\usepackage{url}

\usepackage{appendix}

\usepackage{matlab-prettifier}
\usepackage{listings}

\usepackage{multicol}

\usepackage{hyperref}
\usepackage{amssymb}
\usepackage{amsbsy}
\usepackage{amsmath}
\usepackage{amsthm}
\usepackage{graphicx}
\usepackage{subfigure}
\usepackage{calc}
\usepackage{color}
\usepackage{multirow}
\usepackage{multicol}
\usepackage{booktabs}
\usepackage{wasysym}
\usepackage{lmodern} 
\usepackage{textcomp} 
\usepackage{wasysym} 
\usepackage{dashbox} 
\usepackage{mathtools} 
\usepackage{tabto}  
\usepackage{algorithm}
\usepackage{algpseudocode}
\usepackage{esdiff} 

\usepackage{cancel}

\usepackage{empheq}

\usepackage{enumitem}

\usepackage{etoolbox}
\makeatletter
\patchcmd{\@makefntext}{\insertfootnotetext{#1}}{\insertfootnotetext{\scriptsize#1}}{}{}
\makeatother

\usepackage[normalem]{ulem} 

\newcommand{\cut}[1]{{}}

\usepackage{pifont}
\makeatletter
  \@ifclassloaded{article} 
    {}
    {
      \usepackage{natbib}
      \newcommand{\citesmall}[2][]
      {
         \ifthenelse{\equal{#1}{}}
         {\begin{small}\citep{#2}\end{small}}
         {\begin{small}\citep[#1]{#2}\end{small}}
      }
      \newcommand{\citefootnote}[2][]
      {
        \footnote{
          \ifthenelse{\equal{#1}{}}
          {\begin{small}\citep{#2}\end{small}}
          {\begin{small}\citep[#1]{#2}\end{small}}
        }
      }

      \mode<presentation>
      {
        \setbeamertemplate{footline}[frame number]

        \setbeamertemplate{background canvas}[vertical shading][bottom=white!10,top=white!10]
        \usetheme{default}

        \setbeamertemplate{sections/subsections in toc}[sections numbered]

        \setbeamertemplate{navigation symbols}{}

        \setbeamertemplate{itemize items}[circle]
        \setbeamercolor{itemize item}{fg=black}
        \setbeamercolor{enumerate item}{fg=black}
        \setbeamercolor{itemize subitem}{fg=black}
        \setbeamercolor{enumerate subitem}{fg=black}

        \setbeamertemplate{itemize/enumerate body begin}{\normalsize}
        \setbeamertemplate{itemize/enumerate subbody begin}{\normalsize}
        \setbeamertemplate{itemize/enumerate subsubbody begin}{\normalsize}

        \setbeamertemplate{blocks}[rounded][shadow=false]
        \setbeamercolor{block title}{fg=black,bg=gray!40}
        \setbeamercolor{block body}{fg=black,bg=gray!10}

        \setbeamertemplate{frametitle}[default][center]
        \setbeamercolor{frametitle}{fg=black}
        \setbeamerfont{frametitle}{shape=\bfseries}

        \linespread{1.2}

        \setlength{\parskip}{\smallskipamount}

        \AtBeginNote{%
        }

        \setbeamertemplate{note page}[compressed]
        \setbeamercolor{note page}{bg=white}
      }
    } 
\makeatother%

\newcommand{\cO}{{\mathcal{O}}}
\newcommand{\cP}{{\mathcal{P}}}

\newcommand{\RR}{\mathbb{R}}
\newcommand{\CC}{\mathbb{C}}
\newcommand{\NN}{\mathbb{N}}

\DeclareMathOperator*{\argmin}{arg\,min}

\newcommand{\vvvert}{{\vert\kern-0.25ex\vert\kern-0.25ex\vert}}

\makeatletter
  \@ifclassloaded{article} 
    {

    }
    {} 
\makeatother%

\begin{document}
\title{Fast Partial Fourier Transforms for Large-Scale Ptychography}
\author{Ricardo Parada\footnotemark[1]$\:\:$\footnotemark[2] \and
Samy Wu Fung\footnotemark[3]
\and
Stanley Osher\footnotemark[1]}
\date{\today}
\maketitle

\renewcommand{\thefootnote}{\fnsymbol{footnote}}
\footnotetext[1]{Department of Mathematics, University of California, Los Angeles}
\footnotetext[3]{Department of Applied Mathematics and Statistics, Colorado School of Mines}
\footnotetext[2]{Corresponding author. Email: rparadaumanzor@gmail.com}

\begin{abstract}
Ptychography is a popular imaging technique that combines diffractive imaging with scanning microscopy. The technique consists of a coherent beam that is scanned across an object in a series of overlapping positions, leading to reliable and improved reconstructions. Ptychographic microscopes allow for large fields to be imaged at high resolution at additional computational expense. In this work, we explore the use of the fast Partial Fourier Transforms (PFTs), which efficiently compute Fourier coefficients corresponding to low frequencies. The core idea is to use the PFT in a plug-and-play manner to warm-start existing ptychography algorithms such as the ptychographic iterative engine (PIE). This approach reduces the computational budget required to solve the ptychography problem. Our numerical results show that our scheme accelerates the convergence of traditional solvers without sacrificing quality of reconstruction.
\end{abstract}

\section{Introduction}
\label{sec: intro}
\begin{figure}[t]
    \centering
    \includegraphics[width = 0.7\textwidth]{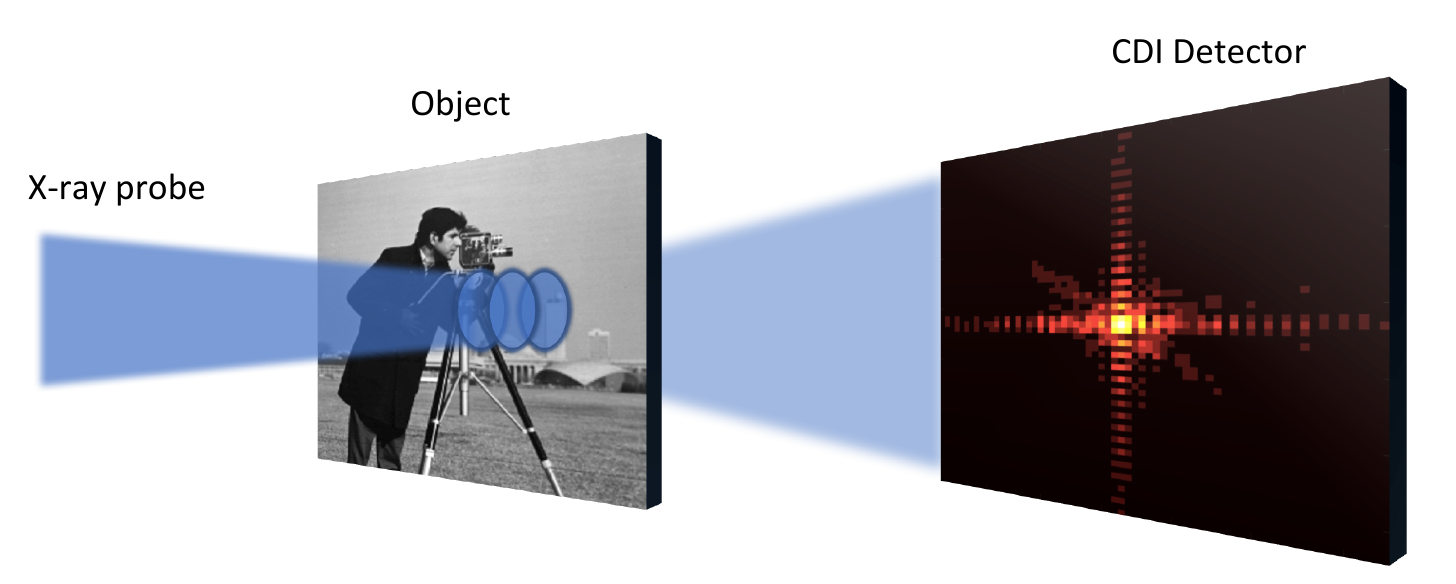}
    \caption{\small{Illustration of the ptychography experiment with three overlapping scans.}}
    \label{fig: ptycho_illustration}
\end{figure}

Ptychography is a coherent diffraction imaging (CDI) technique used across various fields including materials science~\cite{holler2017high,hoppe2013high,pelz2014fly}, biology~\cite{marrison2013ptychography,suzuki2016dark}, and x-ray crystallography~\cite{de2016ptychographic}. Originally developed to enhance resolution in electron or x-ray microscopy, ptychography replaces single-element detectors with two-dimensional array detectors and integrates diffractive imaging with scanning microscopy. 
In this process, a coherent beam is systematically moved over an object in a pattern of overlapping positions, linking information across successive diffraction patterns (refer to Fig.~\ref{fig: ptycho_illustration}). In conventional single-pattern CDI, the application of appropriate finite support constraints is essential for the effectiveness of standard algorithms such as Error Reduction~\cite{gerchberg1972practical}, Hybrid Input Output (HIO)~\cite{fienup1982phase}, gradient-based algorithms~\cite{candes2015phase}, Relaxed Averaged Alternating Reflections~\cite{luke2004relaxed}, and Saddle Point Optimization~\cite{marchesini2007phase,tripathi2015visualizing}. 
In ptychography, prior knowledge of scanning positions inherently provides these constraints, resulting in methods that are quicker and more reliable than those used in single-pattern CDI.

The extensive use of ptychography has spurred significant research into its reconstruction techniques. The extended Ptychographic Iterative Engine (ePIE)~\cite{maiden2009improved} is one of the most favored methods, involving alternating projections onto non-convex modulus constraint sets, and has become increasingly popular in the optics field. 
Mathematically, PIE functions as a projected steepest descent algorithm applied to a specific objective function. 
Other established methods include conventional gradient-based techniques and the Wirtinger Flow~\cite{candes2015phase,xu2018accelerated}, which employs a spectral method for initial setup. 
Another well-known method, called PhaseLift~\cite{candes2015phaseReview}, reformulates the phase retrieval challenge as a convex optimization problem, but requires solving for a significantly larger number of unknowns, becoming impractical for large-scale applications. 
Despite considerable efforts, solving large-scale ptychographic phase retrieval remains a formidable challenge due to the complexity of handling high resolution, high overlaps, and small scanning beams, all of which contribute to vast amounts of intricate data.

\subsection{Our Contribution}
\label{subsec: contribution}
\begin{figure}[t]
    \small
    \centering
    \begin{tabular}{cc}
        FFT $(512 \times 512)$ & PFT $(128 \times 128)$
        \\
        \includegraphics[width=0.35\textwidth]{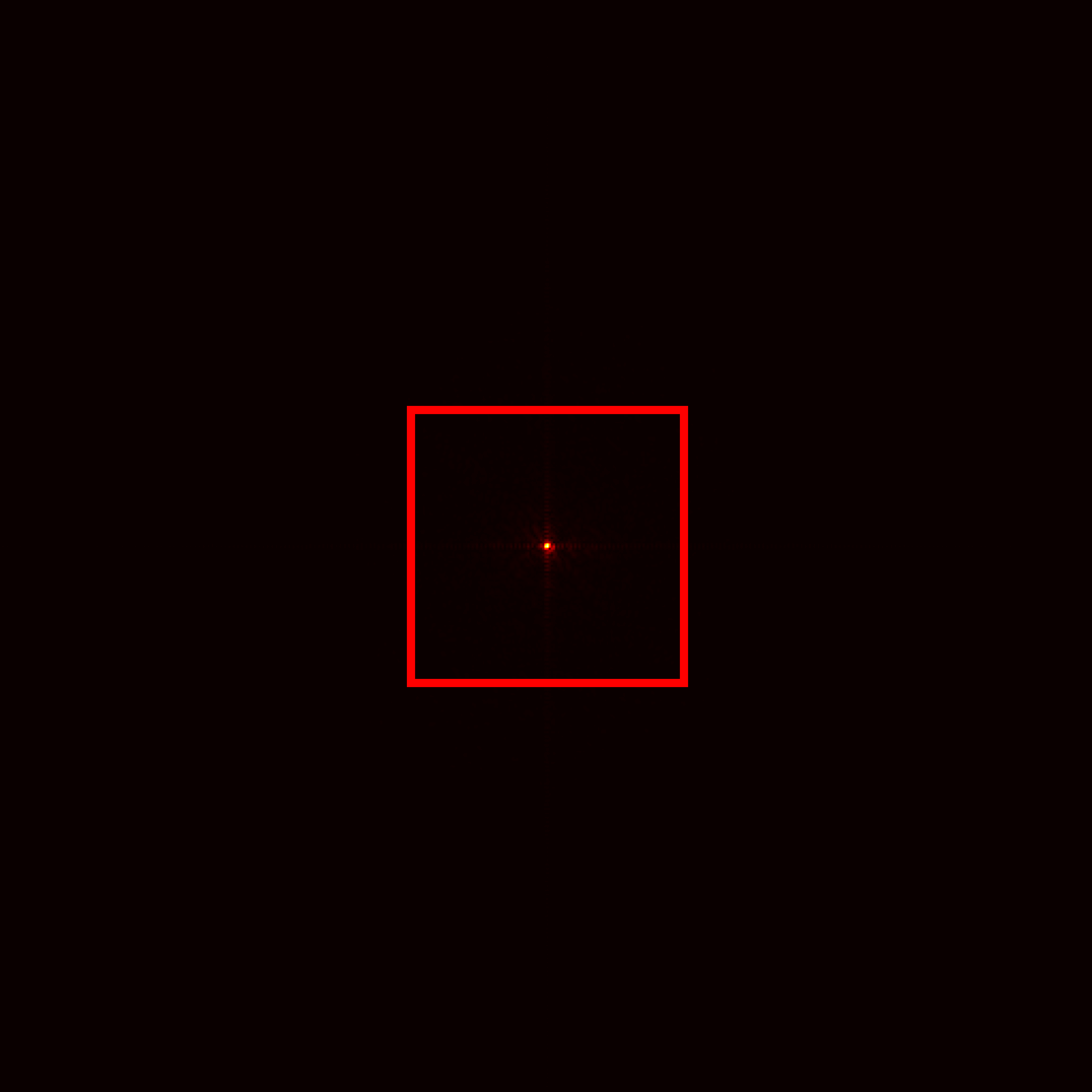}
        &
        \includegraphics[width=0.35\textwidth]{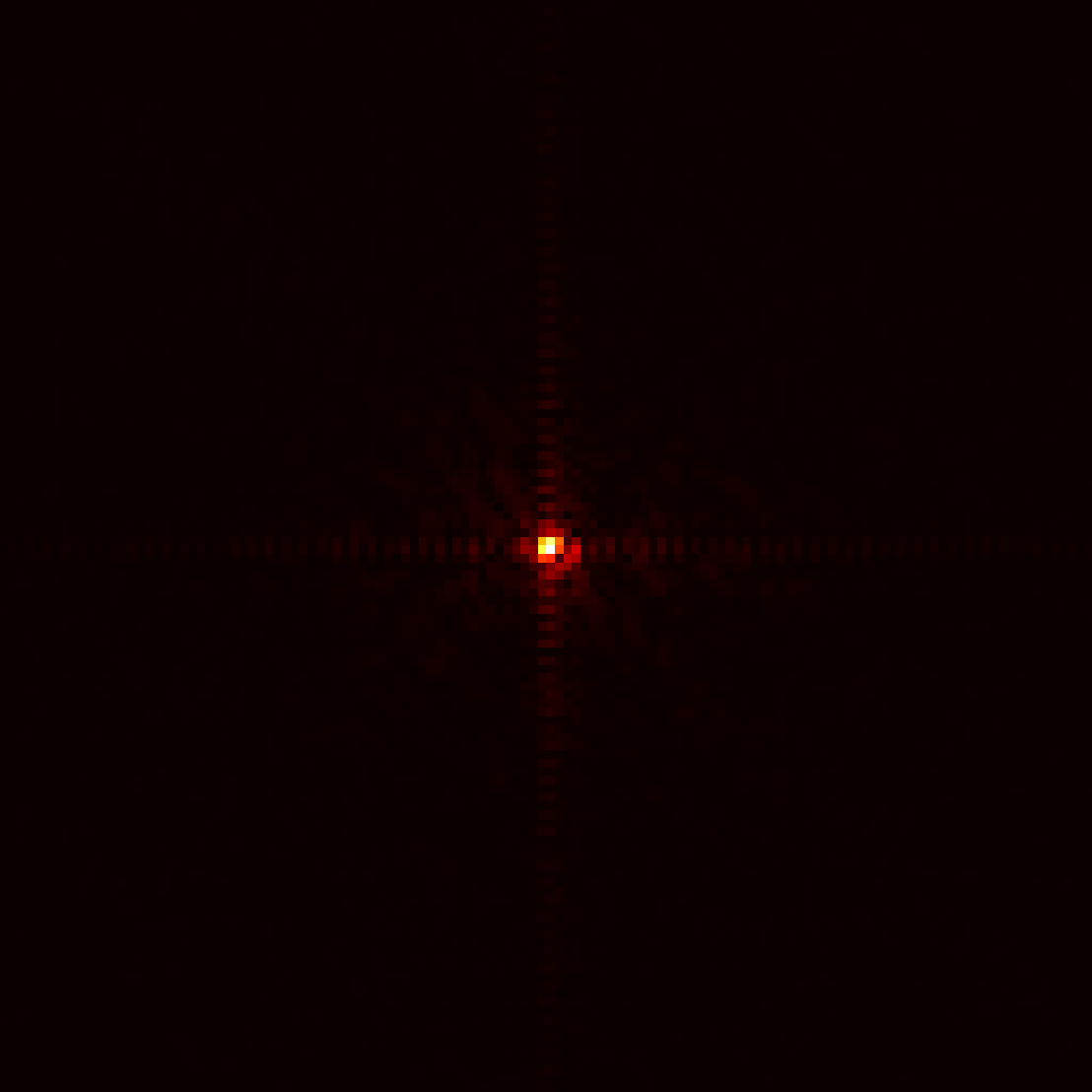}
    \end{tabular}
    \caption{\small{Illustration of coefficients computed by the PFT. On the left, we show the full FFT applied to an image of size $512 \times 512$. The red square shows the frequencies that the PFT computes on the right \emph{without} requiring one to take the FFT of the original image and then cropping as one might naively attempt.}}
    \label{fig: pft_crop}
\end{figure}

In this work, we focus on large scale ptychography problems where even applying the FFT is considered computationally taxing. 
In particular, we explore the use of recent work on fast partial Fourier transforms (PFTs)~\cite{park2021fast} as a warm up mechanism within the ePIE algorithm to accelerate convergence and time-to-solution of large-scale ptychographic phase retrieval problems. 
See Figure~\ref{fig: pft_crop} for an illustration of the PFT and Section~\ref{subsec: pft} for more details.
The core idea is to let the PFT-based ePIE algorithm capture the large features from the low frequencies in the initial iterations of ePIE, followed by standard FFT-based ePIE to capture the fine details of the reconstruction.
Our experiments show that including the PFT within existing algorithms such as the ePIE 1) does not reduce the quality of the reconstruction, and 2) helps accelerate convergence by improving runtimes for large-scale problems.
Importantly, in order to be able to use the PFT in gradient-based algorithms, it is necessary to be able to differentiate through the PFT operator. 
To this end, we also provide a PyTorch~\cite{paszke2019pytorch} implementation of the PFT in order for users to be able to differentiate through the operator using automatic differentiation (AD)~\cite{van2018automatic}. 

\section{Ptychography Background}
\label{sec: ptychography_background}
Let $z \in  \mathbb{C}^{n}$ be the object of interest and $d_j \in  \mathbb{R}^{m}$ be the observed data (or intensities) measured from the $j^{th}$ probe, where $n$ and $m$ are the dimensions of the vectorized object and data resolution images, respectively. 
A ptychography experiment is modeled by
\begin{equation}
    d_j = | \mathcal{F} (\omega \odot Q_j z) |, \quad j=1,\ldots,N,
\end{equation}
where $\mathcal{F} \in \mathbb{C}^{m \times m}$ is the discrete Fourier operator, $\omega \in \mathbb{R}^m$ is the localized probe, $z_j$ is the object after being scanned, $\odot$ represents the Hadamard (or element-wise) product, and $Q_j \in \mathbb{R}^{m \times n}$ is a matrix with
binary elements extracting a patch (with the index $j$ and size $m$) from the entire sample.
The blind ptychography problem (BPP) is an inverse problem and can be stated as:
    \begin{equation}
    \text{Find $\omega$ and $z$ such that $d_j = | \mathcal{F} (\omega \odot Q_j z) |$ $\quad$ for  $\quad j=1,\ldots,N$.}
    \label{eq: BPP_problem}
    \end{equation}
In practice, the observations $d_j$ can contain noise, usually Poisson, as is common in most inverse problems.
In our experiments, we will also consider the (non-blind) ptychography problem where the probes are known. This model can be realistic in some settings where domain scientists know the probing mechanism a priori. We also note that the above is a specific type of ptychography based on CDI; however other forms of ptychography exist, e.g., Fourier Ptychography~\cite{zheng2021concept}, and frequency-resolved optical gating (FROG)~\cite{trebino2012frequency}. While the techniques discussed in this work can be applied to all of them, we will only discuss CDI-based ptychography for ease of presentation.

\subsection{The Ptychographic Iterative Engine}
As previously stated, the ePIE algorithm is perhaps the most widely used algorithm in practice due to its simplicity. It can be viewed as an alternating projection algorithm onto non-convex modulus constraint sets that solve following feasibility problem. Let the $j^{th}$ measurement constraint set and its corresponding projection operator be denoted by
\begin{equation}
    \mathcal{M}_j^\omega = \left\{z \in \mathbb{C}^n : | \mathcal{F}(\omega \odot Q_j z)| = d_j\right\} \quad \text{ and } \quad \mathcal{P}_{\mathcal{M}_j^\omega}(z) = \mathcal{F}^{-1} \left[ d_j \odot \exp(i \theta(\mathcal{F}(\omega \odot Q_j z))) \right],
\end{equation}
where $\theta \colon \mathbb{C} \to [-\pi, \pi]$ returns the argument of a complex number and is applied element-wise. 
The ePIE algorithm generates an approximation to problem by solving the following optimization problem
\begin{equation}
    \min_{\omega, z} \; \Phi(\omega, z) = \frac{1}{N}\sum_{j=1}^N \left\| \omega \odot Q_j z - \mathcal{P}_{\mathcal{M}_j^\omega}(z) \right\|^2.
\end{equation}
In particular, the ePIE iterates are given by 
\begin{align}
    z_j^{k+1} &= z_j^k - \beta \bar{\omega}^k \left( \omega^k \odot z_j^k - \mathcal{P}_{\mathcal{M}_j^{\omega^{k}}}(z_j^k) \right), \label{eq: pie_update}
    \\
    \omega^{k+1} &= \omega^k - \gamma \bar{z}_j^{k+1}\left( \omega^k \odot z_j^{k+1}  - \mathcal{P}_{\mathcal{M}_j^{\omega^k}} (z_j^{k+1})\right), \quad j=1, \ldots, N,
    \label{eq: epie_probe_update}
\end{align}
where $\beta, \gamma > 0$ are positive scalars corresponding to step sizes generally chosen to be small~\cite{rodenburg2008ptychography, thibault2009probe} and $\bar{\omega}$ is the the complex conjugate of $\omega$.
The above iterates correspond to performing gradient descent updates on $\Phi$ with respect to $z_j$ and $\omega$ in an alternating manner~\cite{fung2020multigrid}. Importantly, we note that to recover $z$ from the latest $z_j, j=1,\ldots,N$ iterates, the pixels extracted (or illuminated) by $Q_j$ are updated in $z$. Finally, we remark that when the probe $\omega$ is known, then the problem is known as a non-blind ptychographic retrieval problem, and the corresponding algorithm simply uses~\eqref{eq: pie_update} and is known as PIE.

\section{A Hybrid ePIE algorithm}
In this section, we describe the PFT and present our proposed hybrid ePIE algorithm.

\subsection{The Fast Partial Fourier Transform (PFT)}
\label{subsec: pft}
As seen in Section~\ref{sec: ptychography_background}, the key ingredient in the blind (and non-blind) ptychography problem is the discrete fast Fourier transform (FFT), $\mathcal{F}$. Although the discrete FFT is fast and can be applied in $\mathcal{O}(n \log n)$ complexity, it can be the primary bottleneck of reconstruction algorithms such as ePIE when the size of the image is extremely high ($n \gg 1$). 

Recently, an algorithm for approximating the FFT, called the fast partial Fourier transform (PFT)~\cite{park2021fast}, was introduced to speed up the computation of the FFT with applications to time series. The primary motivation for the PFT arises from applications where the the resulting data from the FFT in the frequency domain is sparse i.e., the Fourier coefficients corresponding to high frequencies are predominantly small or equal to zero, and not all Fourier coefficients are necessary for the task at hand~\cite{park2021fast}. 
Indeed, computing the FFT and then cropping the necessary coefficients still requires a cost of $\mathcal{O}(n \log n)$. 
To this end, the core idea of the PFT is to have a fast mechanism that truncates the high frequencies from the FFT in an efficient manner; in particular, the application of the PFT can be done at a cost of $\mathcal{O}(n + \widetilde{m} \log \widetilde{m})$ where $\widetilde{m} \ll n$. 

For ease of presentation, we describe the one-dimensional PFT; however, a two-dimensional version is a straightforward extension and its description can be found in Appendix~\ref{app: pft_2d_description}. To implement the PFT, there are two phases, an offline phase, which is used to pre-compute the polynomial approximation (see Algorithm~\ref{alg: 1d_pft_offline}), and then an online phase which allows one to apply the PFT on the fly (see Algorithm~\ref{alg: 1d_pft_online}). As we are only concerned with ptychographic reconstructions, we are primarily concerned with the online phase of the PFT application and also leave a more thorough derivation of the offline phase in Appendix~\ref{app: pft_2d_description}. 

\subsubsection{PFT Offline Configuration Phase: Approximating Twiddle Factors}
\label{subsec: pft_offline_derivation}
While a thorough derivation and description of the PFT is presented in ~\cite{park2021fast}, we present parts of the derivation in this work for completeness.
Recall that the Discrete Fourier Transform (DFT) is given by:
\begin{equation}
    \hat{z}_{t} = \sum_{k \in [n]}  z_k e^{\left({-2\pi i t k}/{n}\right)}
\end{equation}
where $z \in \CC^{n}$ is an $n$-dimensional complex data vector and $[n] = \{0, 1, \ldots, n-1\}$.
We assume that $n$ is a composite integer so that there exist $p, q > 1$ such that $n = pq$. The Cooley-Tukey algorithm~\cite{cooley1965algorithm} rearranges the above expression as:
\begin{equation}
    \begin{split}
    \hat{z}_t = \sum_{k \in [p]}\sum_{j \in [q]} z_{qk + j} e^{ {-2\pi i t(qk + j)}/{n}  }
    = \sum_{k \in [p]}\sum_{j \in [q]} z_{qk + j} e^{\left(-2\pi itj/{n}\right)} \cdot e^{\left({-2\pi itk}/{p}\right)}.
    \end{split}
\end{equation}
Further modification of the above expression yields
\begin{equation}
    \begin{split}
        \hat{z}_{t} & = \sum_{k \in [n]} z_{k}e^{-2\pi it(k - q/2)/n} \cdot e^{-\pi it/p} \\
        & = \sum_{k \in [p]}\sum_{j \in [q]} z_{qk + j}e^{-2\pi it(j - q/2)/n} \cdot e^{-2\pi itk/p} \cdot e^{-\pi it/p}.
    \end{split}
    \label{eq: modified_summation}
\end{equation}
Here, $[q]$ corresponds to the set $\{0,1,\ldots,q-1\}$.
The key idea behind the PFT is to use a polynomial to approximate the exponential $e^{\pi i x}$. Afterwards, one can re-scale the polynomials and use exponent laws to get an approximation of each of the \emph{twiddle factors} in the collection $\left\{ e^{-2 \pi i t (j - q/2) / n} \right\}_{j=0}^{q - 1}$.

To choose the approximating polynomial of the exponential function $e^{\pi ix}$, we consider the choice of hyper-parameter $p$ given $n$ and desired output size $\widetilde{m}$. Let $\|\cdot\|_{R}$ be the uniform norm (supremum norm) restricted to a set $R \subseteq \RR$, that is, $\|f\|_{R} = \sup\left\{|f(x)|: x \in \mathbb{R}\right\}$. Then, given non-negative integer $\alpha$, and non-zero real number $\xi$, we define the polynomial $\cP_{\alpha, \xi}$ as the best approximation to $e^{\pi ix}$ out of the space $P_{\alpha}$ of polynomials on $\mathbb{R}$ of degree at most $\alpha$ under the restriction $|x| \leq |\xi|$ as:
\begin{equation}\cP_{\alpha, \xi} = \argmin_{P \in P_{\alpha}} \left\| P(x) - e^{\pi ix}\right\|_{|x| \leq |\xi|}.
\end{equation}
We remark that such polynomials exist and are unique and there are minimax algorithms that can be used to compute them~\cite{fraser1965survey}.
The goal is to choose $\xi \in \mathbb{R}$ such that $|\widetilde{m}| \leq |\xi|$. Given a tolerance $\varepsilon > 0$ and a positive integer $r \geq 1$, we define $\xi(\varepsilon, r)$ to be the scope about the origin such that the exponential function $e^{\pi ix}$ can be approximated by a polynomial of degree less than $r$ with approximation bound $\varepsilon$:
\begin{equation}
    \xi(\varepsilon, r) := \sup\left\{\xi \geq 0 : \left\|\cP_{r - 1, \xi}(x) - e^{\pi ix}\right\|_{|x| \leq \xi} \leq \varepsilon\right\},
\end{equation}
where we express the corresponding polynomial as $\cP_{r - 1, \xi}(x) = \sum_{j\in[r]} w_{\varepsilon,r,j} x^j$. 
Using a minimax approximation algorithm~\cite{fraser1965survey}, precompute $\xi(\varepsilon, r)$ and $\{w_{\varepsilon, r, j}\}_{j \in [r]}$ for several tolerance $\varepsilon$'s (e.g. $10^{-1}, 10^{-2}, \dots$) and positive integer $r$'s (the authors of ~\cite{park2021fast} choose values in the range $1 \leq r \leq 25$). 
When $n, \widetilde{m}, p,$ and $\varepsilon$ are given, we choose the minimum $r$ satisfying $\xi(\varepsilon, r) \geq \widetilde{m}/p$. Following the preceding argument, one can then show that the re-scaled polynomial $\cP_{r-1, \xi(\varepsilon, r)}(-2t(j - q/2)/n)$ approximates $e^{-2 \pi it (j - q/2) / n}$ on $|t| \leq \left|\frac{n}{2(j - q/2)} \cdot \frac{\widetilde{m}}{p}\right|$ for each $j \in [q]$~\cite{park2021fast}. 
Noting that 
$
\left|\frac{n}{2(j - q/2)} \cdot \frac{\widetilde{m}}{p}\right| = \left|\frac{q}{2j - q} \cdot \widetilde{m}\right| \geq \widetilde{m}
\quad \text{for all} \quad j \in [q],
$
we have a polynomial approximation on $|t| \leq \widetilde{m}$ for each twiddle factor in the collection $\left\{ e^{-2 \pi i t (j - q/2) / n} \right\}_{j=0}^{q - 1}$, namely the the re-scaled polynomials $\left\{\cP_{r-1, \xi(\varepsilon, r)}(-2t(j - q/2)/n)\right\}_{j = 0}^{q - 1}$. 
Algorithm~\ref{alg: 1d_pft_offline} shows how to build the polynomial approximation\footnote{The algorithm written here centers the PFT at 0 - this is all we need for our purposes. However, with a minor modification, one may also compute the PFT centered at another coordinate~\cite{park2021fast}}. The polynomial coefficients $w_{\epsilon, r-1, j}$ are precomputed, and we obtain them from the code database of~\cite{park2021fast}.
\begin{algorithm}[t]
    \small{
    \caption{: Configuration (Offline) Phase of 1D PFT}
    \label{alg: 1d_pft_offline}
    \hspace*{\algorithmicindent} \textbf{Input}: Input size $n \in \mathbb{N}$, crop size $\widetilde{m} \in \mathbb{N}$, divisor $p \in \mathbb{N}$, and tolerance $\varepsilon$ \\
    \hspace*{\algorithmicindent} \textbf{Output}: Matrices $B \in \mathbb{C}^{q \times r}, W \in \mathbb{C}^{(2\widetilde{m} + 1) \times r}$, divisor $p$, and number of rows and columns, $q$ and $r$
    \begin{algorithmic} [1]
        \State $q = n / p$
        \State $r = \min\{r \in \NN : \xi(\varepsilon, r) \geq \widetilde{m}/p\}$
        \Comment{degree of polynomial $\cP$ approximating $e^{\pi ix}$ within tolerance $\varepsilon$}
        \For{$(l, j) \in [q] \times [r]$} 
            \State $x = (1 - 2l/q)$
            \State $B[l, j] = w_{\varepsilon, r-1, j} \cdot x^{j}$ \Comment{Using precomputed $w_{\varepsilon, r-1, j}$}
        \EndFor

        \For{$k = -\widetilde{m}, -\widetilde{m}+1, \ldots, \widetilde{m}-1, \widetilde{m}$}
            \State $W[k, j] = (k/p)^{j} \cdot e^{-\pi i k/p}$
            \Comment{Precompute remaining terms}
        \EndFor
    \end{algorithmic}
    }
\end{algorithm}
\subsubsection{Online PFT Computation}
Substituting the approximating polynomial described in Section~\ref{subsec: pft_offline_derivation} for the twiddle factors in~\eqref{eq: modified_summation} and performing some algebraic manipulations, we can represent the summation as a matrix-matrix multiplication $C = Z \times B$ where $Z \in \CC^{p \times q}$ is the reshaped data vector $z$ and $B \in \CC^{q \times r}$ is defined in line 5 of Algorithm~\ref{alg: 1d_pft_offline} followed by a series of FFT computations applied to the columns of the resulting matrix $C$.
Worth noting, the primary cost in Algorithm~\ref{alg: 1d_pft_online} arises from a) the matrix-matrix multiplication between a matrix $Z \in \mathbb{C}^{p \times q}$ and $B \in \mathbb{C}^{q \times r}$ in line 2 and b) the application of the FFT on a vector of size $\mathbb{C}^{p}$, which is shown to have complexity $\mathcal{O}(n + \widetilde{m}\log \widetilde{m})$~\cite[Theorem 3]{park2021fast}. For thorough details, we refer the reader to~\cite[Section 3.4]{park2021fast}.

A hyper-parameter that needs to be chosen for the PFT computation is the integer $p$ such that $n = pq$. In our experiments, we are interested in the setting where $\widetilde{m} \ll n$ as we wish to economically capture low frequency features when deploying the PFT. In this case, the primary cost in Algorithm~\ref{alg: 1d_pft_online} occurs in the matrix-matrix multiplication in line 2, which has complexity {$\mathcal{O}(nr)$}. Thus, it is recommended that  $\widetilde{m}/p$ be small so that the number $r$ of approximating polynomial terms decreases~\cite[Section 4.3]{park2021fast}. As previously mentioned, a 2D implementation of the PFT is a straightforward extension of the 1D PFT, where given an image of size $n_1 \times n_2$, we have to decide on $\widetilde{m}_1, \widetilde{m}_2, p_1, p_2$. In our experiments, we find that choosing $p_1 = p_2 = \widetilde{m}_1 = \widetilde{m}_2 = 64$ is adequate for improved performance. In this case, the resulting cropped image, has size $128 \times 128$ since the cropping is performed from $[-\widetilde{m}_1, \widetilde{m}_1] \times [-\widetilde{m}_2, \widetilde{m}_2]$ from the center; we note that an FFT on an image of size $128 \times 128$ is trivial to compute. 

\begin{algorithm}[t]
    \small{
    \caption{: Computation (Online) Phase of 1D PFT}
    \label{alg: 1d_pft_online}
    \hspace*{\algorithmicindent} \textbf{Input}: Vector $z \in \mathbb{C}^n$, crop size $\widetilde{m}$, and configuration results $B \in \mathbb{C}^{q \times r}, W \in \mathbb{C}^{2\widetilde{m}+1 \times r}, p, q, r$ \\
    \hspace*{\algorithmicindent} \textbf{Output}: Vector $\hat{z}_{\text{PFT}}$ of estimated Fourier coefficients of $z$ in range $[-\widetilde{m}, \widetilde{m}]$
    \begin{algorithmic} [1]
        \State{\begin{tabular}{p{0.500\textwidth}r}
         \hspace*{-8pt} 
        $Z = z\text{.reshape}(p, q)$ 
        & 
        \Comment{reshape $z$ into $p \times q$ matrix}
        \end{tabular}
         }
        \State{
        \begin{tabular}{p{0.492\textwidth}r}
         \hspace*{-13pt} 
        $C = Z \times B$ 
        & 
        \Comment{matrix multiply $Z$ by $B$}
        \end{tabular}
        }
        \For{$j \in [r]$}
            \State{\begin{tabular}{p{0.465\textwidth}r}
         \hspace*{-13pt} 
         $\hat{C}[:, j] = $ FFT($C[:, j]$)
         & 
         \Comment{apply FFT to each column of $C$}
         \end{tabular}
         }
        \EndFor
       \For{$k = -\widetilde{m}, -\widetilde{m}+1, \ldots, \widetilde{m}-1, \widetilde{m}$}
            \State{\begin{tabular}{p{0.465\textwidth}r}
         \hspace*{-13pt}  
         $\hat{z}_{\text{PFT}}[k] = \sum_{j=0}^{r-1} \left(\hat{C}[k\%p,j] \cdot W[k,j]\right)$ 
         &
         \Comment{sum values of $\hat{C}$ (modulus  $p$) times $W$}
         \end{tabular}
         }
        \EndFor
    \end{algorithmic}
    }
\end{algorithm}

\subsection{Hybrid ePIE Algorithm}
We propose a hybrid ePIE algorithm that consists of using a PFT-based ePIE as a ``warm up'' followed by the standard FFT-based ePIE. 
The idea is to use the PFT-based ePIE to capture large features arising from the low frequencies in a cheap manner. 
Indeed, while one may use the PFT in a plug-and-play manner, we empirically observe that it is better to initially capture the large features early on, and let the standard FFT-based ePIE algorithm capture the fine details in later iterations. 

To run the PFT-based ePIE, we convert the observed data $d_j, \; j=1,\ldots,N$ to match the output dimension of $\mathcal{F}_{\text{PFT}}$. Considering the reshaped version of $d_j$ as a matrix, this can be easily done by cropping the centers of $d_j \in \mathbb{R}^{m_1 \times m_2}$ to generate $d_j^{\text{crop}} \in \mathbb{R}^{\widetilde{m}_1 \times \widetilde{m}_2}$.

\begin{algorithm}[H]
    \small{
    \caption{: Hybrid ePIE Algorithm}
    \label{alg: hybrid_ePIE}
    \hspace*{\algorithmicindent} 
    \textbf{Input}: 
    Initial Guess $z \in \mathbb{C}^{n_1 n_2 \times 1}$, 
    crop size parameters $\widetilde{m}_1$, $\widetilde{m}_2$, 
    observed data $d_j \in \mathbb{R}^{m_1 \times m_2}$ for $ j=1,\ldots,N$, stopping tolerances $\epsilon_{PFT}, \epsilon$, maximum number of iterations $n_{\text{maxiters}, PFT}, n_{\text{maxiters}}$
    \\
    \hspace*{\algorithmicindent} 
    \vspace{1mm}
    \textbf{Output}: Solution $z_{\text{opt}}$
    \begin{algorithmic} [1]
        \State{\begin{tabular}{p{0.500\textwidth}r}
         \hspace*{-8pt} 
        Construct cropped data $d_j^{\text{crop}} \in \mathbb{R}^{\widetilde{m}_1\widetilde{m}_2 \times 1}$ from $d_j$
        & 
        \Comment{crop centers of $d_j$, $j=1\ldots,N$}
        \end{tabular}
         }
        \vspace{1mm}
        \State{
        run PFT-based ePIE iterates~\eqref{eq: pie_update} - \eqref{eq: epie_probe_update} using $d_j^{\text{crop}}, \; j=1,\ldots,N$ until convergence to obtain $z_{\text{PFT}}$
        }
        \vspace{1mm}
        \State{
        run FFT-based ePIE iterates~\eqref{eq: pie_update} - \eqref{eq: epie_probe_update} using $d_j, \; j=1,\ldots,N$ and $z_{\text{PFT}}$ as an initial condition until convergence to obtain $z_{\text{opt}}$.
        }
    \end{algorithmic}
    }
\end{algorithm}

\section{Numerical Experiments}
\begin{figure}[t]
    \centering
    \begin{tabular}{cc}
        True Magnitude & True Phase
        \\
        \includegraphics[width=0.32\textwidth]{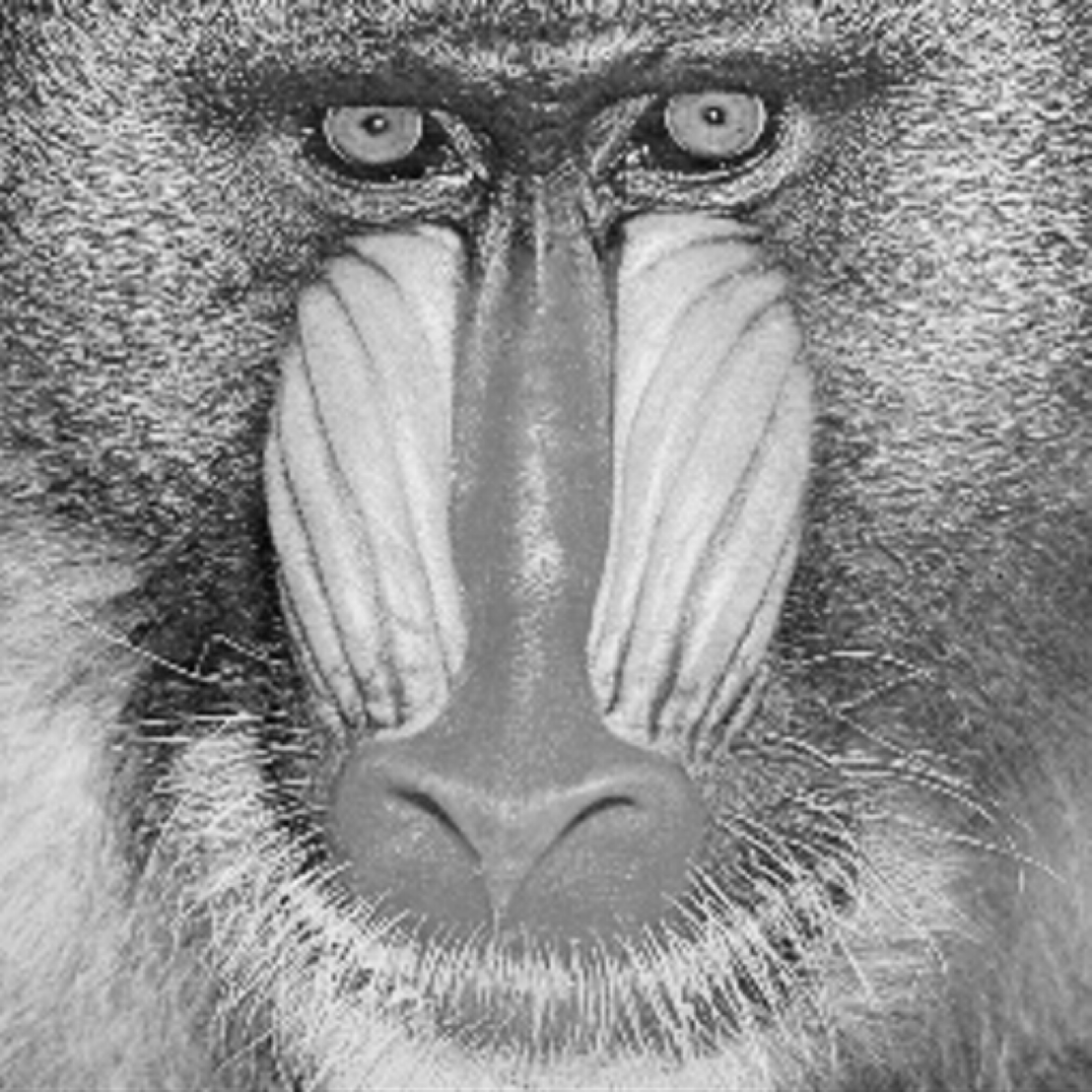}
        &
        \includegraphics[width=0.32\textwidth]{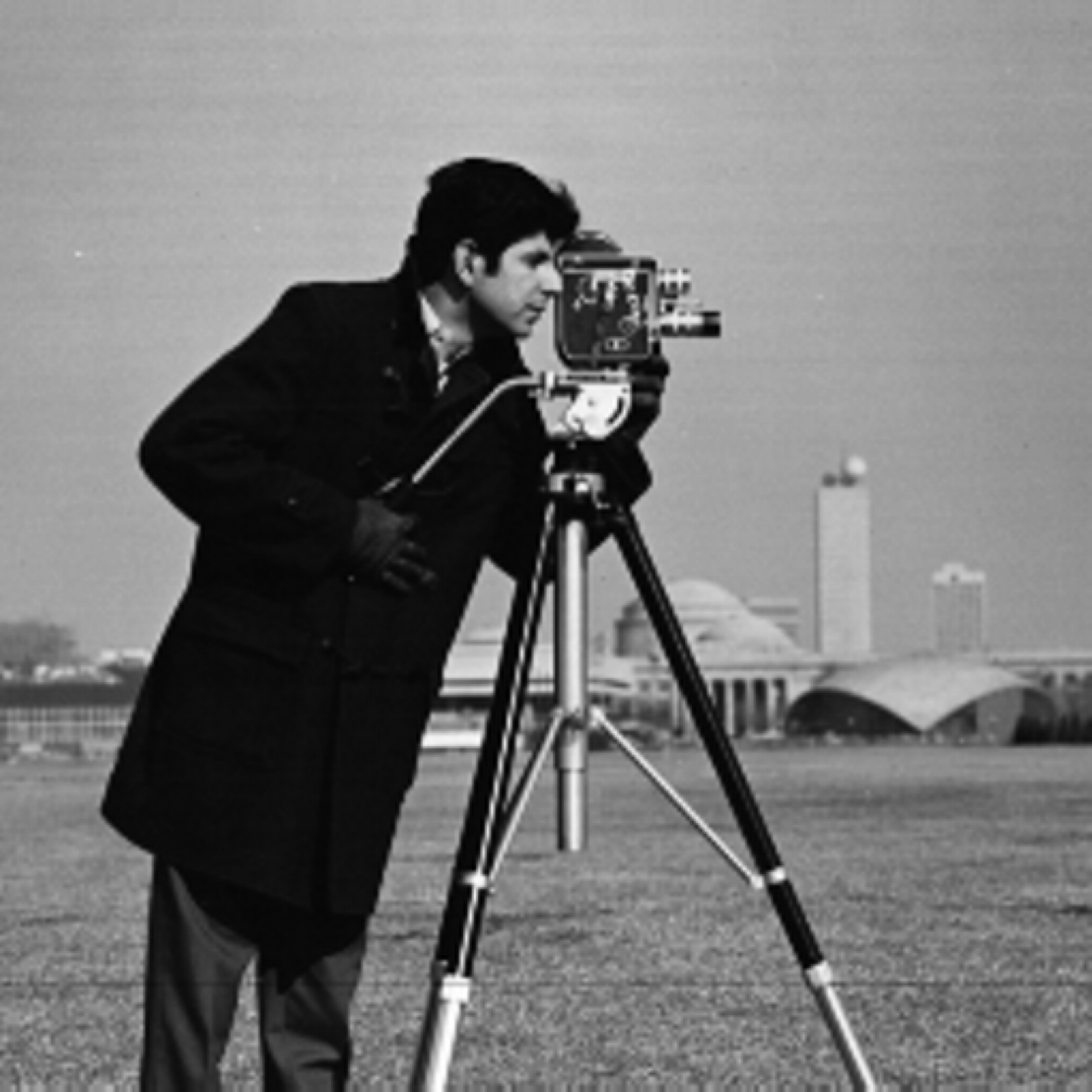}
    \end{tabular}
    \caption{\small{The ground truth used to simulate data in numerical experiments. The baboon image is used as the magnitude and the cameraman image is used as the phase of the object of interest.}}
    \label{fig: true_object}
\end{figure}

We consider the non-blind and the blind ptychography problem. For each problem type, we investigate the quality of reconstructions on a small experiment with image size $512\times512$ and study the quality of local minima. 
Indeed, the PFT is only worth employing when the runtime over the FFT is reduced - this is not as obvious in the small experiments. To this end, we demonstrate the computational benefits of using the PFT, and in particular the proposed hybrid PIE, on a large experiment with image size $16384 \times 16384$.
For all experiments, we assume $m = n$, i.e. square images. We work with two images, the baboon and cameraman. We set the baboon image to be the magnitude and the cameraman to be the phase of the ground truth, as shown in Figure~\ref{fig: true_object}. Here, the magnitude is chosen to have range $[0,1]$ and the phase is chosen to have range $[0, \pi/2]$. 
As explained in Section~\ref{subsec: pft}, one wants $p_1$ and $p_2$ such that $\tilde{m}_1/p_1$ and $\widetilde{m}_2/p_2$ is small; in our case, we find that choosing $\tilde{m} = p = 64$ is good enough for all of our experiments. More details on the choice of hyperparameter $p$ can be found in~\cite[Section 4.3]{park2021fast}. Finally, the stopping criteria considered in all algorithms is $\| z^{k+1} - z^k \| / \|z^k\|$. Following~\cite{park2021fast}, we choose a tolerance of $10^{-7}$ in the computation of the approximating polynomial in Algorithm~\ref{alg: 2d_pft_offline}. Our numerical results are also based on a PyTorch implementation of the PFT that was translated from that presented in~\cite{park2021fast}. The small-scale experiments (Sections~\ref{subsubsec: nonblind_distribution_relerrs} and~\ref{subsubsec: blind_distribution_relerrs}) are run on a machine equipped with Intel Core i5-9400 CPU, 2.90GHz and 16GB of RAM while large-scale experiments (Sections~\ref{subsubsec: nonblind_large_scale} and~\ref{subsubsec: blind_large_scale}) are run on a Lambda Vector machine with a AMD Threadripper Pro 3955WX processor containing 16 cores, 3.90 GHz, 64 MB cache, PCIe 4.0.

\begin{figure}[t]
    \setlength\tabcolsep{1 pt}
    \centering
    \begin{tabular}{ccccccccc}
        $Q_1$ & $Q_2$ & $Q_3$ & 
        $Q_4$ & $Q_5$ & $Q_6$ &
        $Q_7$ & $Q_8$ & $Q_9$ 
        \\
        \includegraphics[width=0.1\textwidth]{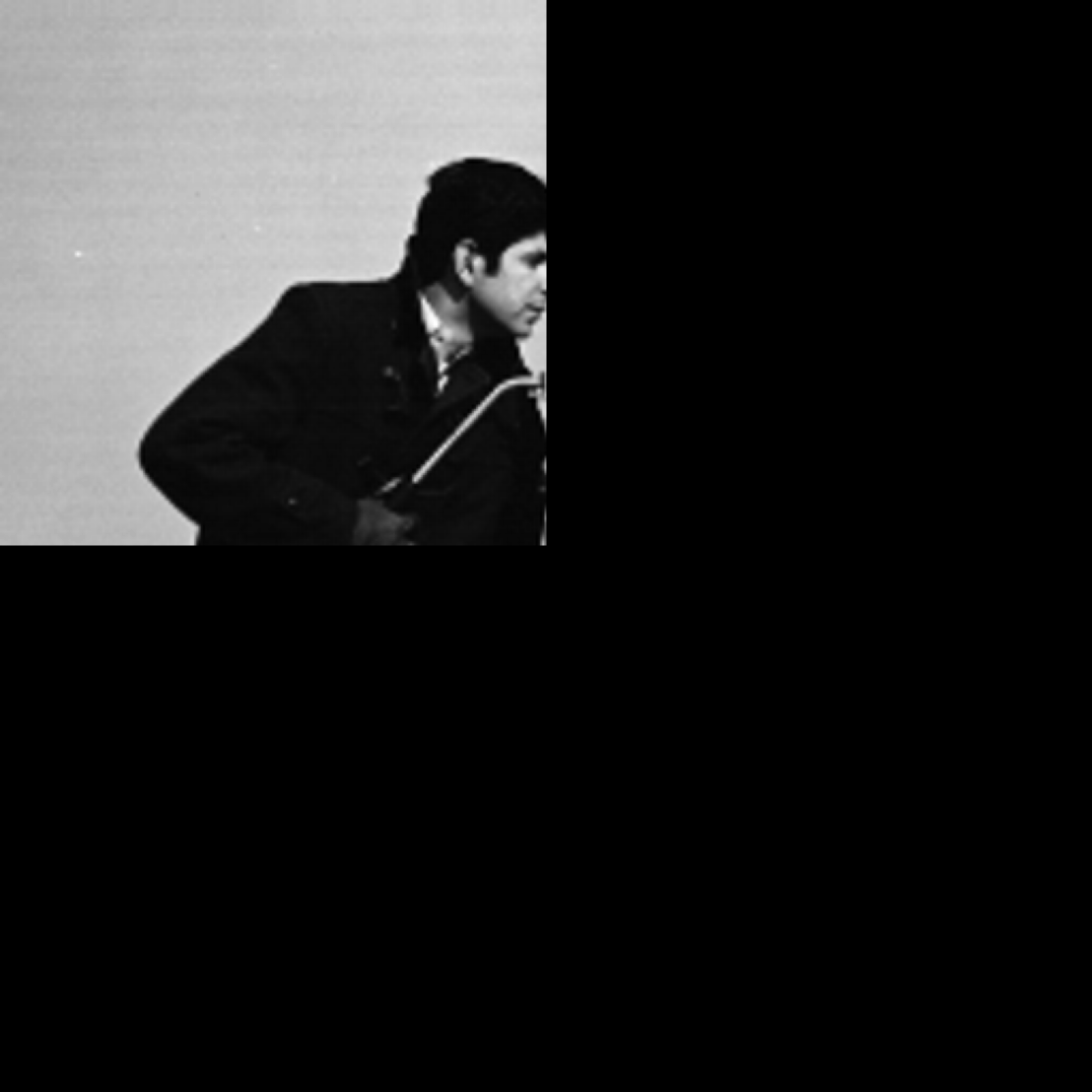}
        &
        \includegraphics[width=0.1\textwidth]{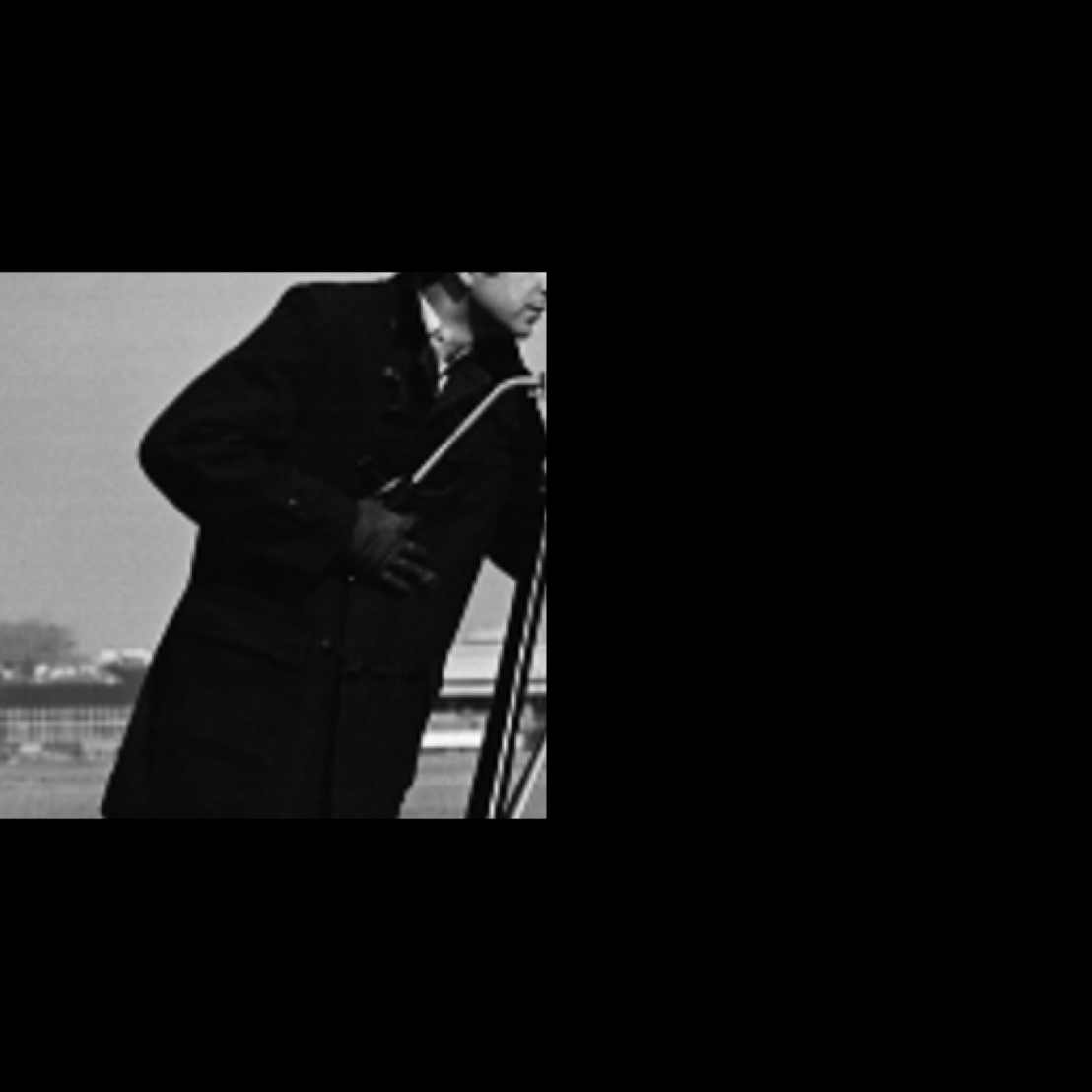}
        &
        \includegraphics[width=0.1\textwidth]{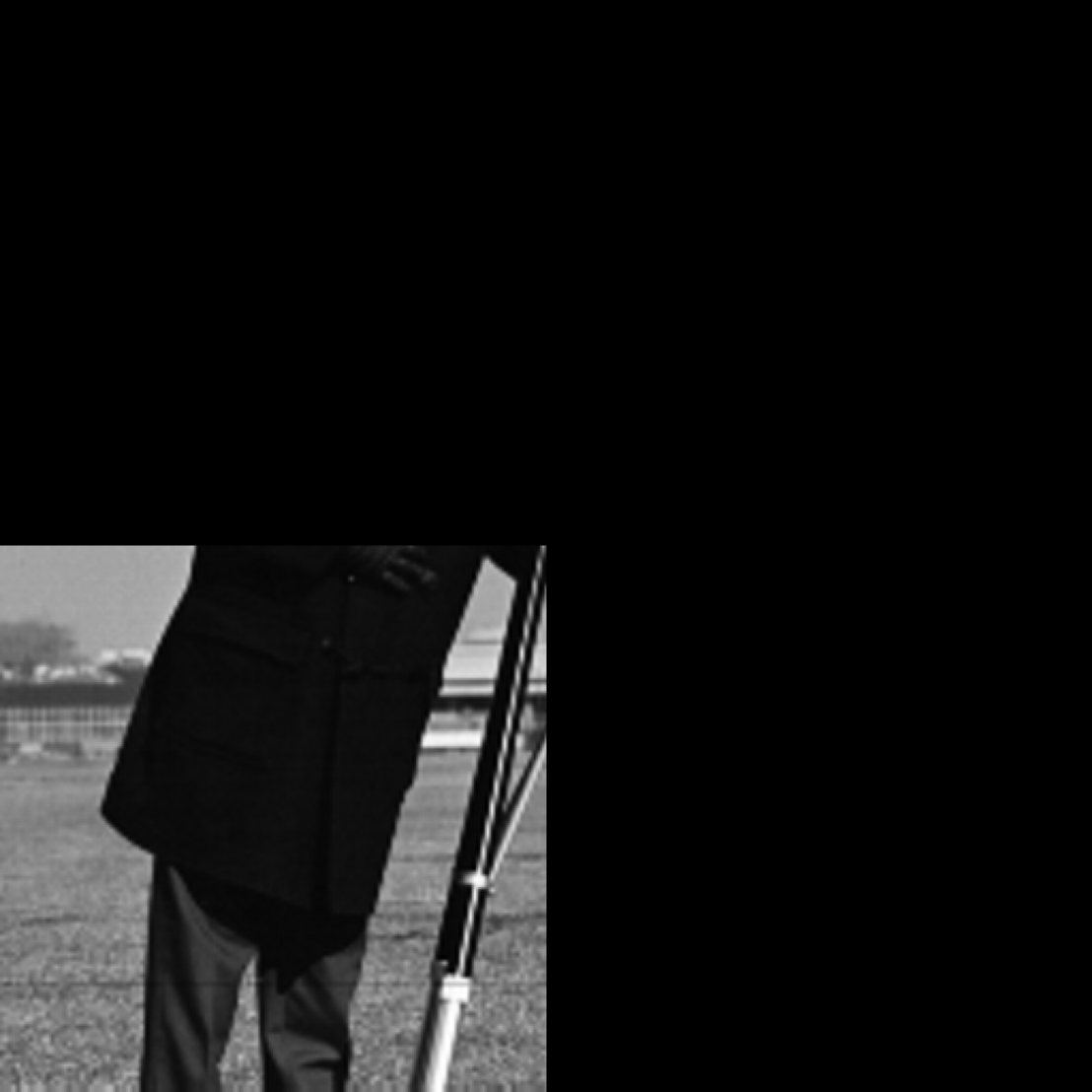} 
        &
        \includegraphics[width=0.1\textwidth]{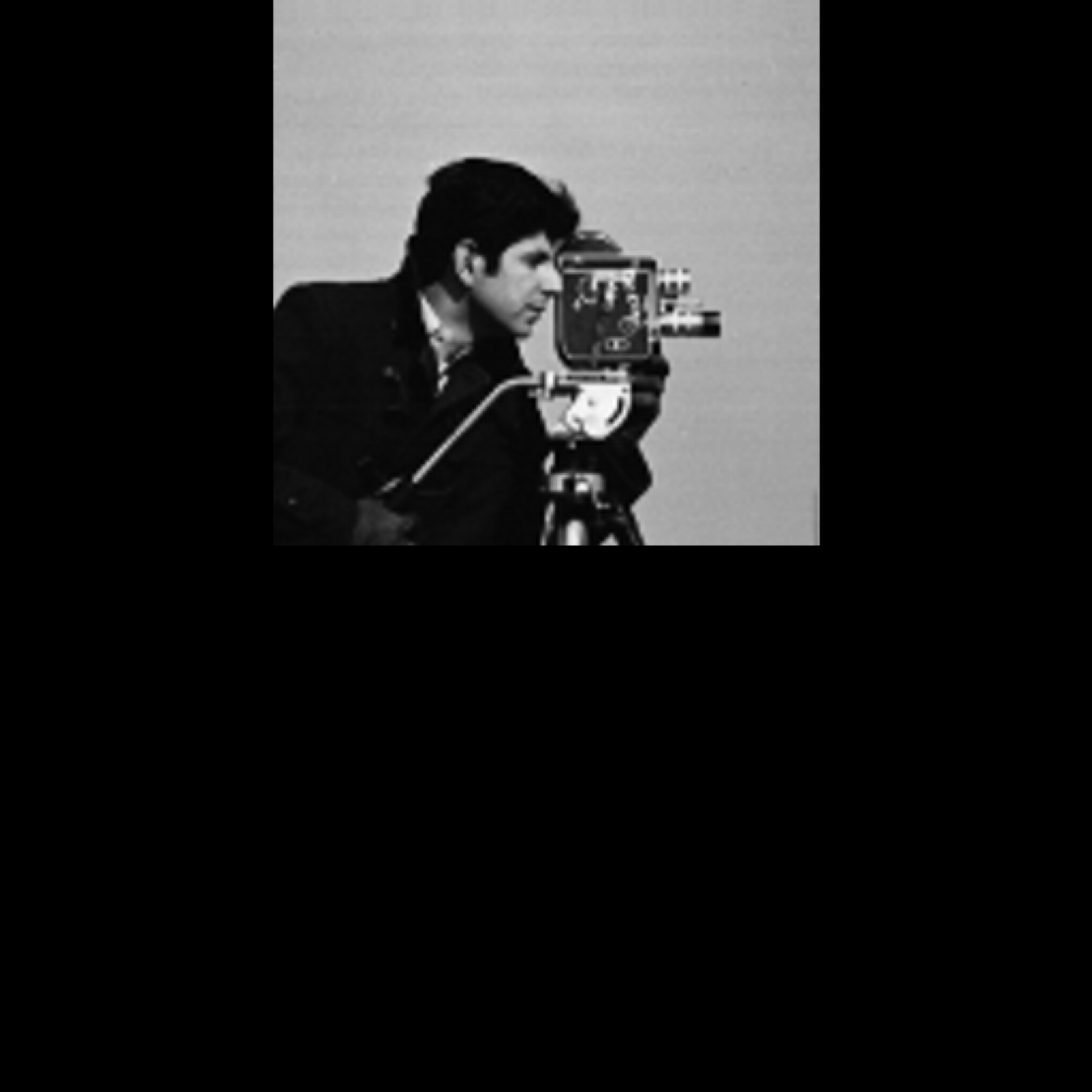}
        &
        \includegraphics[width=0.1\textwidth]{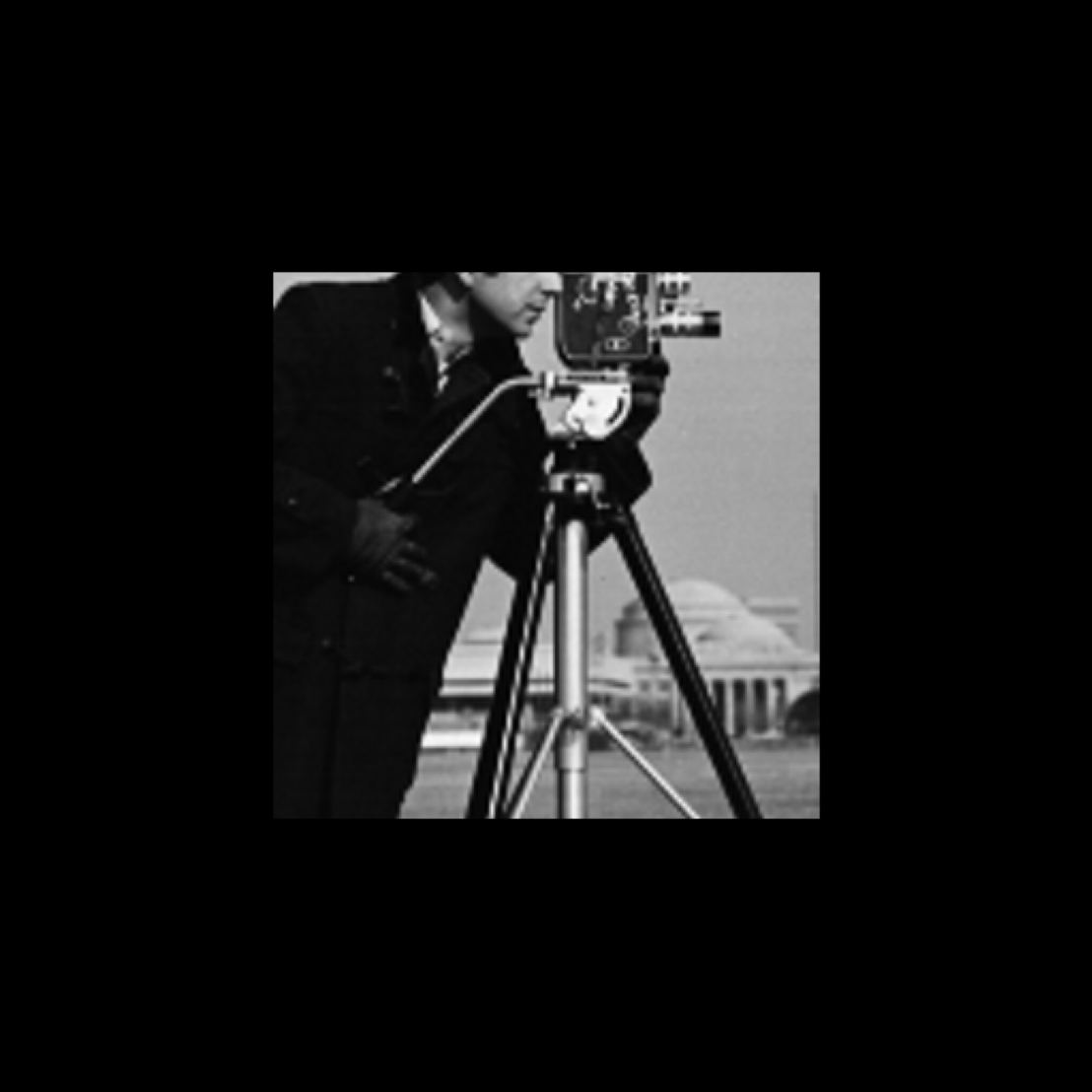}
        &
        \includegraphics[width=0.1\textwidth]{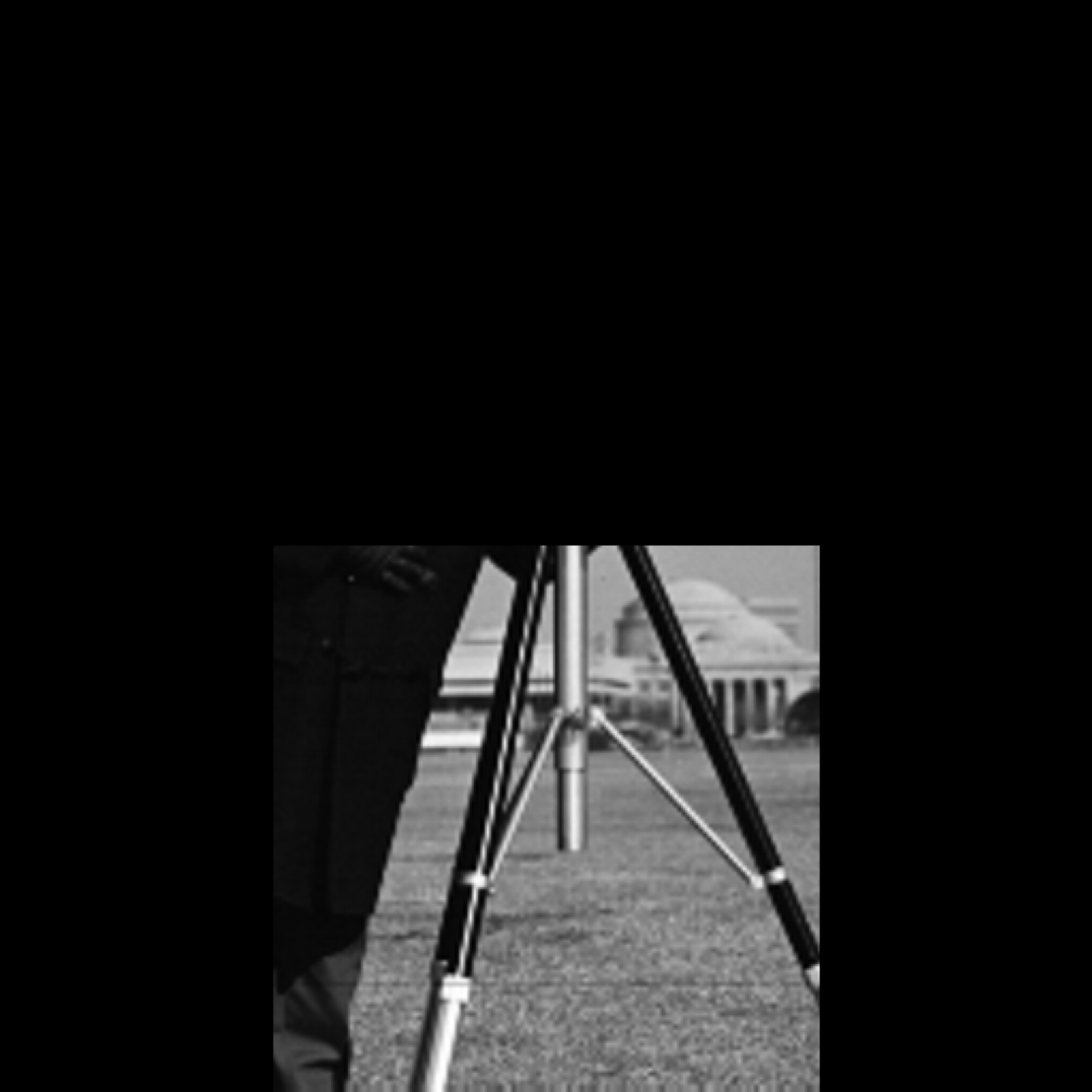} 
        &
        \includegraphics[width=0.1\textwidth]{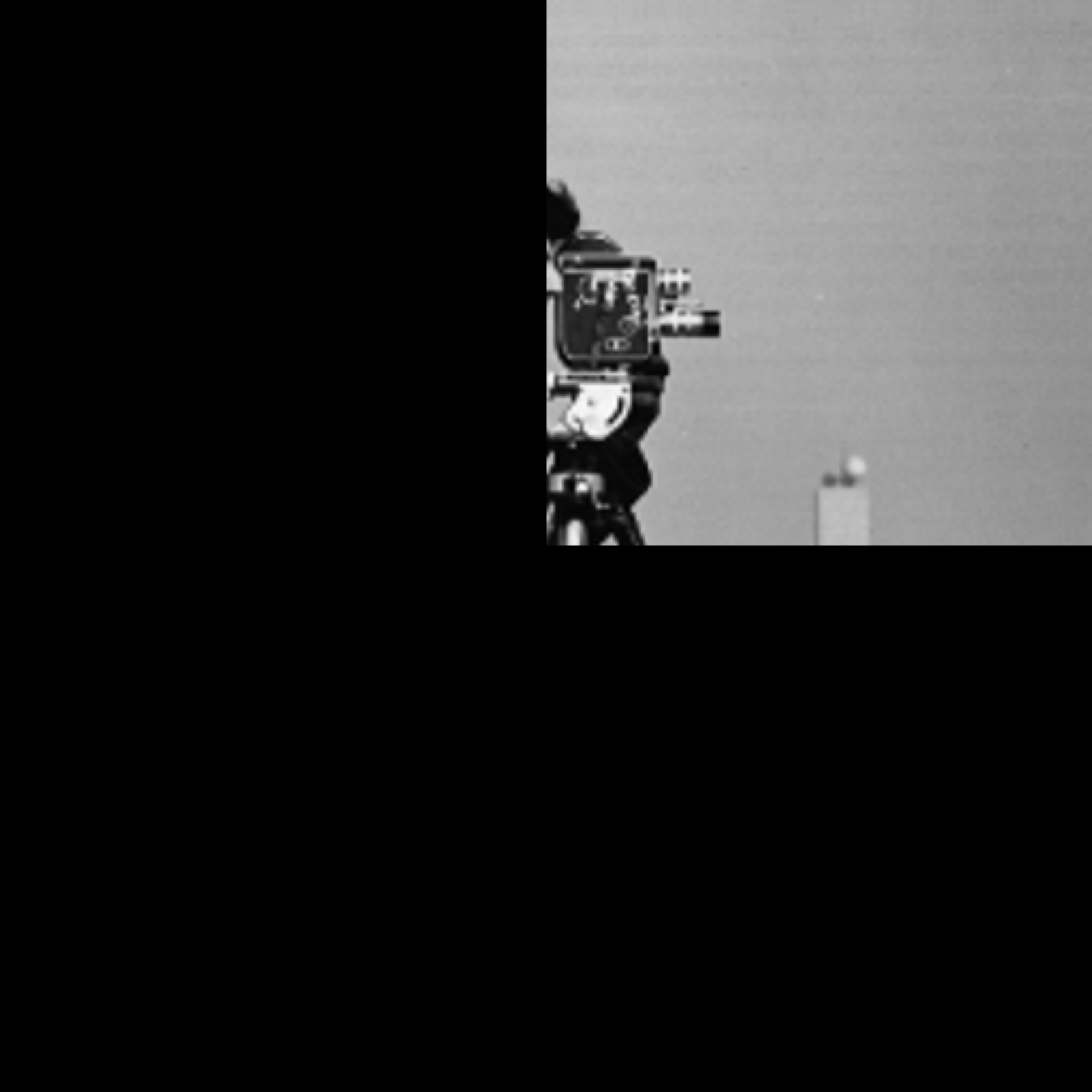}
        &
        \includegraphics[width=0.1\textwidth]{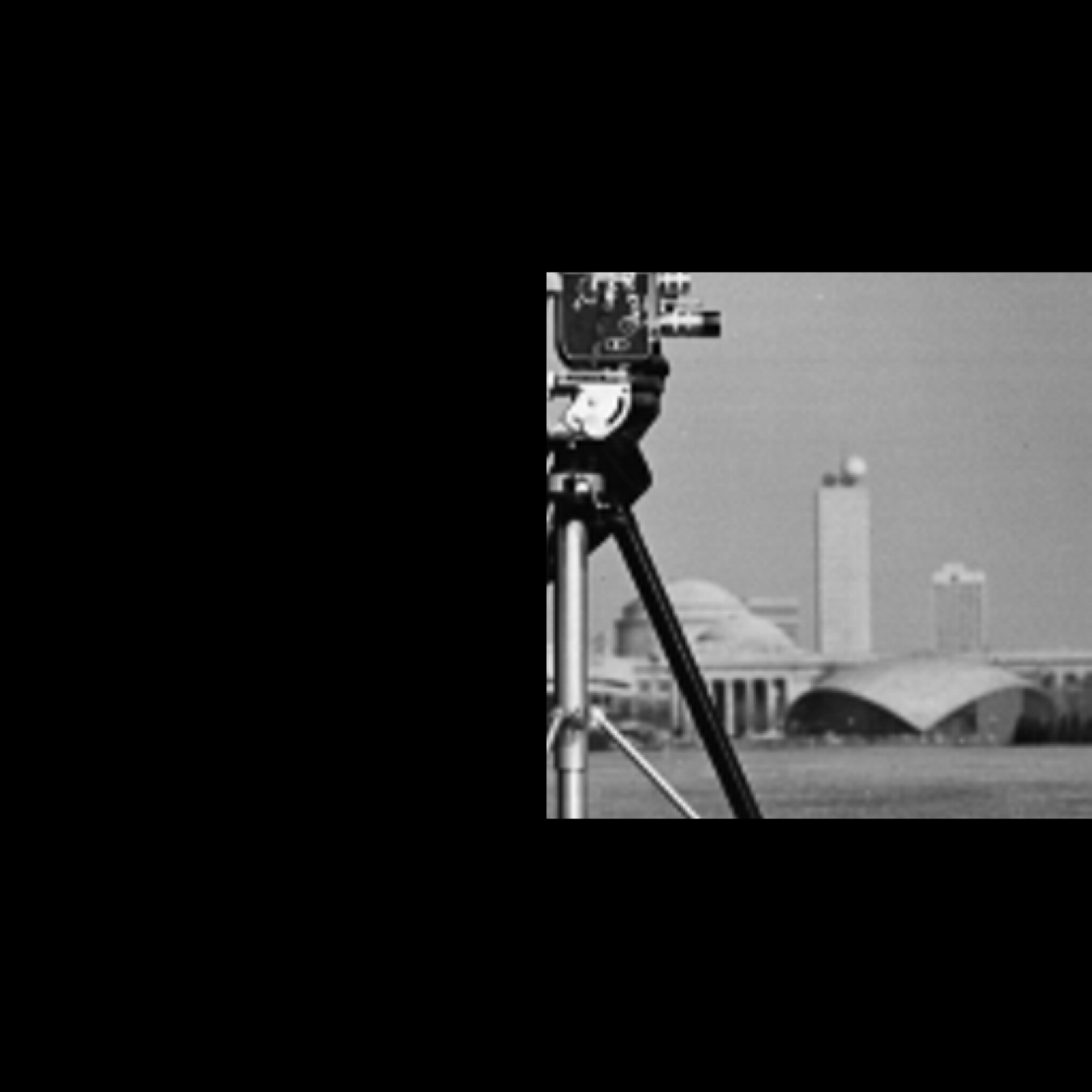}
        &
        \includegraphics[width=0.1\textwidth]{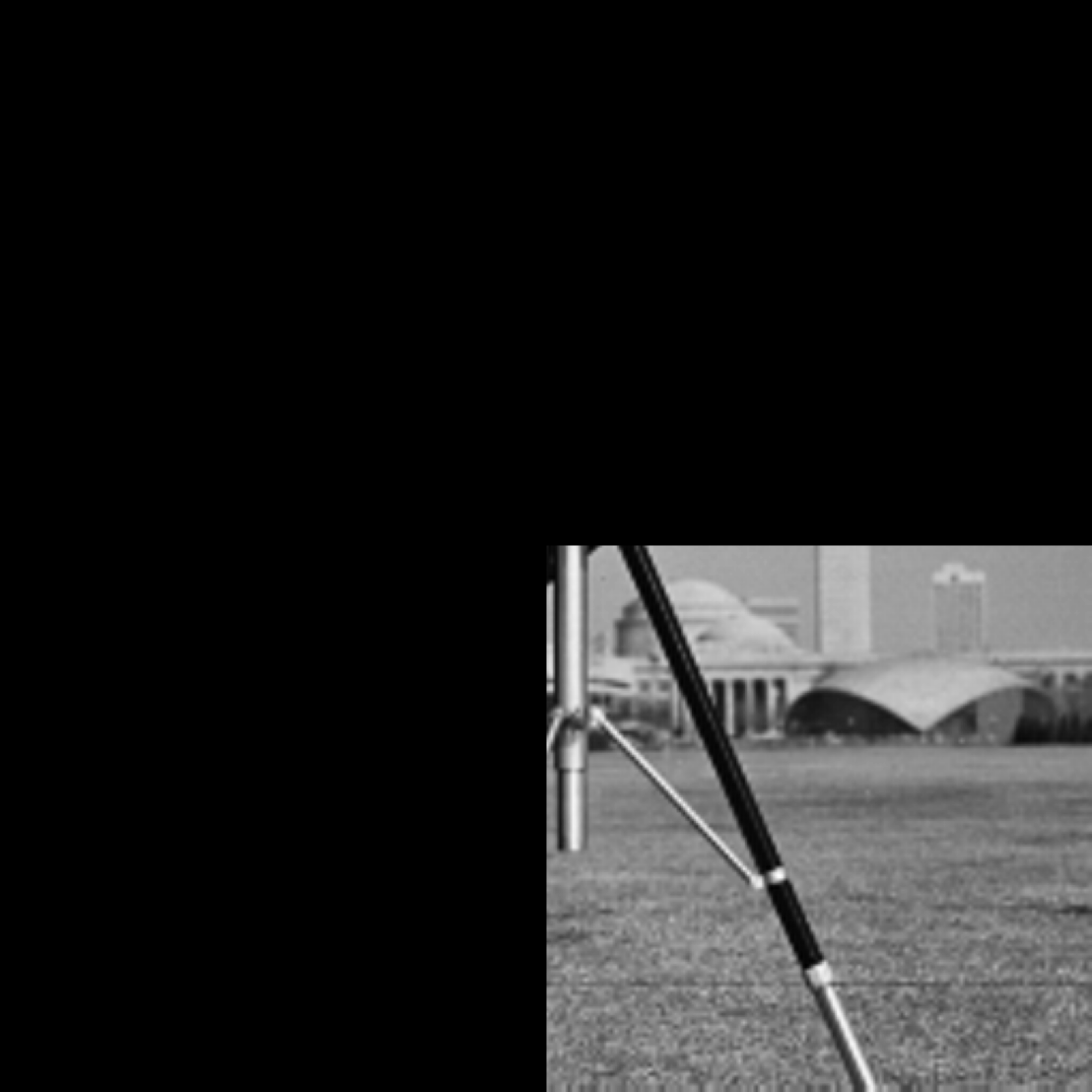}
    \end{tabular}
    \caption{\small{Scanning positions generated by the illumination windows described in Section \ref{subsubsec: nonblind_distribution_relerrs}}. For the small problem where $n_1 = n_2 = 512$ (Section~\ref{subsubsec: nonblind_experimental_setup}), each illumination window $Q_i$ has size $256 \times 256$ and shifts $128$ pixels at a time. Similarly, for the large-scale problem where $n_1 = n_2 = 16384$ (Section~\ref{subsubsec: nonblind_large_scale}, each illumination window $Q_i$ has size $8192 \times 8192$ and shifts $4096$ pixels at a time. In both setups, there is a $50\%$ overlap between consecutive probes.}
    \label{fig: binary_probes}
\end{figure}

\subsection{Non-Blind Ptychography} 
\label{subsec: non_blind_experiment}

\subsubsection{Experimental Setup} 
\label{subsubsec: nonblind_experimental_setup}
As explained in Section~\ref{sec: ptychography_background}, the non-blind ptychography experimental setup assumes $\omega$ is known. In this case, we use probes that resembles an identity operator, that is, $\omega \odot Q_k z = Q_k z, k=1,\ldots,N$. Moreover, each probe illuminates $\frac{n}{2} \times \frac{n}{2}$ pixels. That is, the illumination $Q_{k} \in \RR^{n \times n}$ is a matrix that satisfies
\begin{equation} 
[Q_k]_{i, j} = \begin{cases} 1 & \text{if pixel } (i, j) \text{ is illuminated} \\ 0 &  \text{otherwise }\end{cases}.
\end{equation}
The illumination window is shifted $\frac{n}{4}$ pixels at a time starting from the top-left corner downward until it reaches the bottom-left corner. The probe is then shifted upwards and then to the right until it reaches the top edge again. This scanning procedure is continued until we have covered the entire image, and results in a total of $9$ probes with $50\%$ overlap between adjacent probes (see Figure~\ref{fig: binary_probes} for an illustration).

\begin{figure}[t]
    \centering
    \begin{tabular}{ccc}
        a) $z$ & b) Magnitude of $z$ & c) Phase of $z$
        \\
        \includegraphics[width=0.3\textwidth]{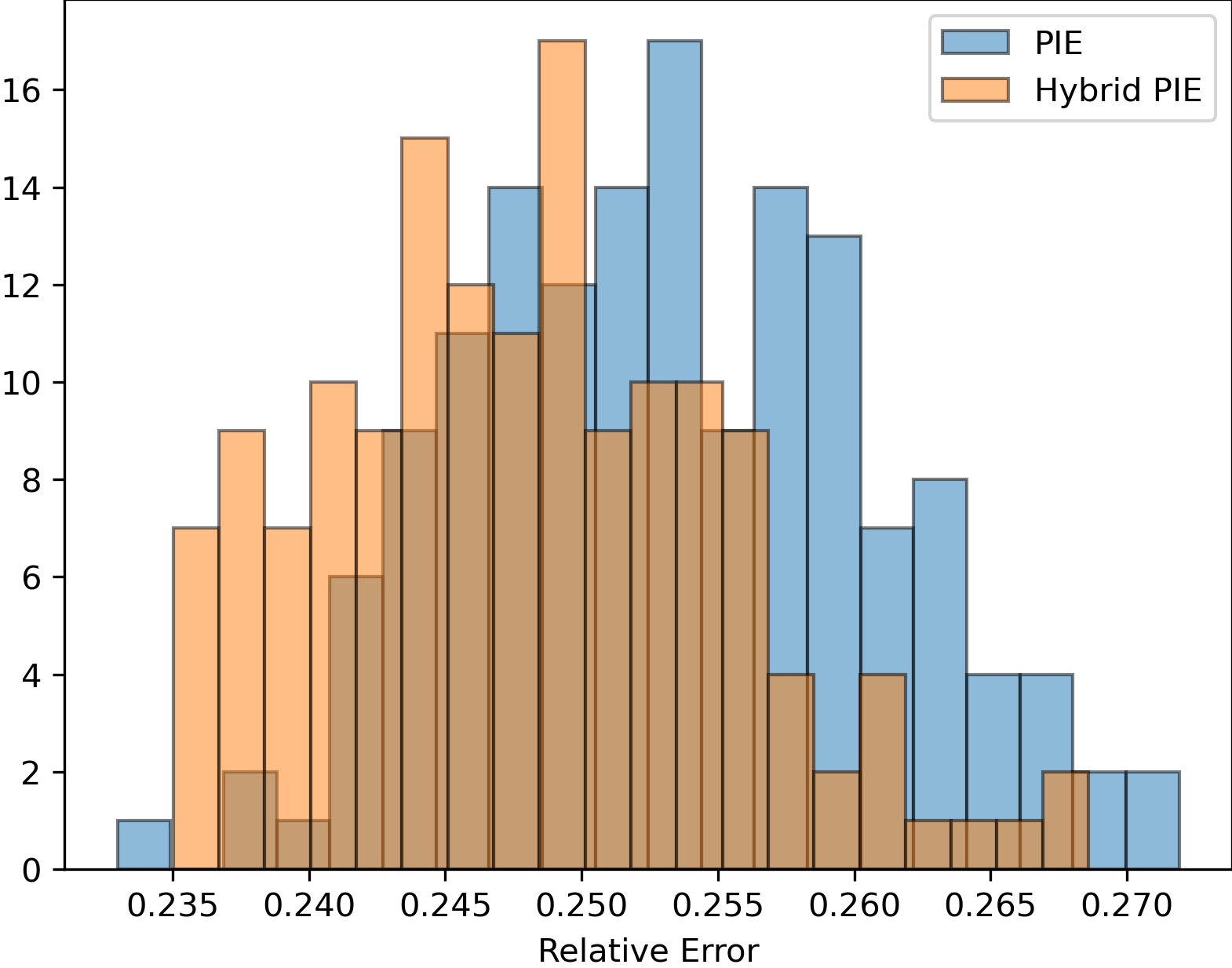}
        &
        \includegraphics[width=0.3\textwidth]{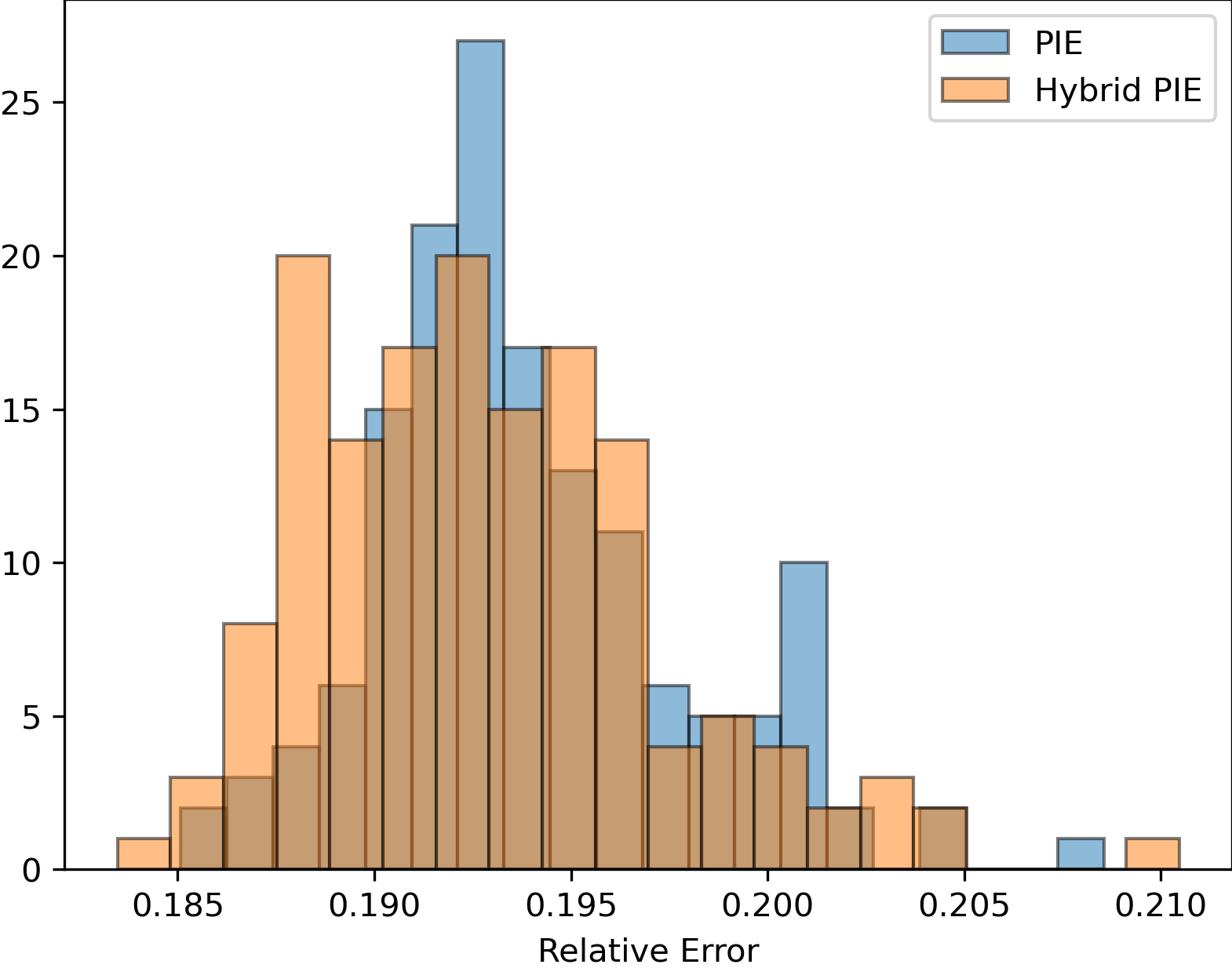}
        &
        \includegraphics[width=0.3\textwidth]{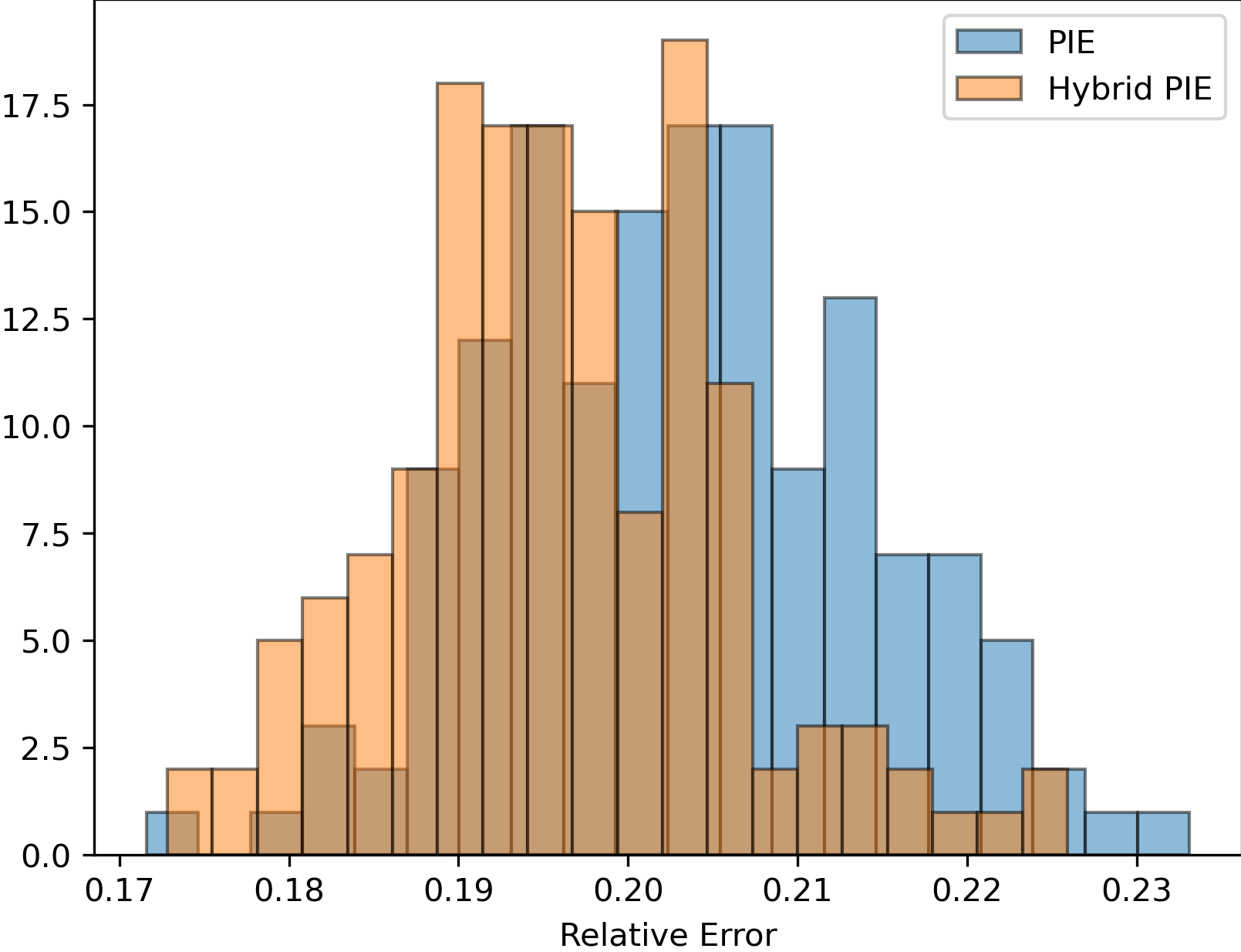}
    \end{tabular}
    \caption{\small{Histogram of the final reconstruction relative errors. The blue histogram shows the relative error frequency for PIE and the orange histogram shows the relative error frequency for hybrid PIE. a) shows the relative error of the reconstructed object, b) shows relative errors of only the magnitude of the object, and c) shows the relative errors of the phase of the object.}}
    \label{fig: pie_rel_err}
\end{figure}

\begin{figure}[t]
    \small
    \centering
    \begin{tabular}{cccc}
        Magnitude SSIM & Phase SSIM & Magnitude PSNR & Phase PSNR
        \\
        \includegraphics[width=0.23\textwidth]{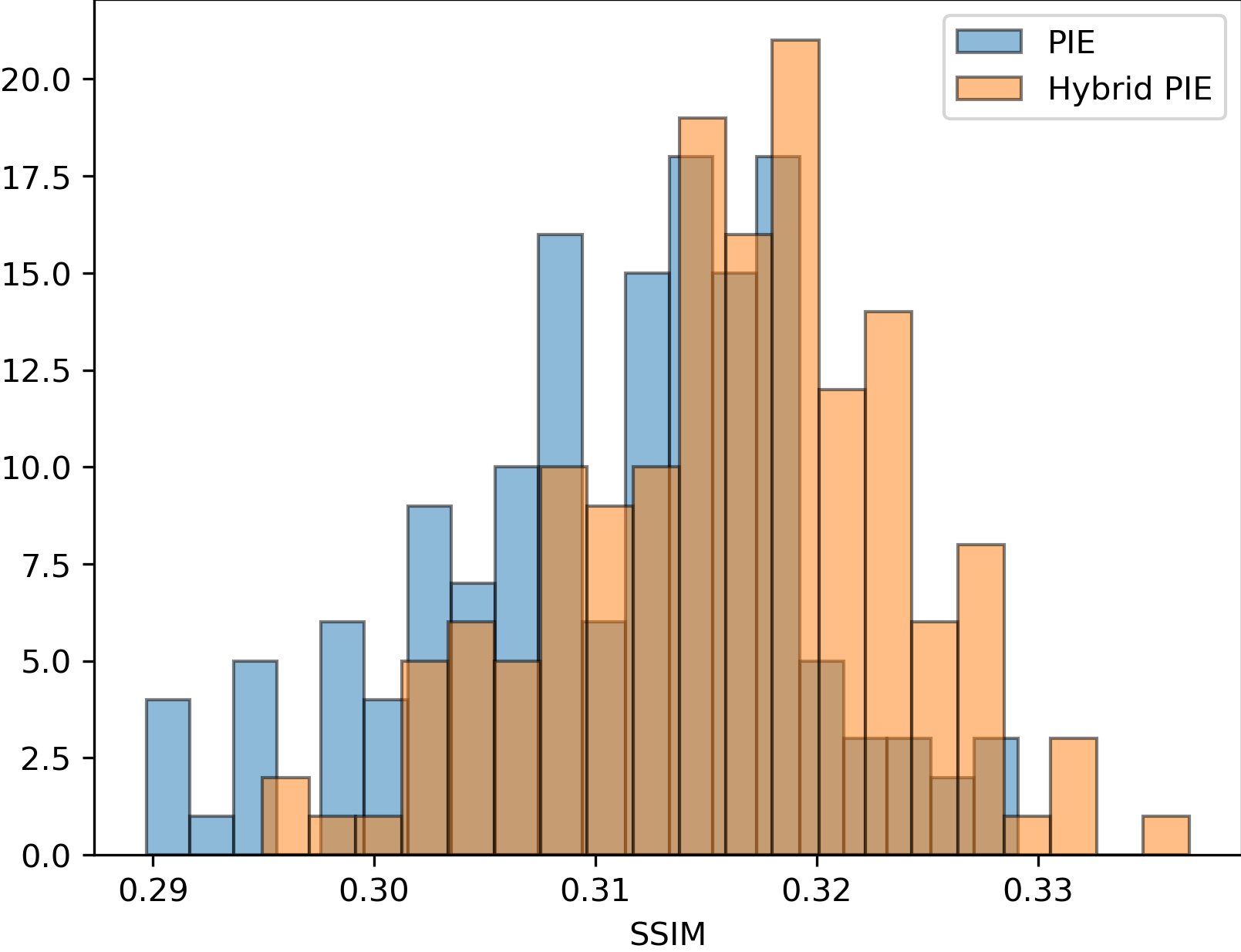}
        &
        \includegraphics[width=0.23\textwidth]{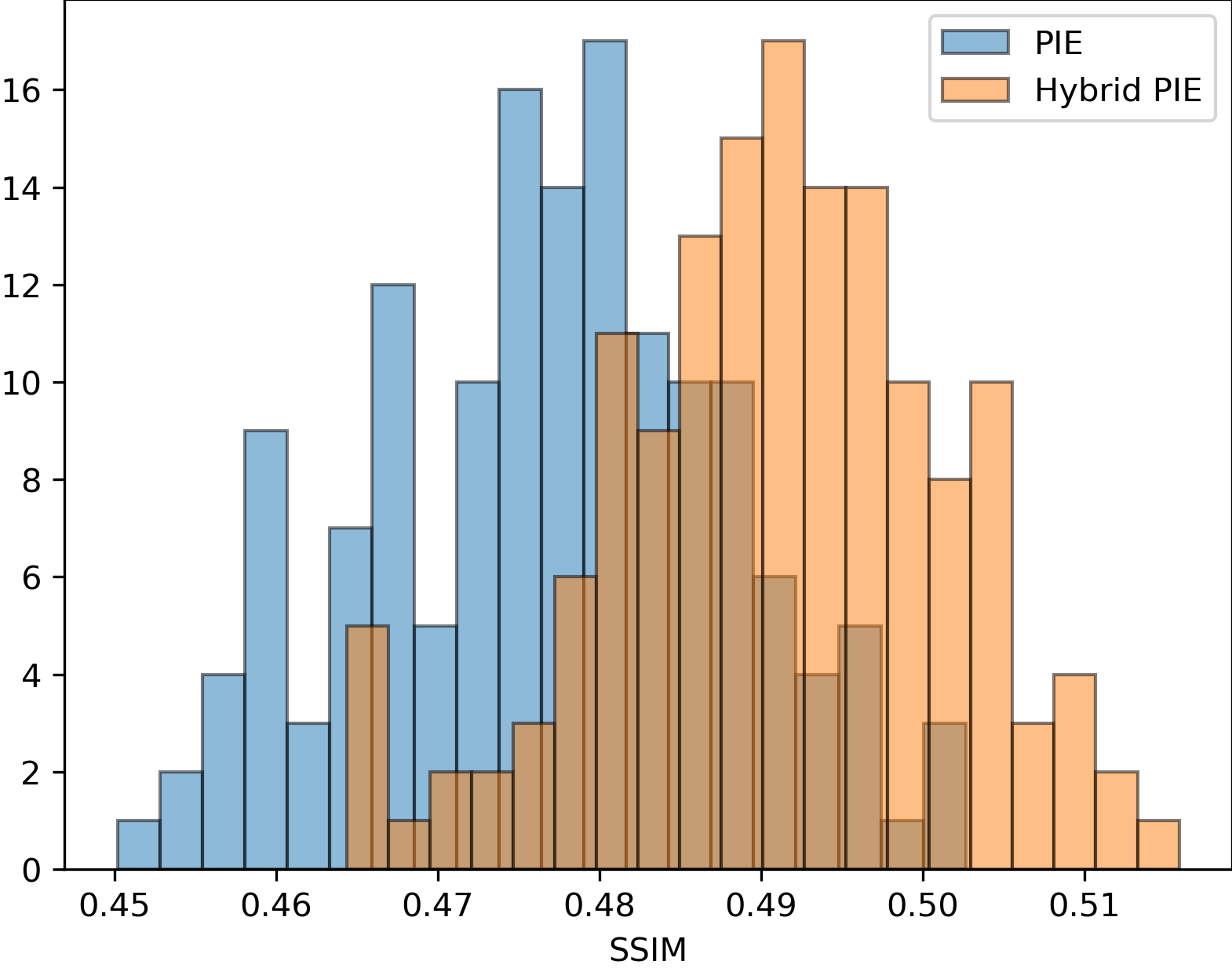}
        &
        \includegraphics[width=0.23\textwidth]{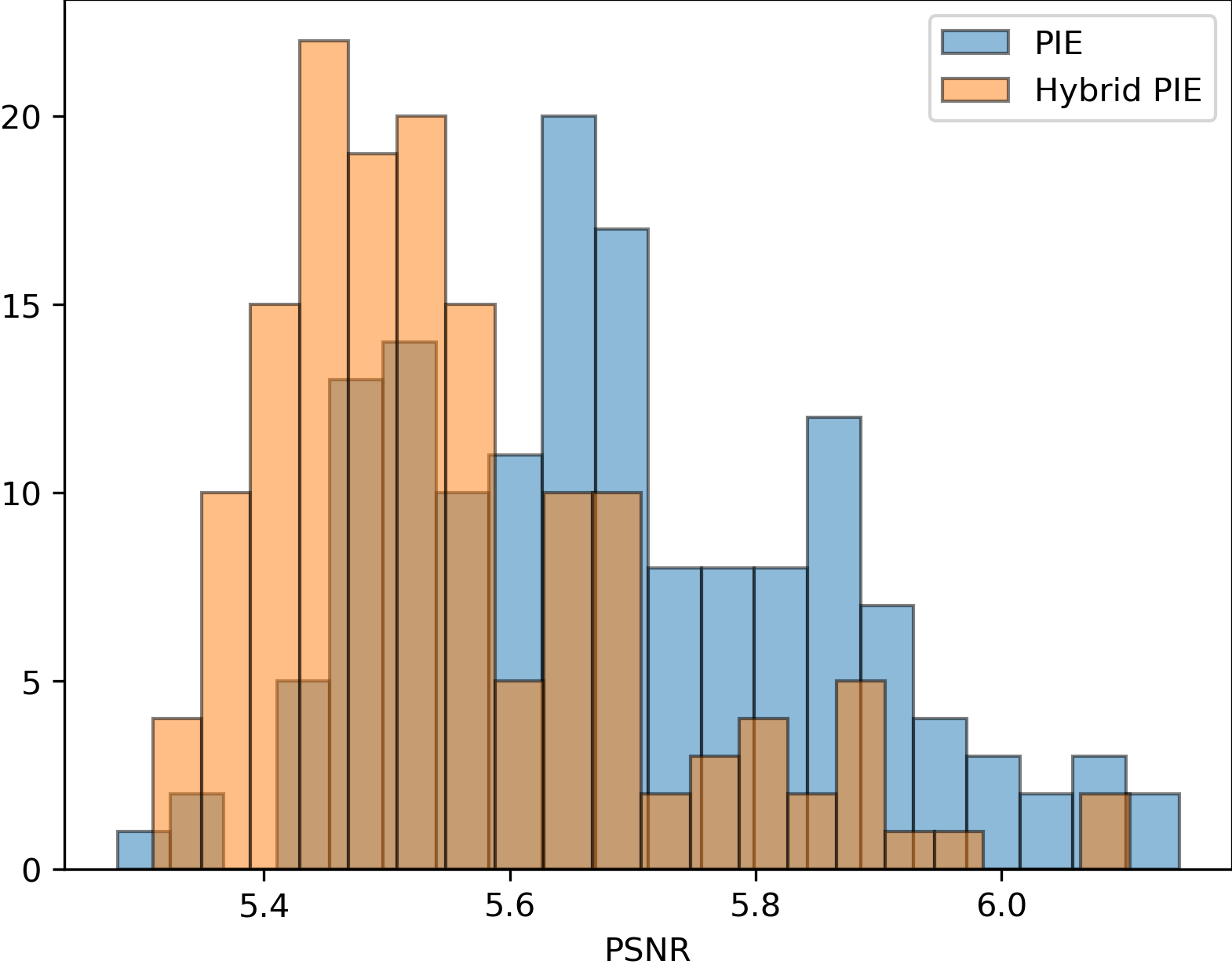}
        &
        \includegraphics[width=0.23\textwidth]{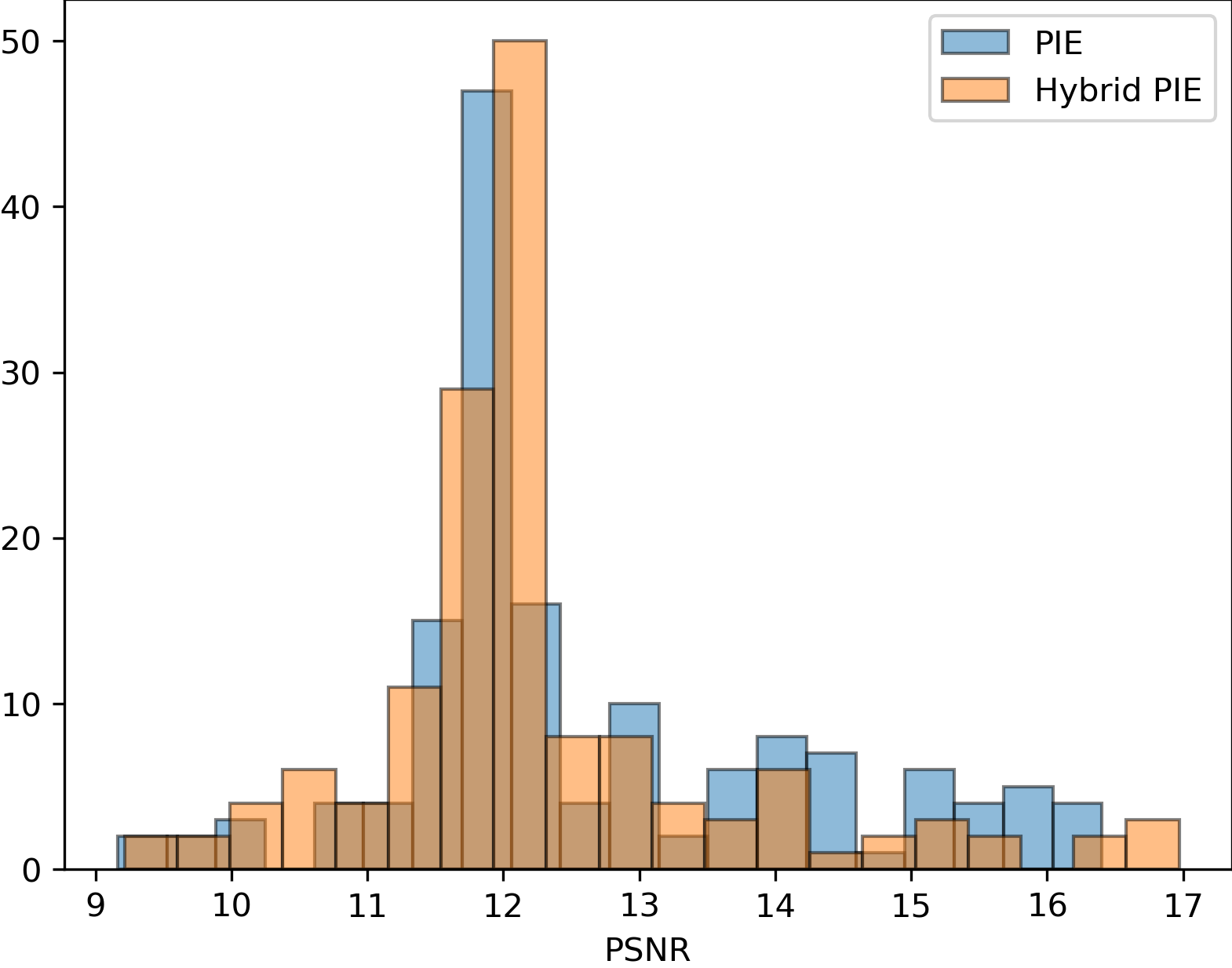}
    \end{tabular}
    \caption{\small{Histogram of the final SSIM and PSNR values for the reconstructed images. The blue histogram shows the SSIM/PSNR for PIE and the orange histogram shows the SSIM/PSNR for hybrid PIE.}}
    \label{fig: pie_ssim_psnr}
\end{figure}

\subsubsection{Distribution of Reconstruction Quality}
\label{subsubsec: nonblind_distribution_relerrs}
To test whether the proposed algorithm~\ref{alg: hybrid_ePIE} leads to improved reconstructions, we compare the distribution of relative errors of the reconstructed images between the PIE and the proposed hybrid PIE algorithms on small-scale images of size $512 \times 512$. 
We generate $150$ reconstructions using both algorithms with uniformly random initial guesses (here, we choose uniformly random in $[0,1]$ for the magnitude and $[0,\pi/2]$ for the phase).
We then compute the relative error of the final reconstructions
\begin{equation}
\label{eq: rel_error}
    \frac{||z - z_{\text{true}}||}{||z_{\text{true}}||}
\end{equation} where $z$ is our current reconstruction estimate and  $z_{\text{true}}$ is the ground truth. Similarly, we also compute the relative errors, peak signal-to-noise ratio (PSNR)~\cite{padmavathi2009performance}, and structural similarity index  SSIM~\cite{wang2004image} of the magnitude and phase to better understand the quality of reconstructions.

\begin{figure}[t]
    \centering
    \begin{tabular}{ccc}
        PIE Phase & PFT Warmup Phase & Hybrid PIE Phase
        \\
        \includegraphics[width=0.3\textwidth]{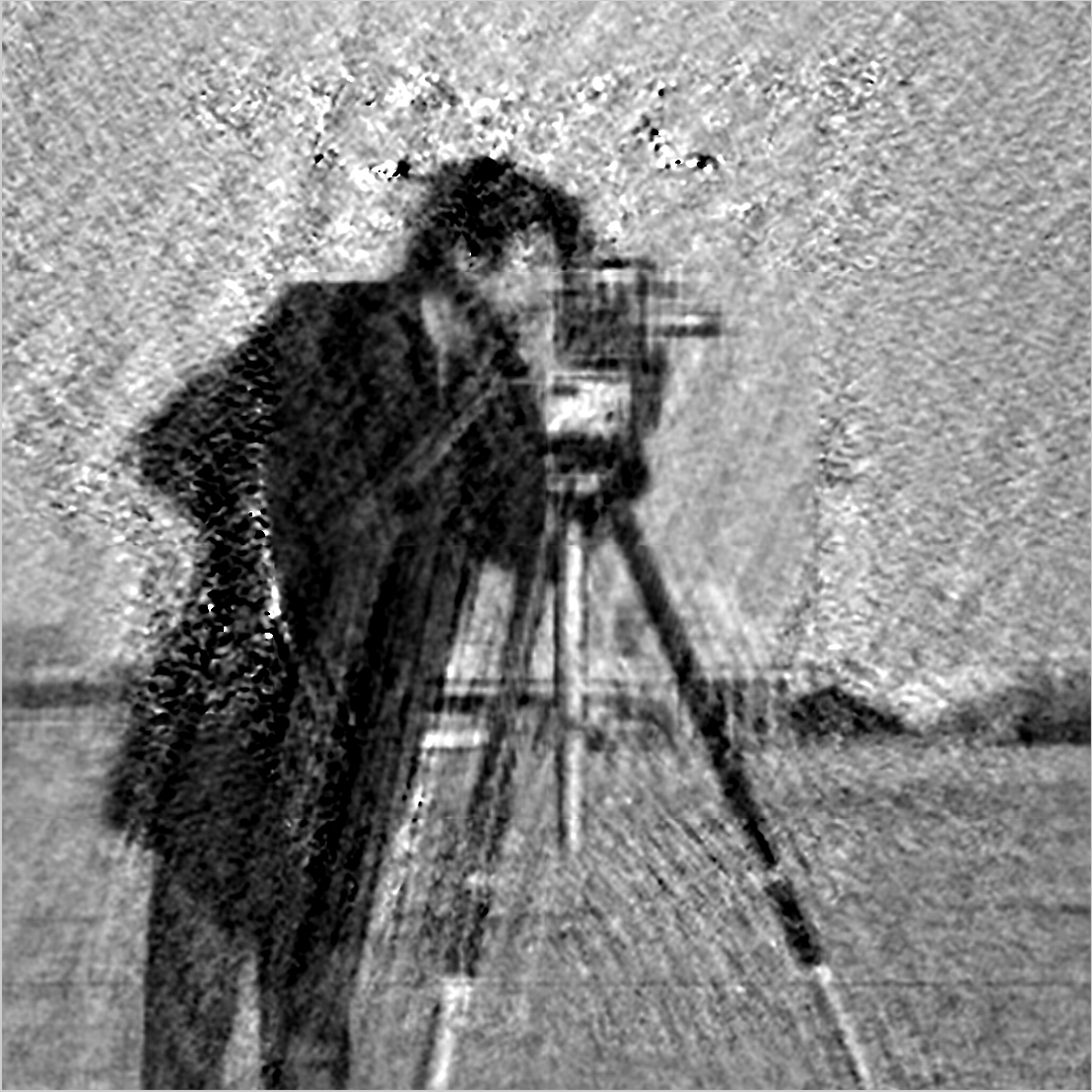}
        &
        \includegraphics[width=0.3\textwidth]{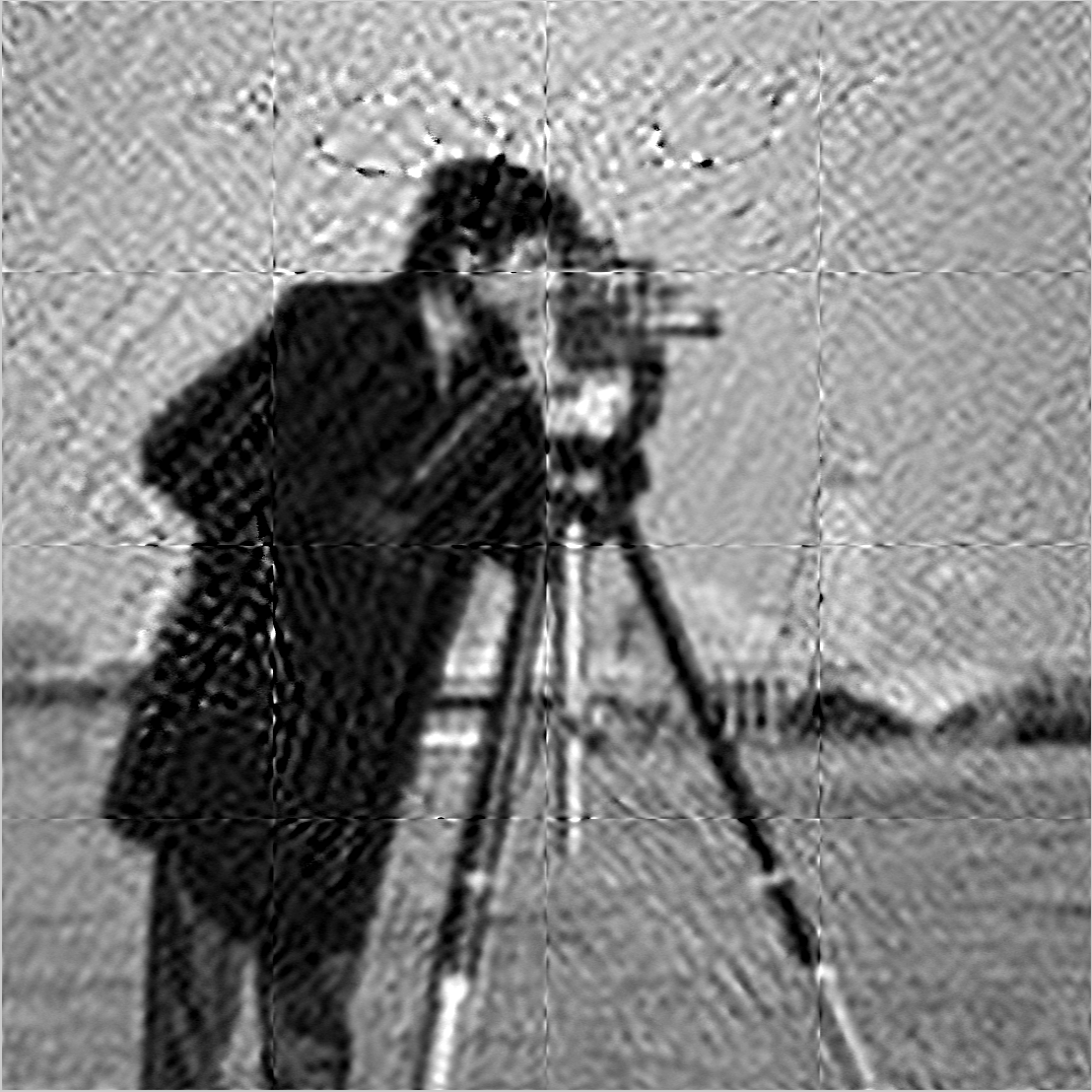}
        &
        \includegraphics[width=0.3\textwidth]{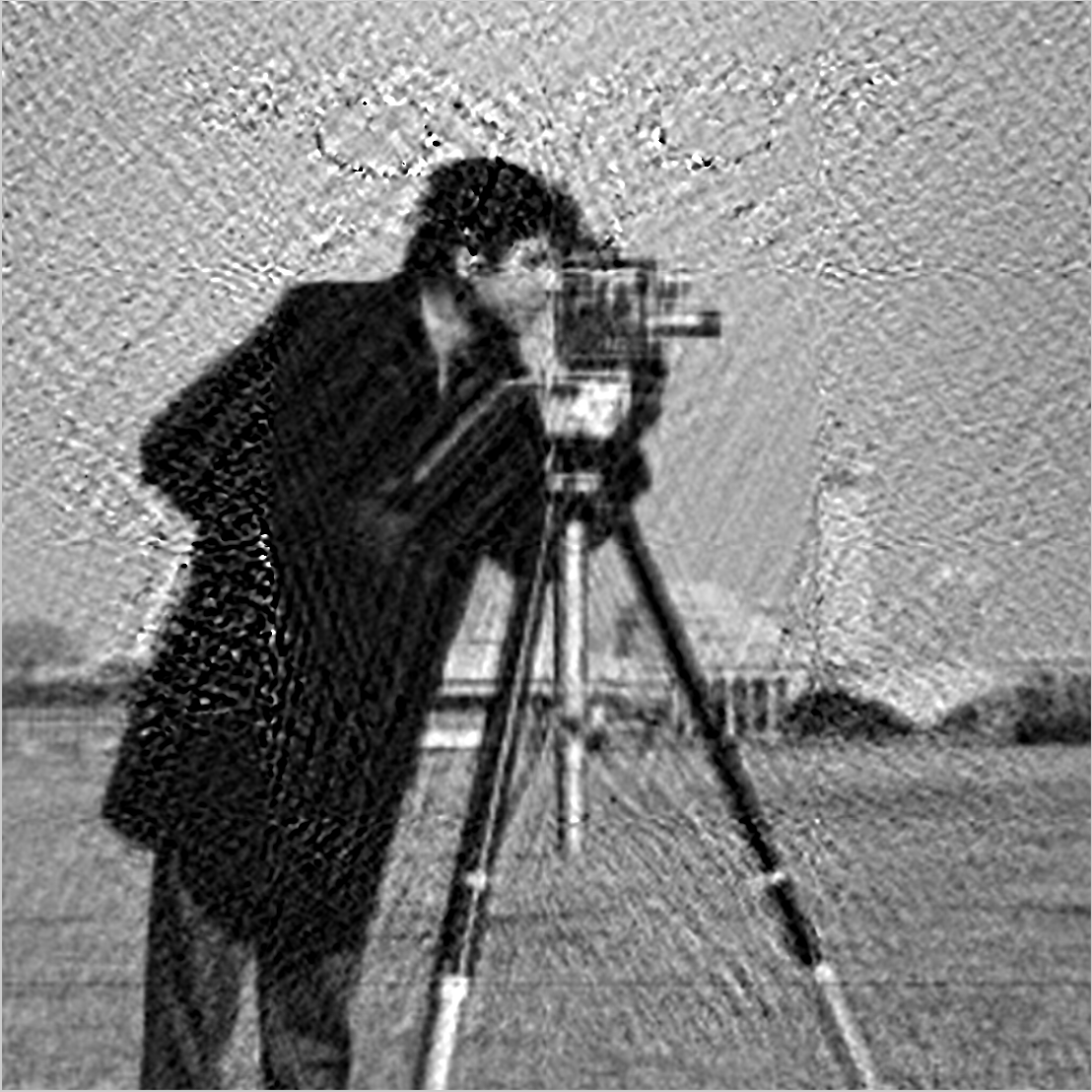}
        \\
        PIE Magnitude & PFT Warmup Magnitude & Hybrid PIE Magnitude
        \\
        \includegraphics[width=0.3\textwidth]{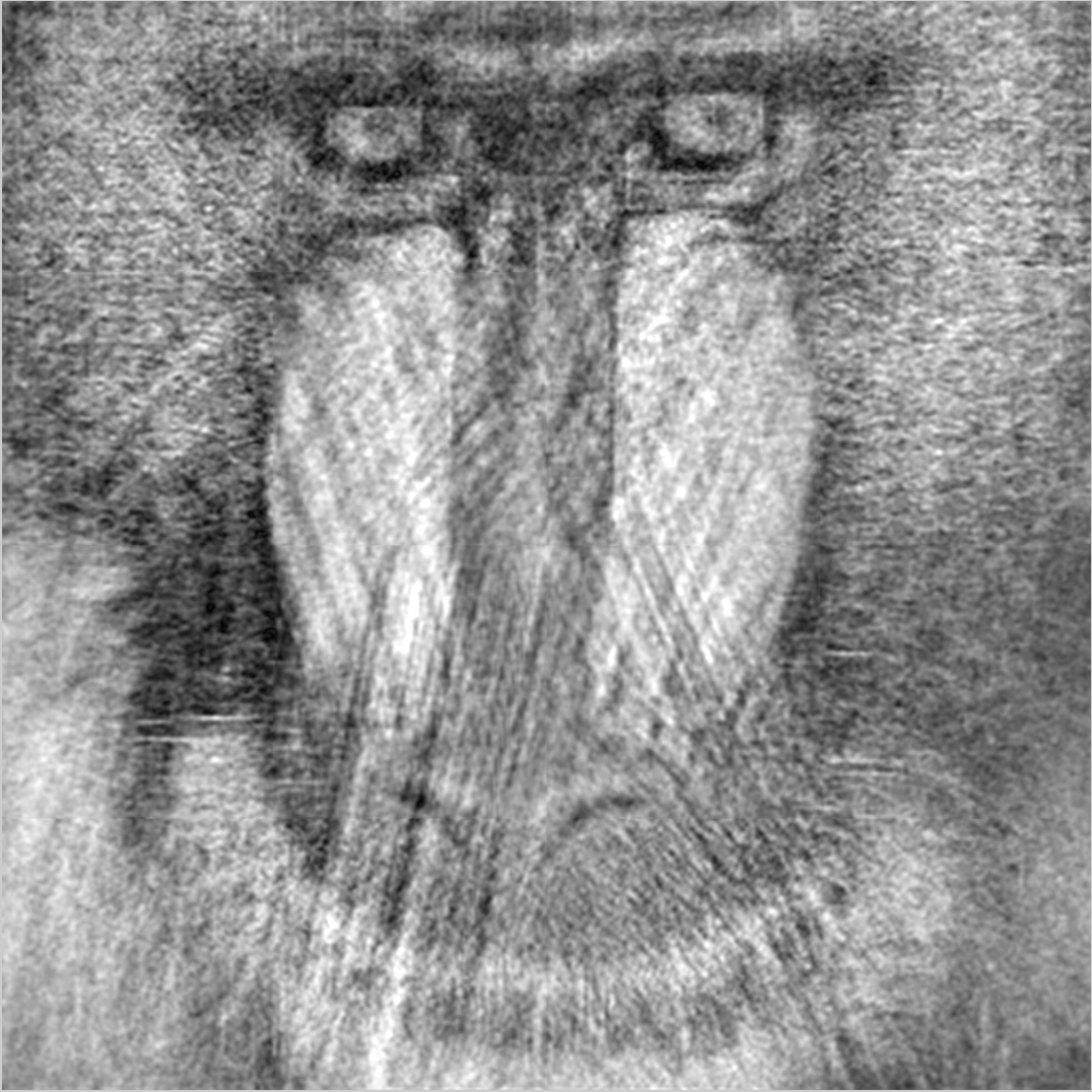}
        &
        \includegraphics[width=0.3\textwidth]{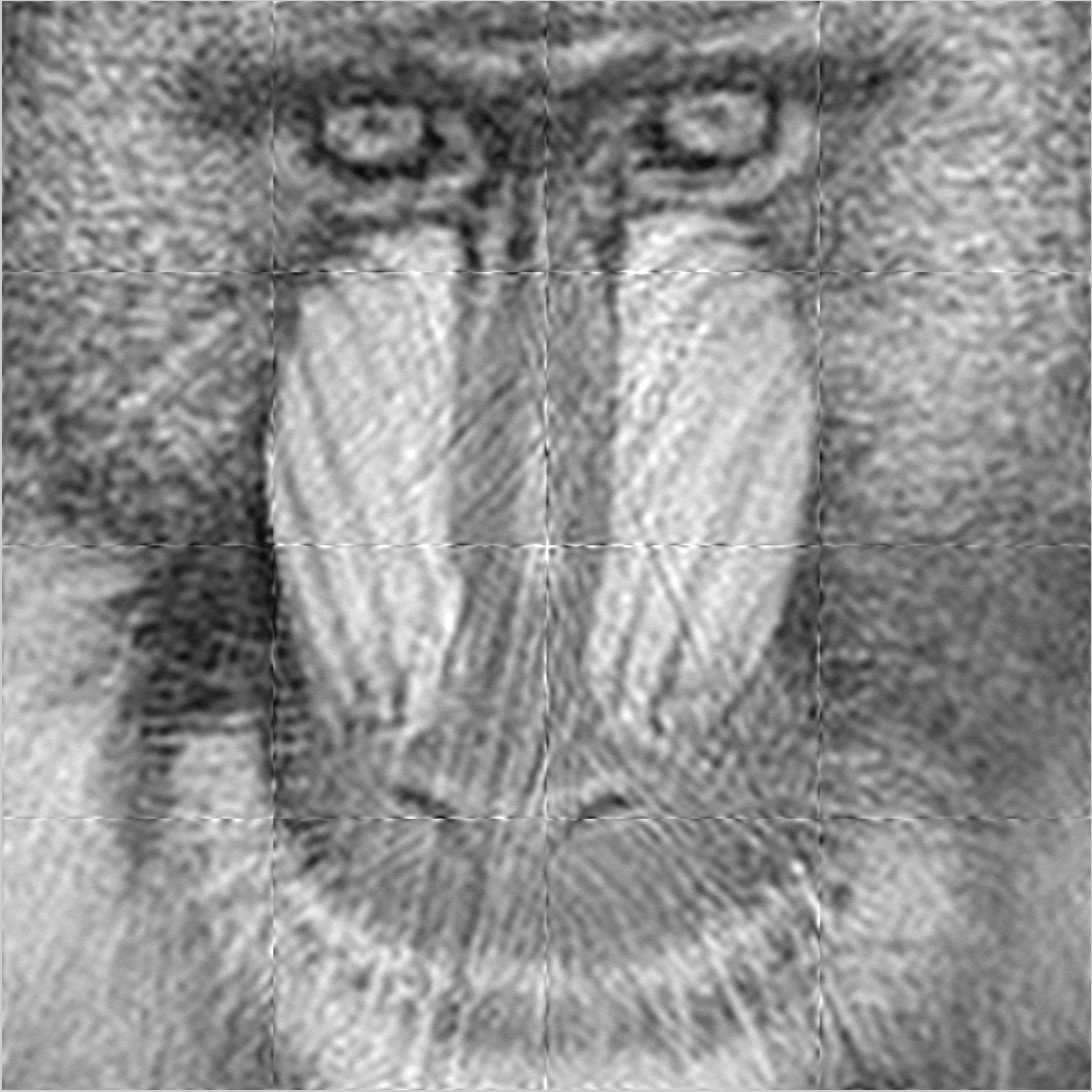}
        &
        \includegraphics[width=0.3\textwidth]{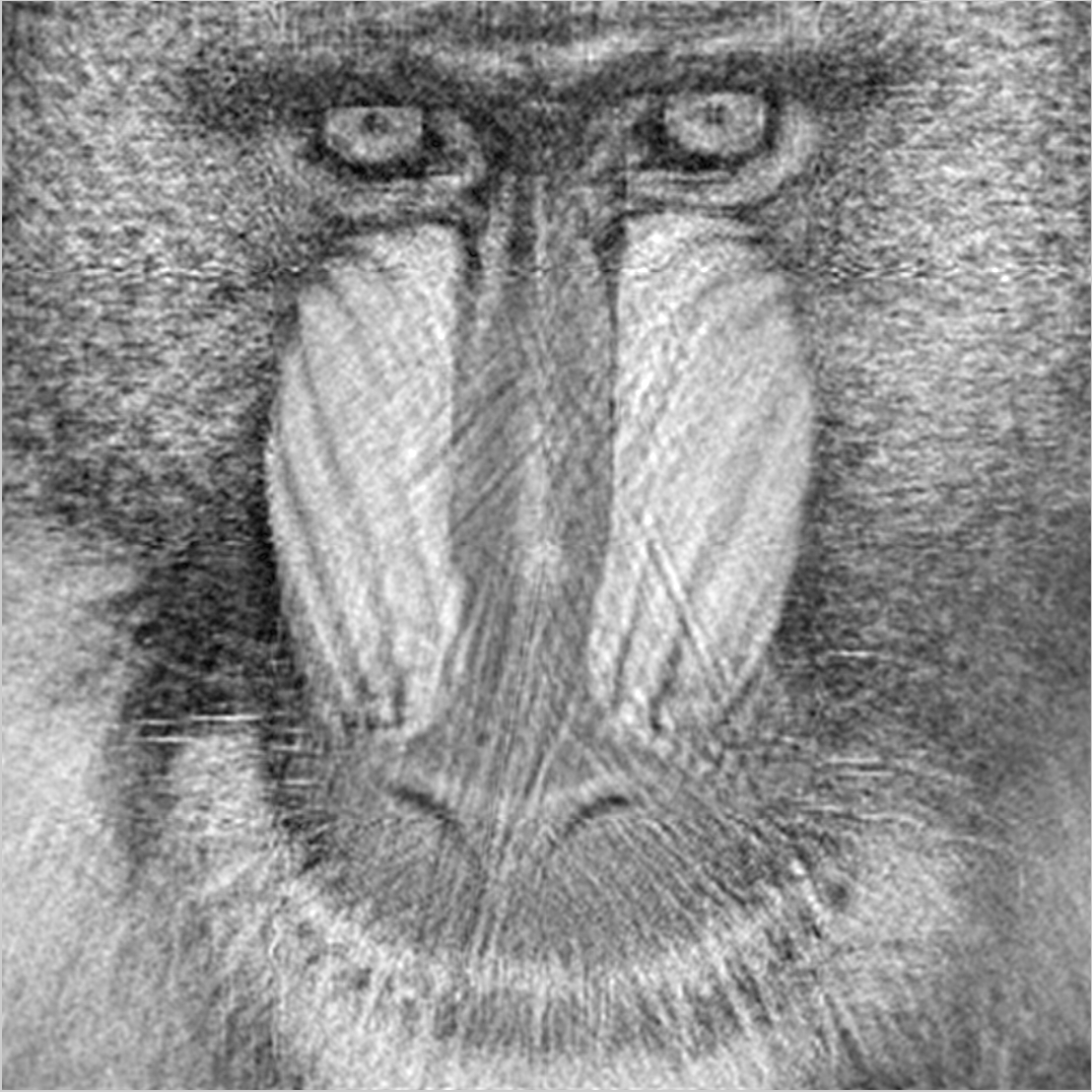}
    \end{tabular}
    \label{fig: non_blind_reconstructions}
    \caption{\small{Large-scale phase and magnitude reconstructions of the FFT-based PIE (first column), PFT-based PIE warmup (second column), and Hybrid PIE (third column. These reconstructions are for the large-scale ptychography problem where $n_1 = n_2 = 16384$.}}
    \label{fig: nonblind_reconstruction_images}
\end{figure}

Figure~\ref{fig: pie_rel_err} shows the relative errors between a) the full object, b) the magnitude of the object, and c) the phase of the object. The results show that the hybrid PIE approach leads to a similar, if not slightly better, distribution of reconstructions as the distribution of the traditional PIE. Similarly, Figure~\ref{fig: pie_ssim_psnr} shows a comparisons of the distributions of magnitude and phase SSIM and PSNRs (note for this figure that higher SSIM and PSNR is better). 

\subsubsection{Large-Scale Reconstruction}
\label{subsubsec: nonblind_large_scale}
To demonstrate the computational benefits of using a PFT-based PIE in the proposed hybrid algorithm, we run a larger experiment where $n_1 = n_2 = m_1 = m_2 = 16384$, and $\widetilde{m}_1 = \widetilde{m}_2 = p_1 = p_2 = 64$, leading to a cropped image of size $128 \times 128$. The degree of the approximating polynomial in this setting is given by degrees $r_1 = r_2 = 13$, which is obtained from the precomputed values $\xi(\varepsilon, r)$ defined in Section~\ref{subsec: pft_offline_derivation}.
Here, we run the PFT-based PIE algorithm until a tolerance of $\epsilon_{\text{pft}} = 10^{-2}$ or a maximum number of $50$ iterations is reached. Afterwards, we run the standard FFT-based PIE until a tolerance of $\epsilon = 5 \times 10^{-4}$ or a maximum number of $100$ iterations is reached. Here, we choose step sizes $\beta = 10$ for the FFT-based PIE and $\beta_{PFT} = 10^{-3}$ for the PFT-based PIE. These were chosen tuned using a logarithmic gridsearch over the set $\{10^{-6}, 10^{-5}, \ldots, 10^2, 10^{3}\}$, and afterwards using a standard gridsearch starting from the currently chosen parameter, e.g., $10^{-2}$, to the next order of magnitude, e.g., $10^{-1}$.
The intuition for choosing a higher tolerance for the PFT-based PIE is based on the fact that the PFT-based PIE will primarily capture large features, whereas the FFT-based PIE needs more iterations to capture fine details. For all reconstructions, we use total variation regularization~\cite{rudin1992nonlinear} with parameters $\lambda = 10^{-6}$ for the FFT-based PIE and $\lambda = 10^2$ for the PFT-based PIE. These are the parameters that led to the lowest relative error and were found over a logarithmic gridsearch over the set~$\{10^{-6}, 10^{-5}, \ldots, 10^2, 10^{3}\}$. Interestingly, the PFT-based PIE benefits from larger total variation.

\begin{figure}[t]
    \small
    \centering
    \begin{tabular}{cccc}
        \includegraphics[width=0.23\textwidth]{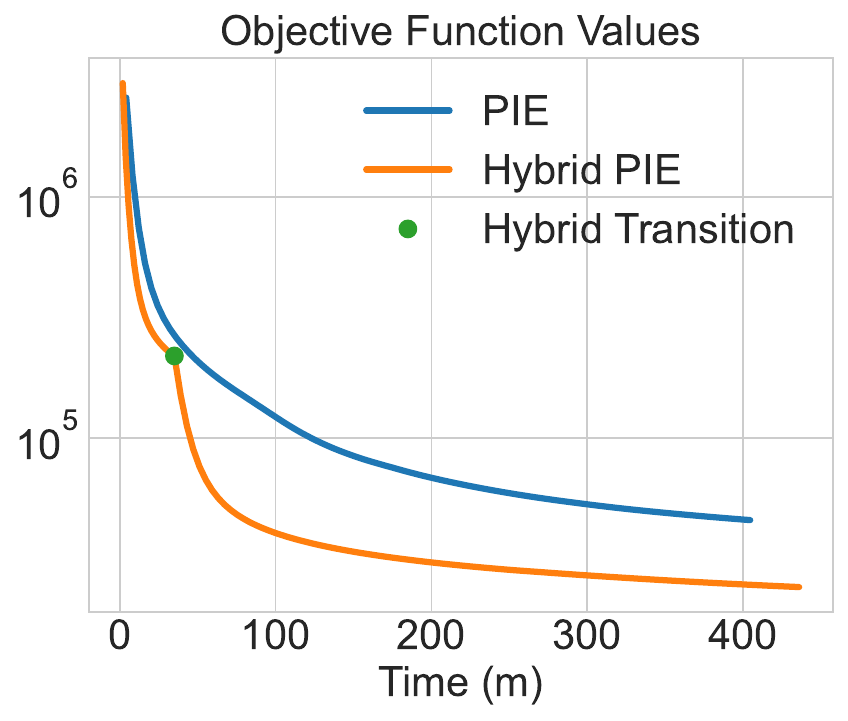}
        &
        \includegraphics[width=0.23\textwidth]{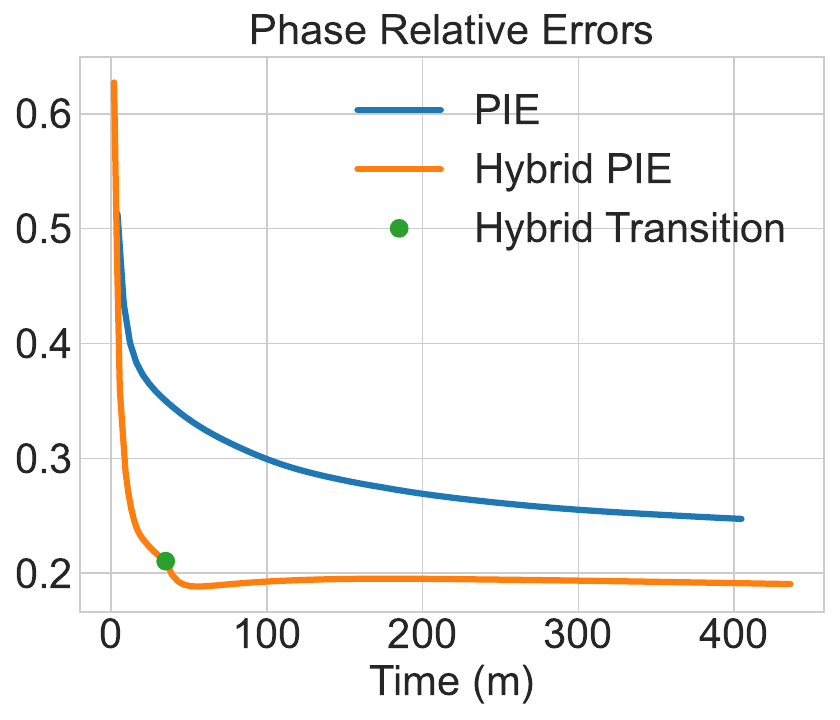}
        & 
        \includegraphics[width=0.23\textwidth]{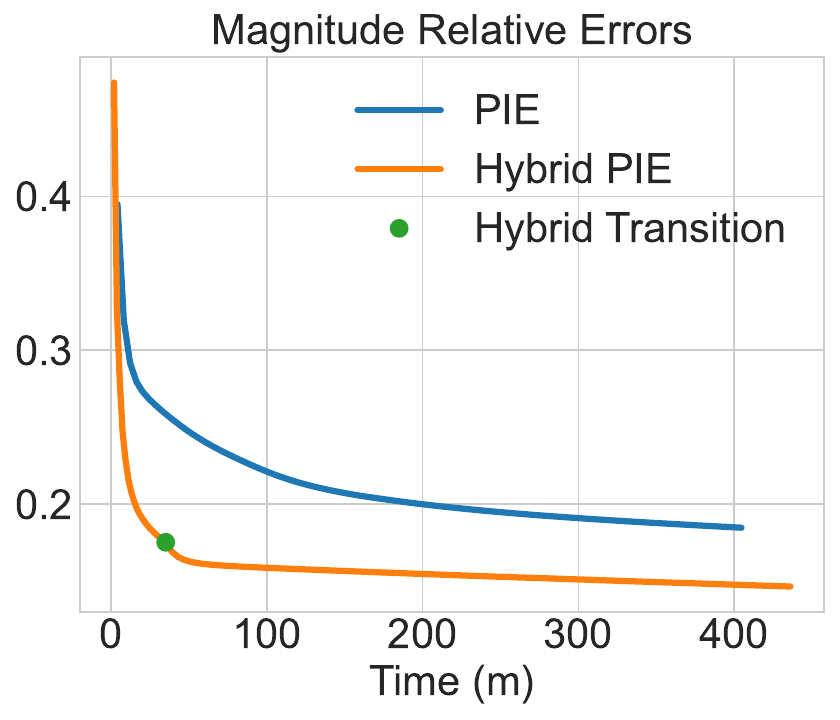}
        &
        \includegraphics[width=0.23\textwidth]{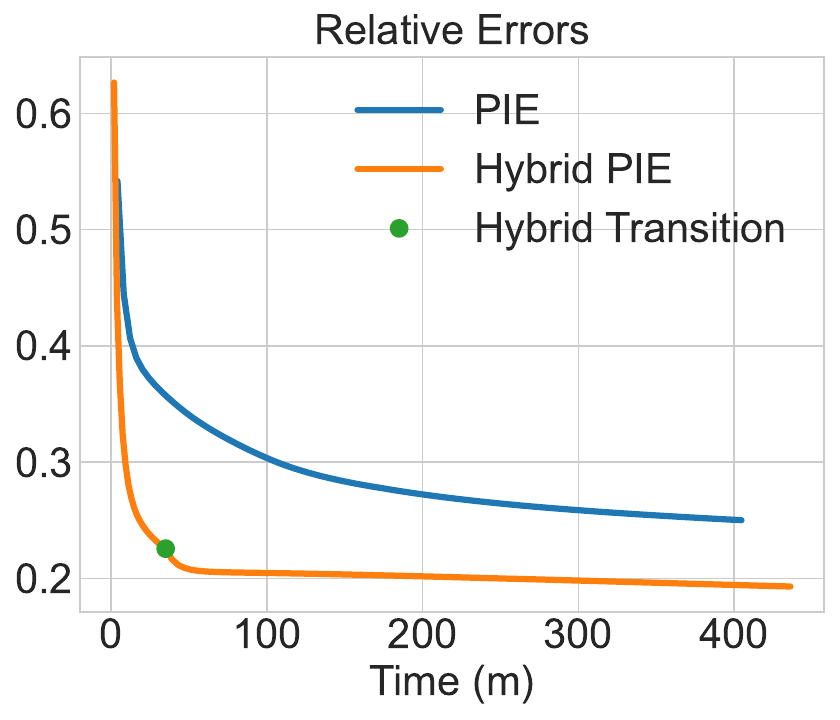}
    \end{tabular}
    \caption{\small{Objective function values (first column), phase relative errors (second column), magnitude relative errors (third column), and relative errors (fourth column) for a large-scale ptychographic reconstruction explained in Section~\ref{subsubsec: nonblind_experimental_setup}. The x-axis represents time in minutes.
    Green dot represents the transition from the PFT-based PIE to the FFT-based PIE in the hybrid algorithm.
    }
    }
    \label{fig: non_blind_rel_errors}
\end{figure}

\begin{figure}[t]
    \small
    \centering
    \begin{tabular}{cccc}
        \includegraphics[width=0.23\textwidth]{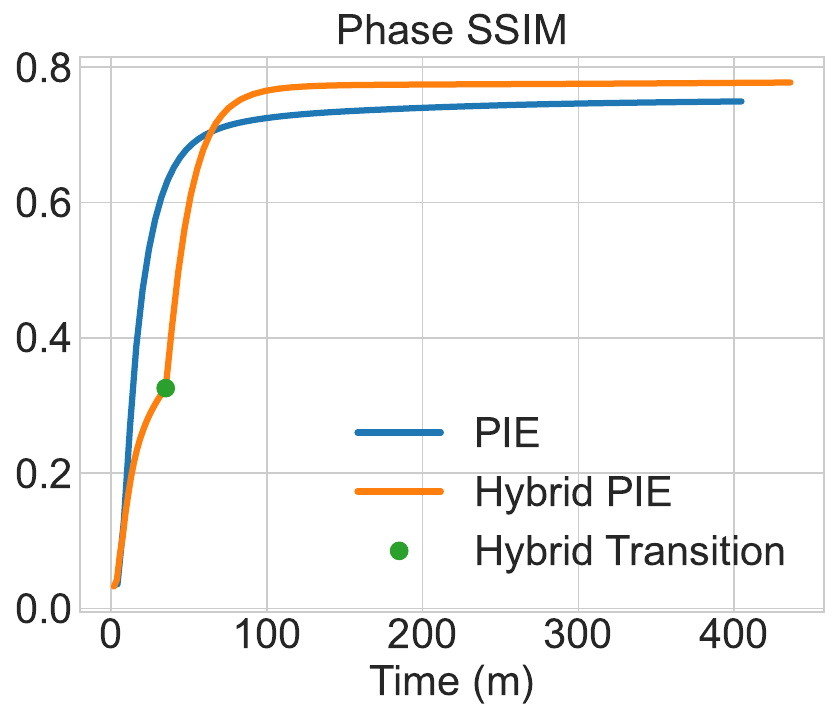}
        &
        \includegraphics[width=0.23\textwidth]{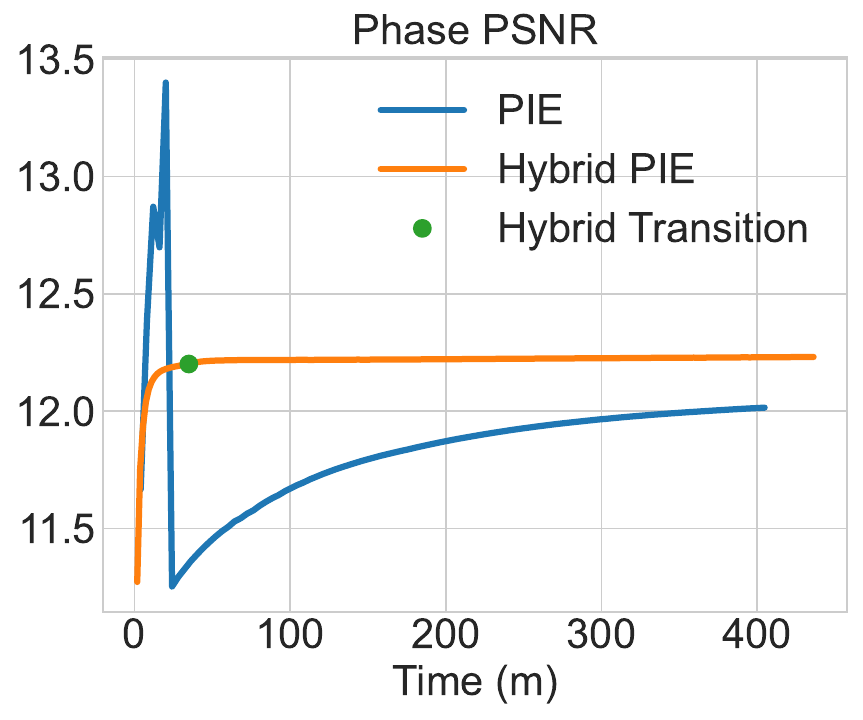}
        &
        \includegraphics[width=0.23\textwidth]{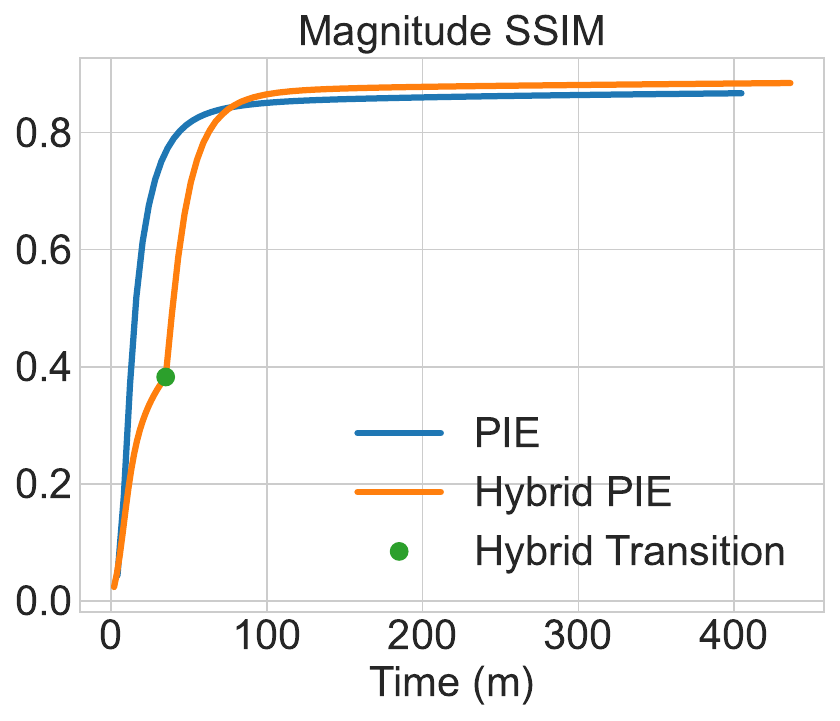}
        &
        \includegraphics[width=0.23\textwidth]{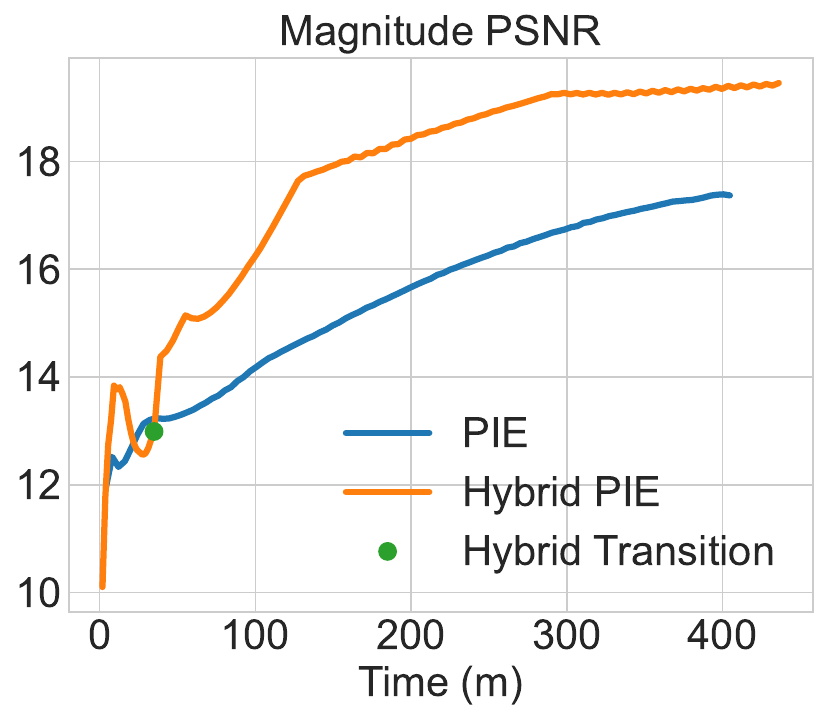}
    \end{tabular}
    \caption{SSIM and PSNR for PIE (blue) and Hybrid PIE (orange) for the phase and magnitude over time (in minutes). The green dot represents the transition from the PFT-based PIE to the FFT-based PIE in the hybrid algorithm.}
    \label{fig: nonblind_ssim_psnr}
\end{figure}

In Figure~\ref{fig: nonblind_reconstruction_images}, we show the reconstructions of the magnitude and phase generated by PIE, the initial PFT-based PIE, and the proposed hybrid PIE.
In Figure~\ref{fig: non_blind_rel_errors}, we show the objective function values of the full reconstruction (first column), the phase (second column), the magnitude (third column), and the relative errors (fourth column). Here, the x-axis represents time in minutes.
Similar results are shown in Figure~\ref{fig: nonblind_ssim_psnr}, where we only show the phase and magnitude structural similarity index (SSIM)~\cite{wang2004image} and peak signal to noise ratio (PSNR)~\cite{poobathy2014edge}. We observe that using the PFT-based PIE as warm up leads to improved SSIM and PSNR on the phase and magnitudes of the reconstruction.

\subsection{Blind Ptychographic Retrieval}
\label{subsec: blind_retrieval_experiment}

\begin{figure}[t]
    \setlength\tabcolsep{1 pt}
    \centering
    \begin{tabular}{cccccccc}
        $Q_1$ & $Q_2$ & $Q_3$ & $Q_4$ & $Q_5$ & $Q_6$ & $Q_7$ & $Q_8$
        \\
        \includegraphics[width=0.12\textwidth]{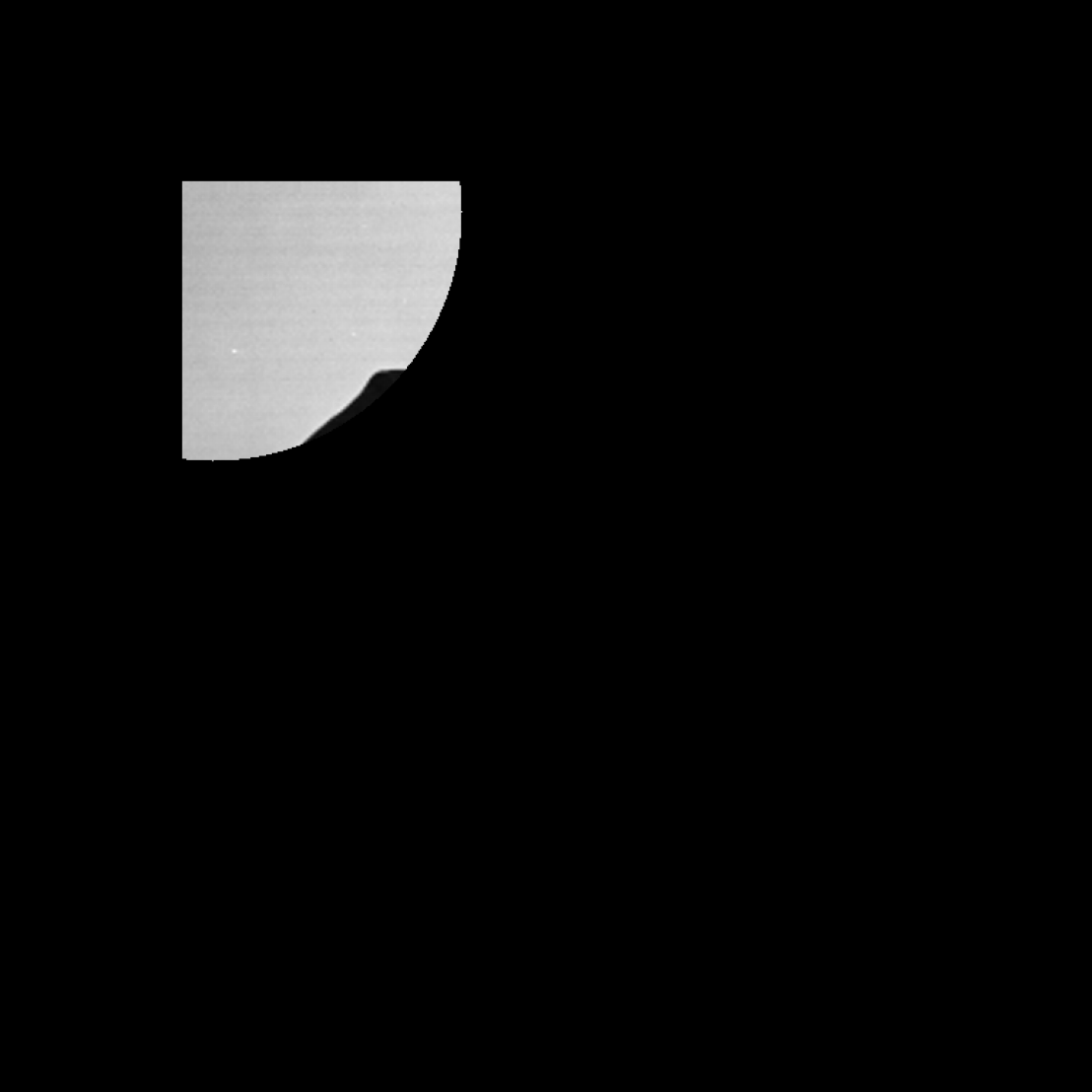}
        &
        \includegraphics[width=0.12\textwidth]{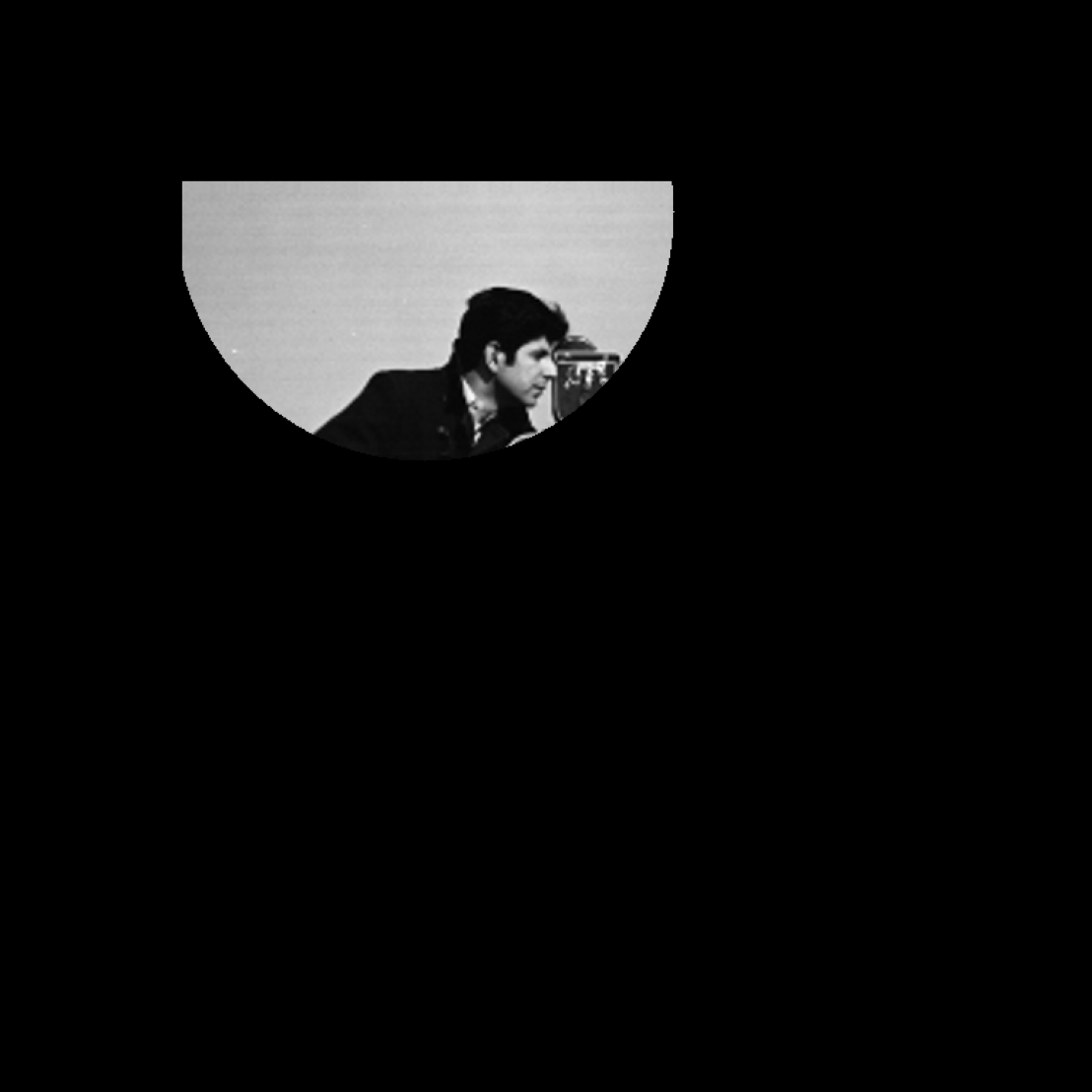}
        &
        \includegraphics[width=0.12\textwidth]{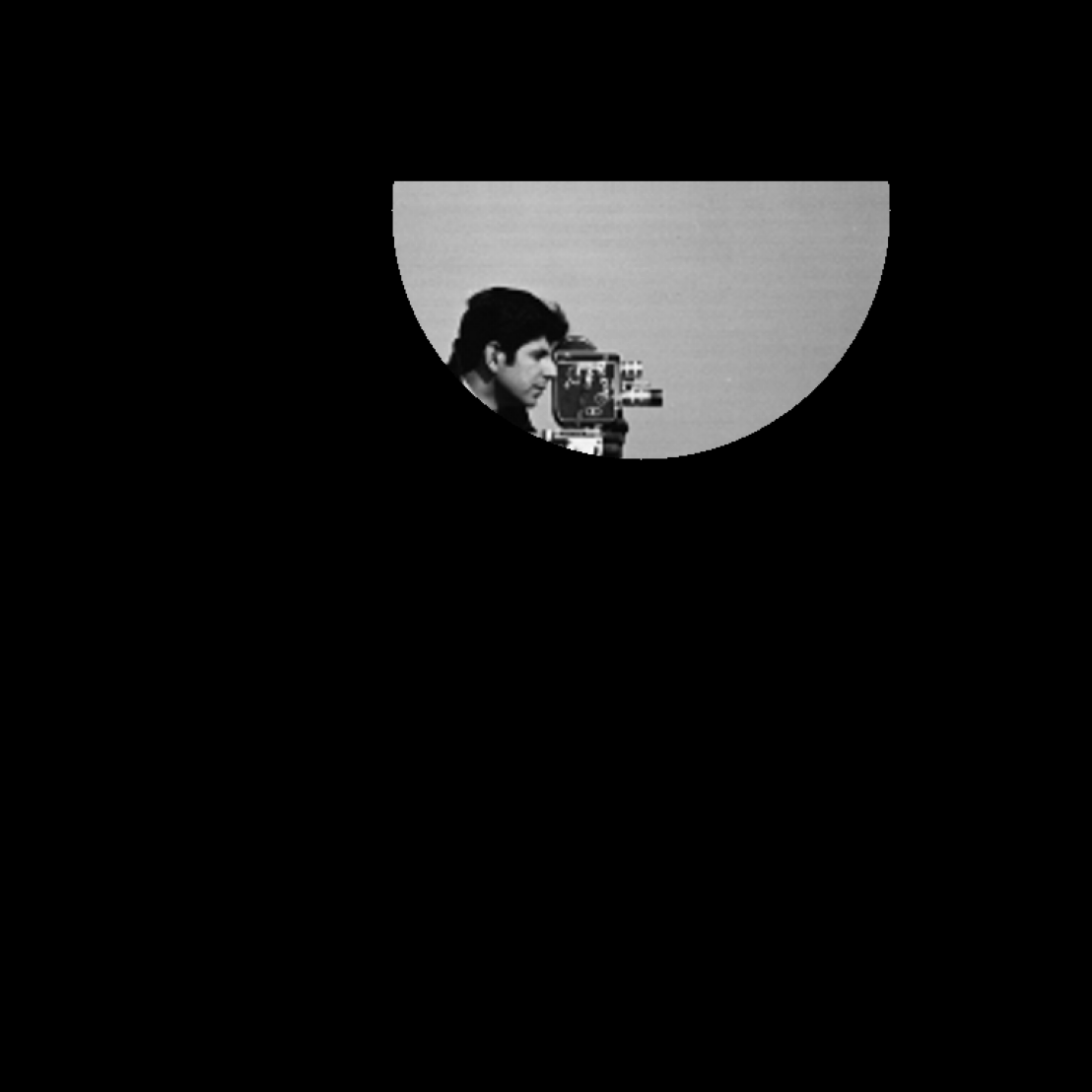}
        &
        \includegraphics[width=0.12\textwidth]{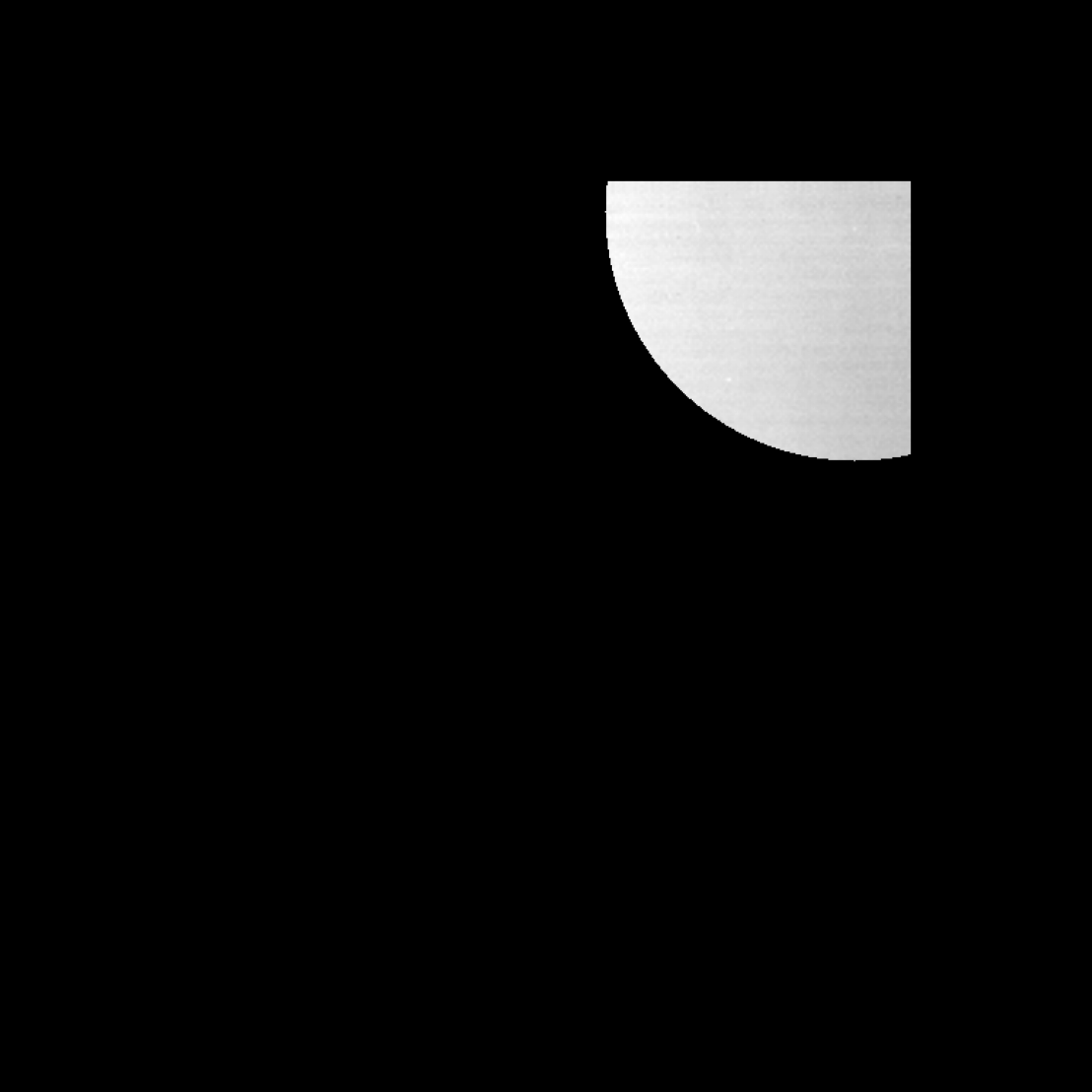} 
        &
        \includegraphics[width=0.12\textwidth]{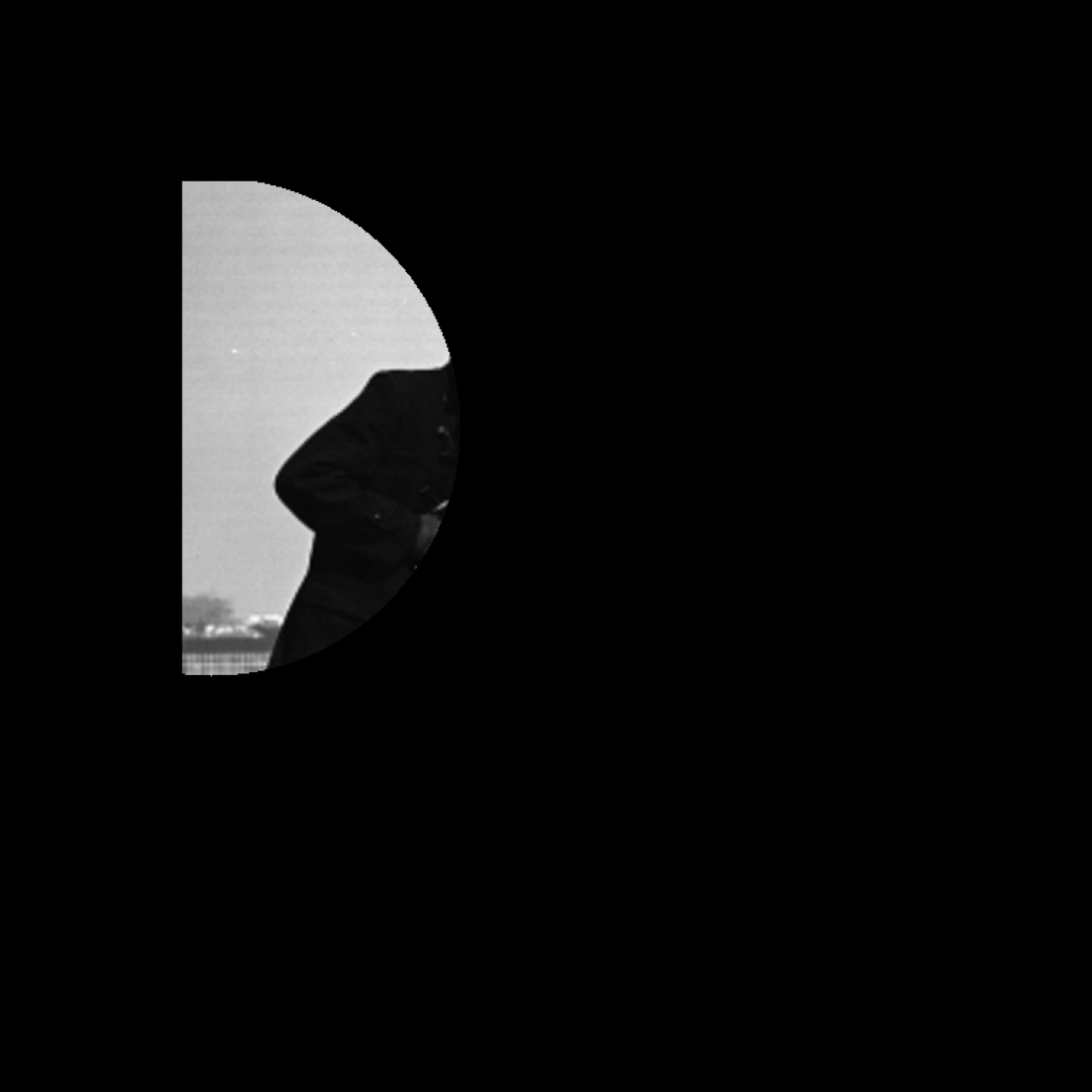}
        &
        \includegraphics[width=0.12\textwidth]{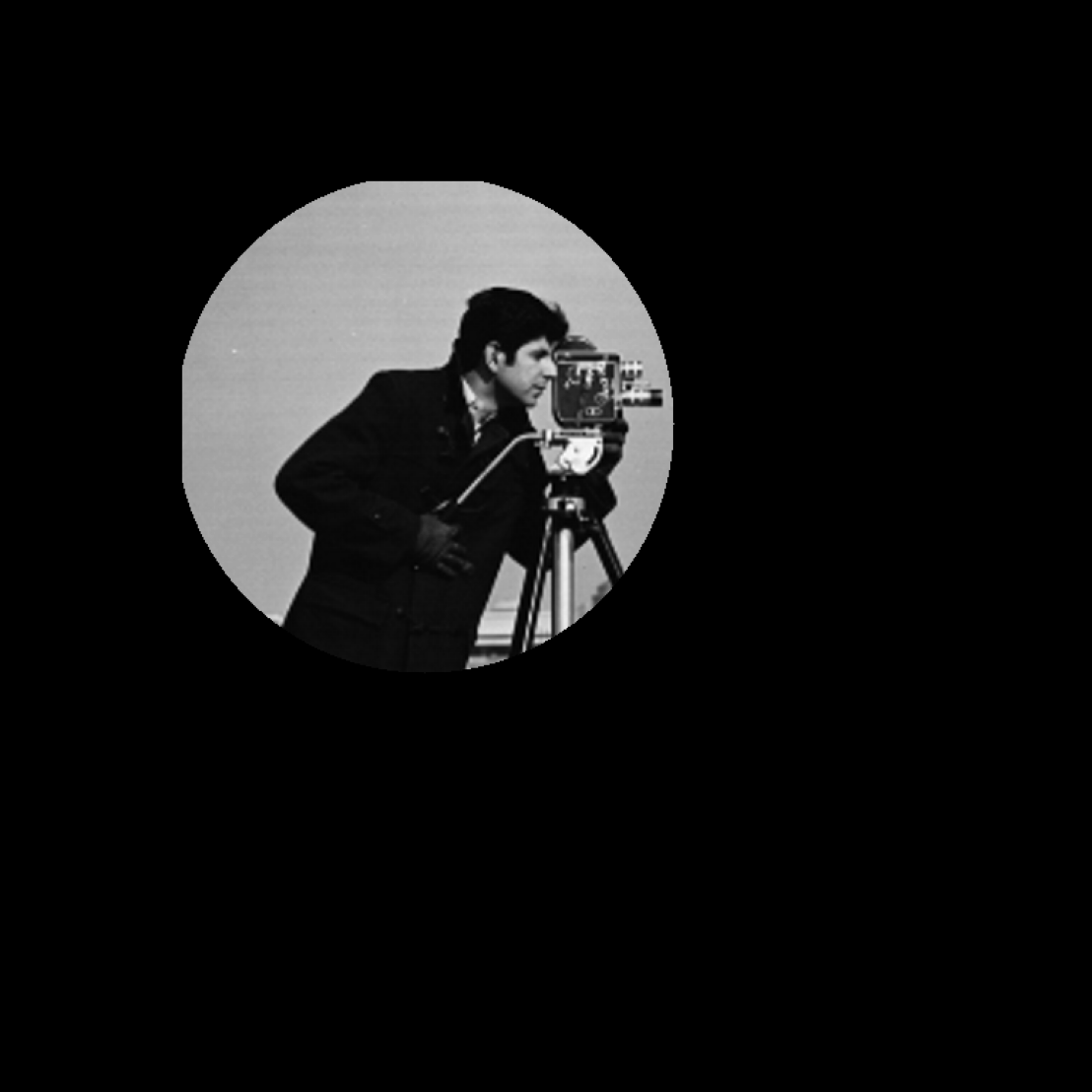}
        &
        \includegraphics[width=0.12\textwidth]{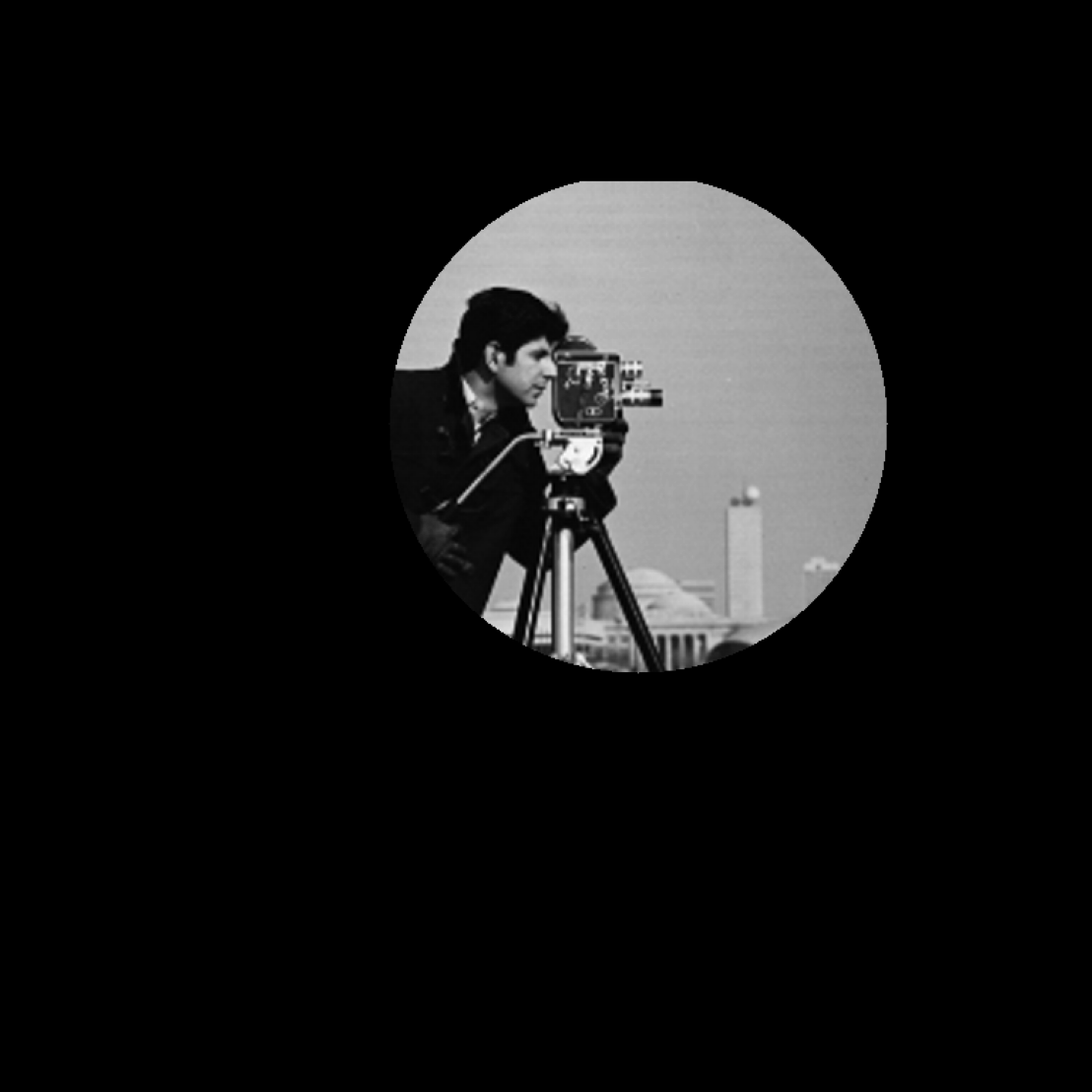}
        &
        \includegraphics[width=0.12\textwidth]{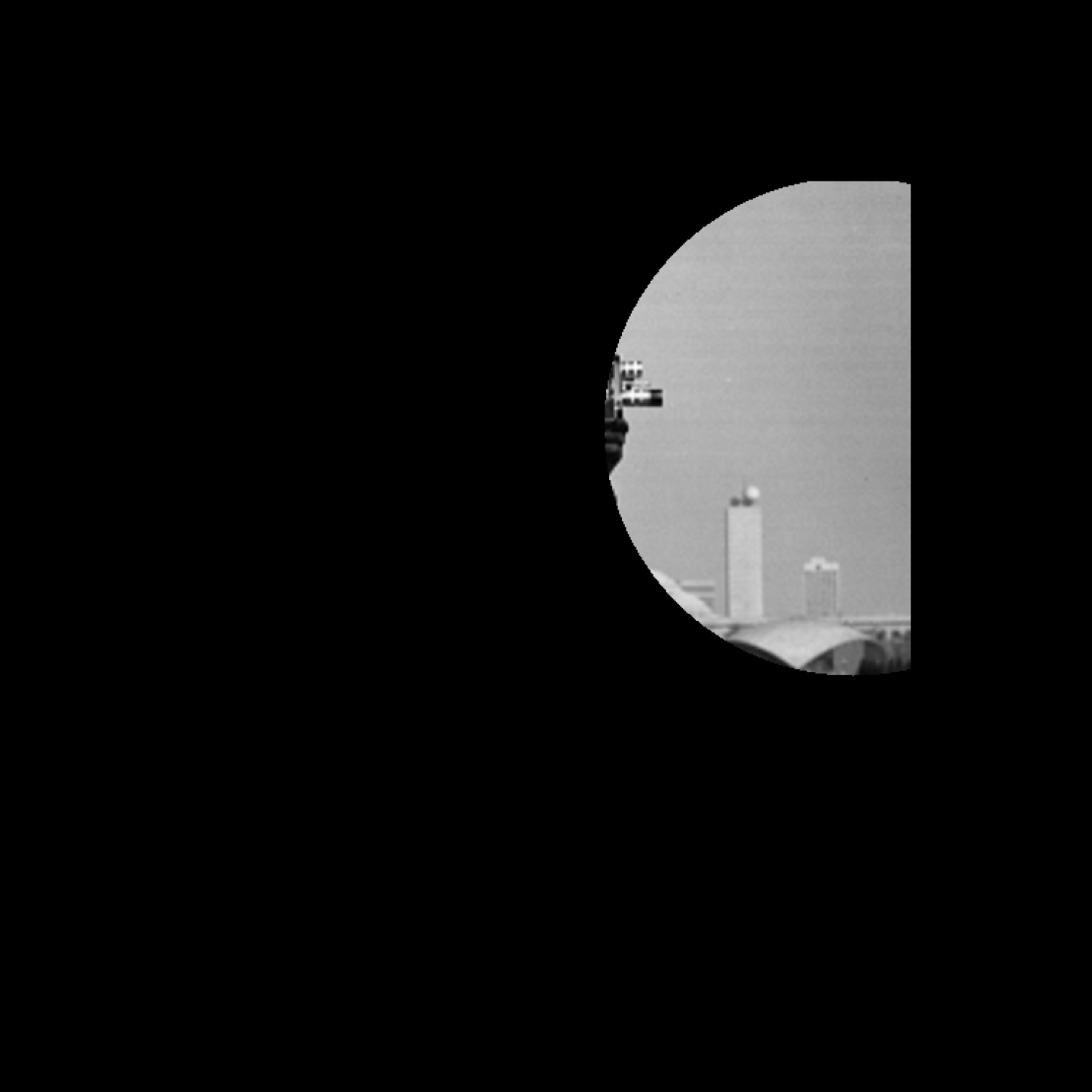} \\
        $Q_9$ & $Q_{10}$ & $Q_{11}$ & $Q_{12}$ & $Q_{13}$ & $Q_{14}$ & $Q_{15}$ & $Q_{16}$
        \\
        \includegraphics[width=0.12\textwidth]{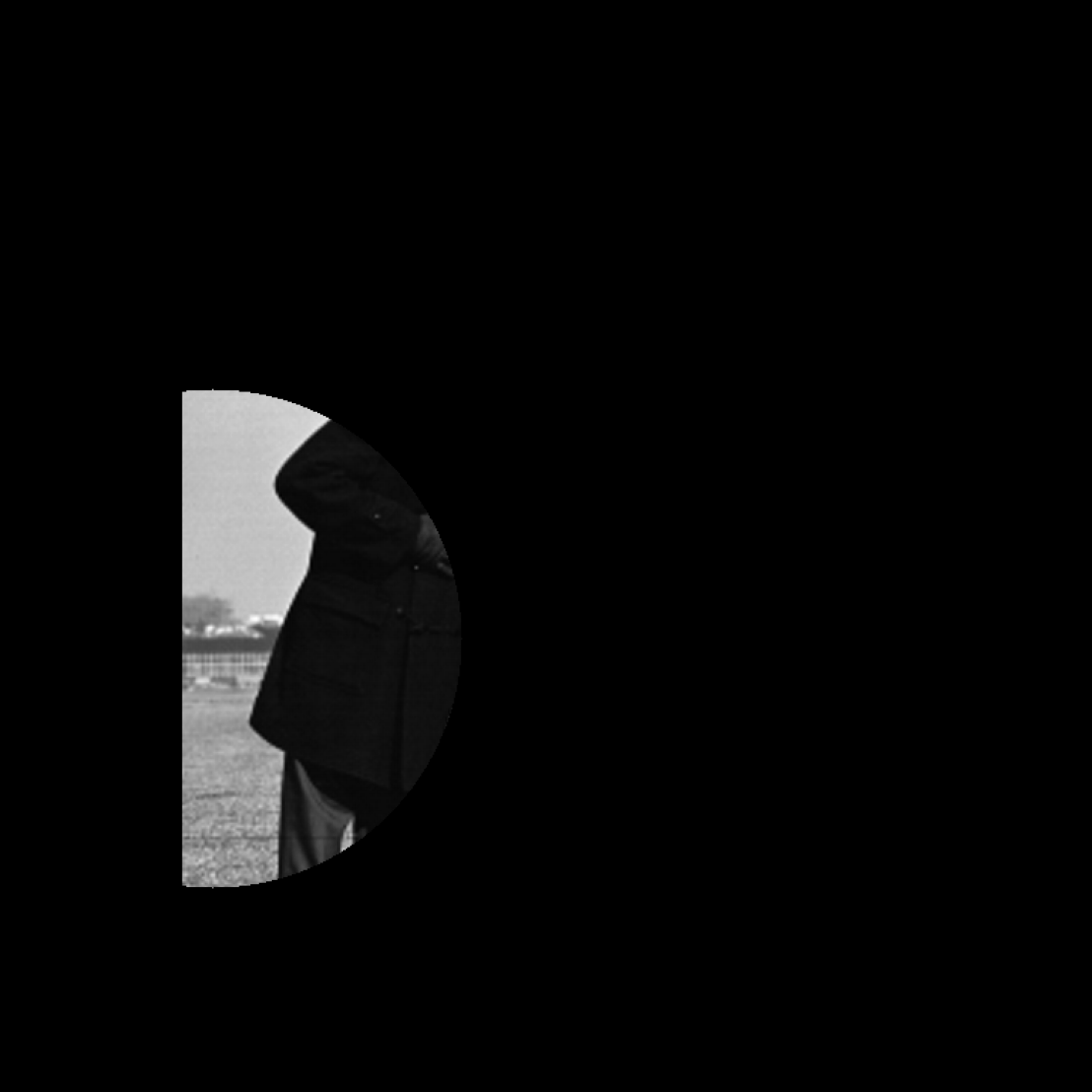}
        &
        \includegraphics[width=0.12\textwidth]{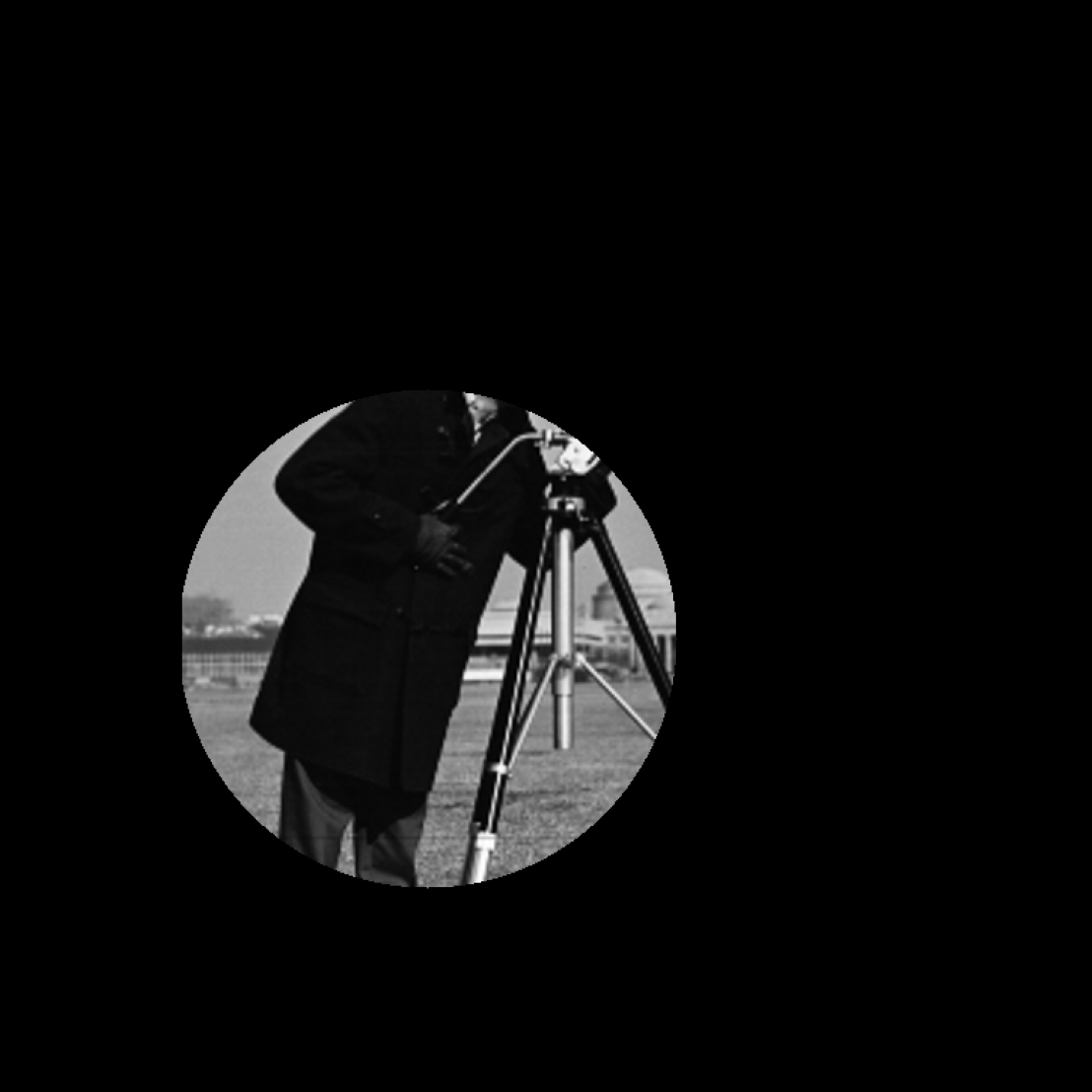}
        &
        \includegraphics[width=0.12\textwidth]{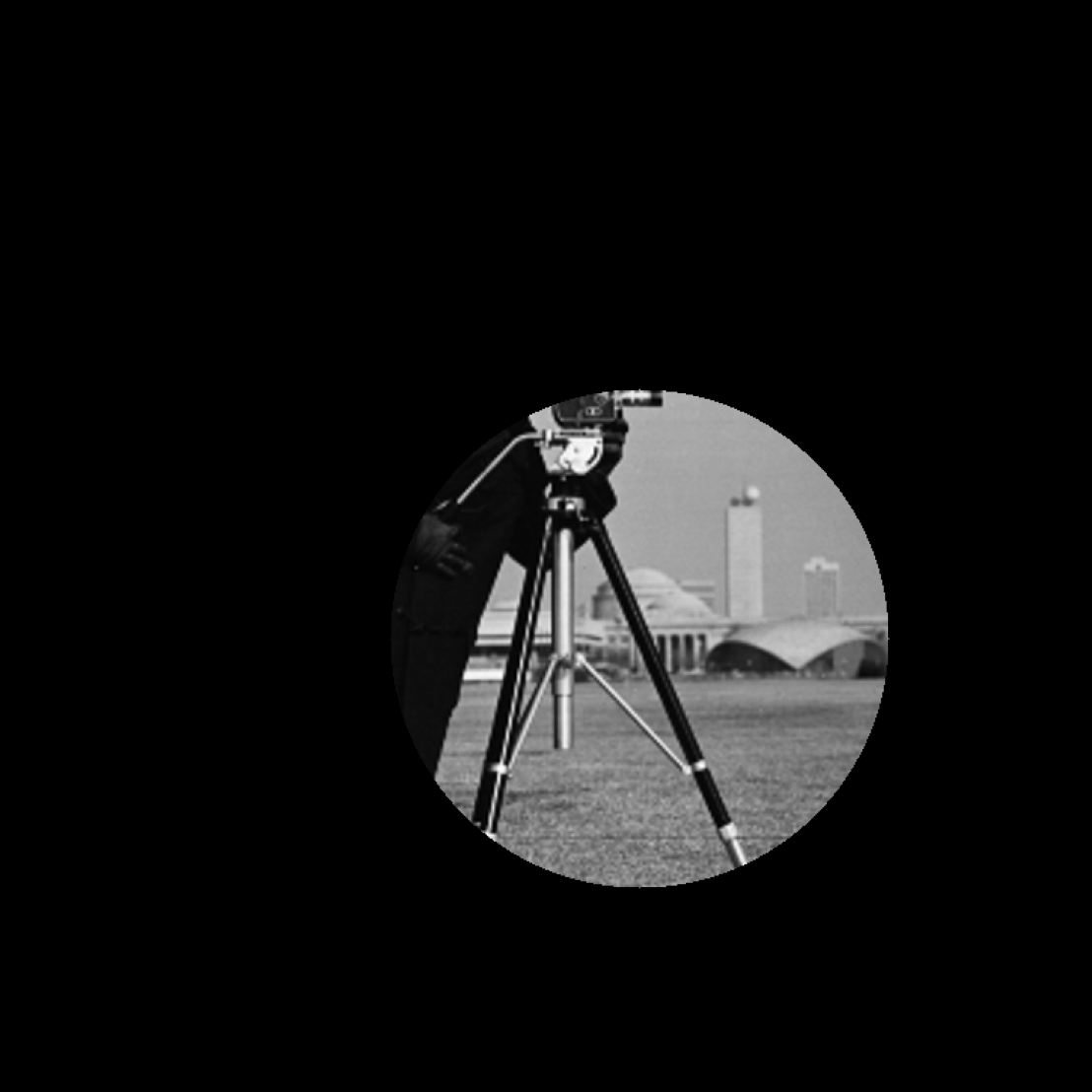}
        &
        \includegraphics[width=0.12\textwidth]{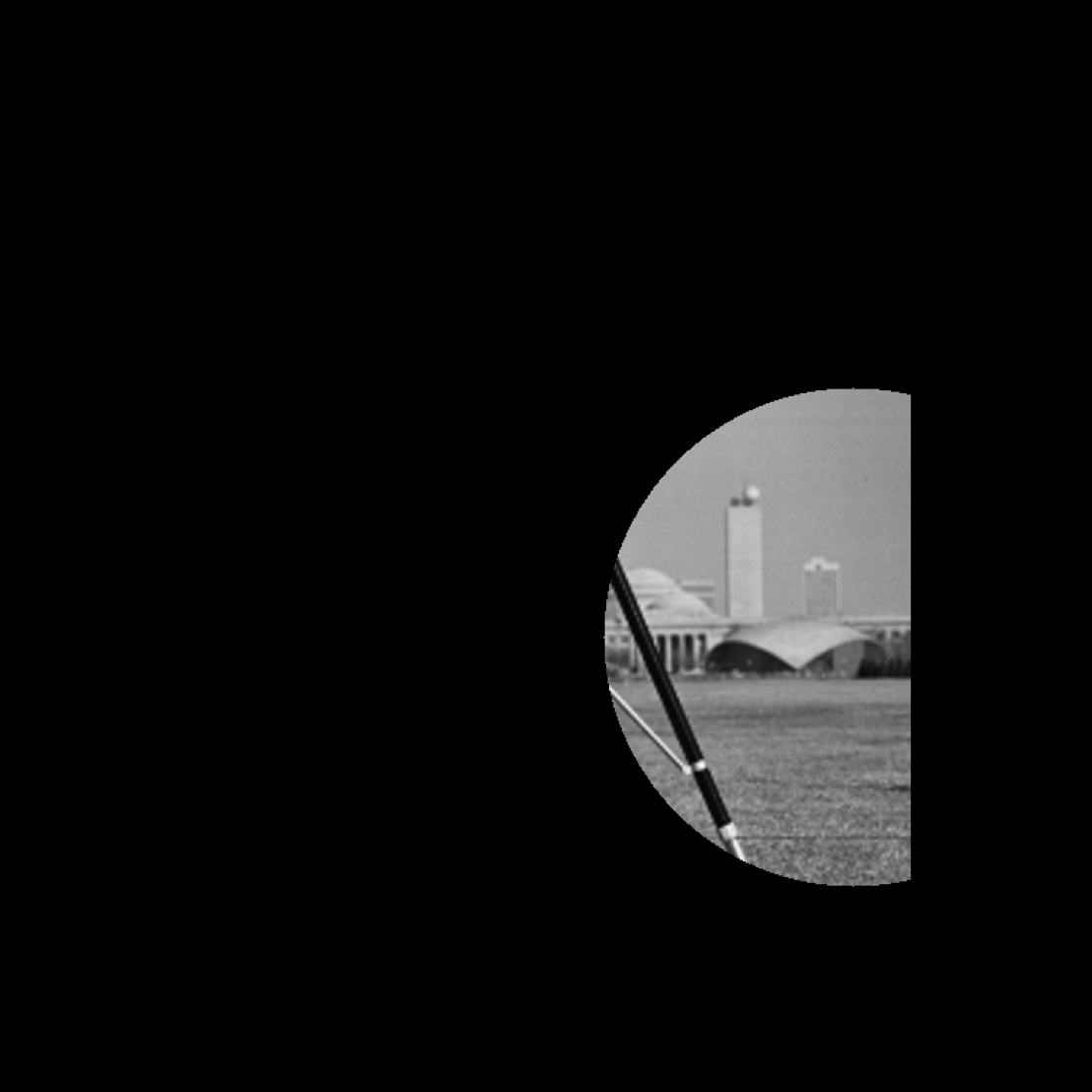}
        &
        \includegraphics[width=0.12\textwidth]{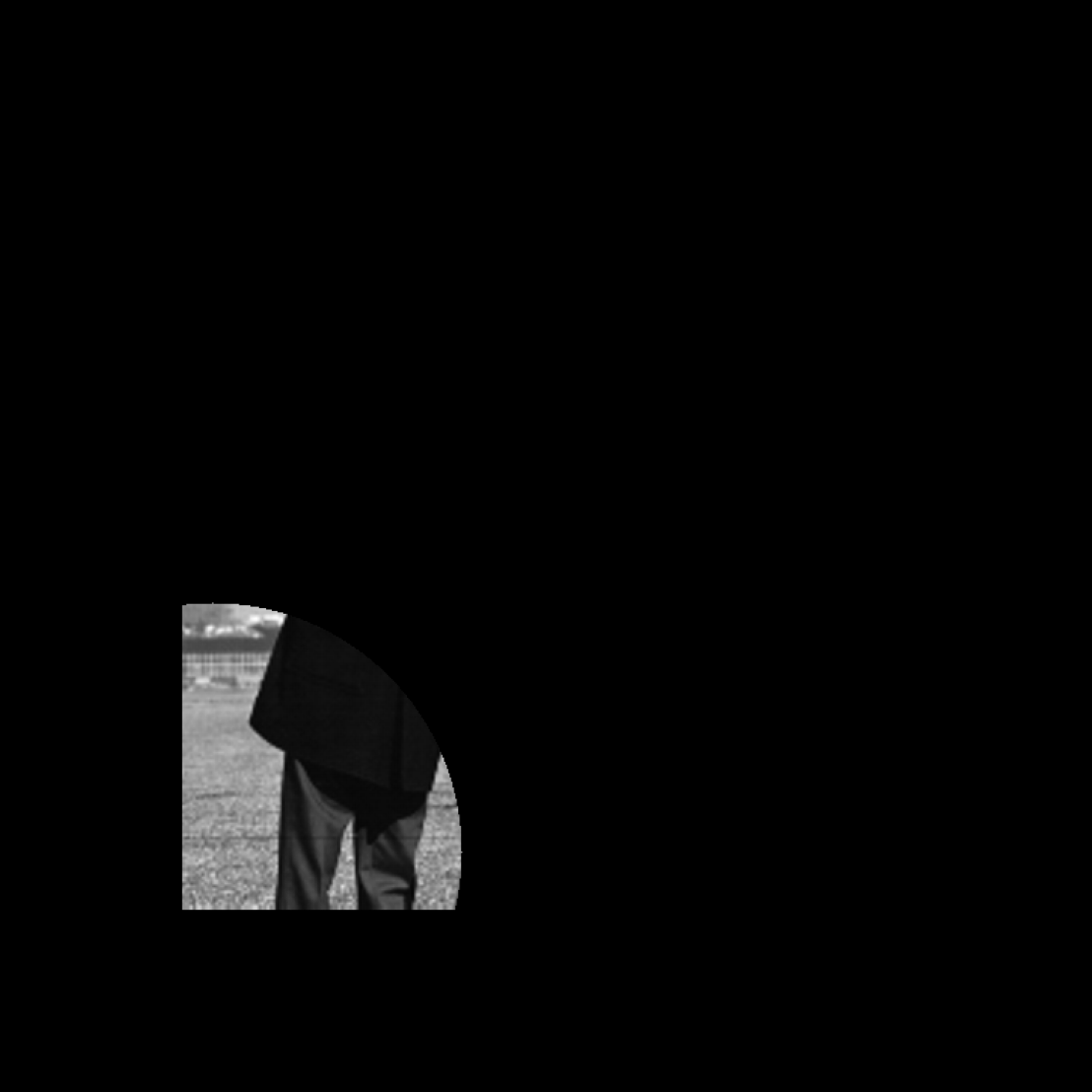}
        &
        \includegraphics[width=0.12\textwidth]{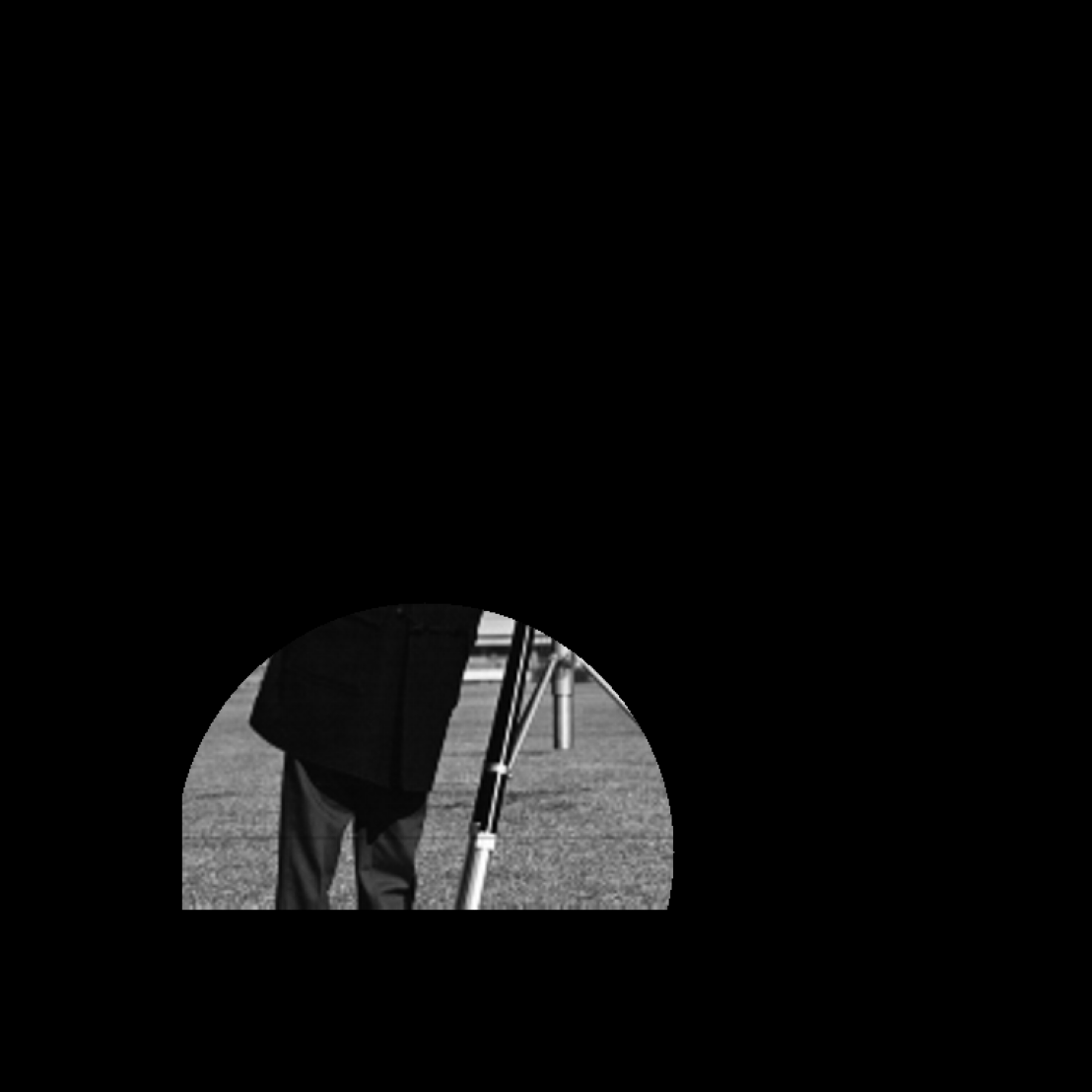}
        &
        \includegraphics[width=0.12\textwidth]{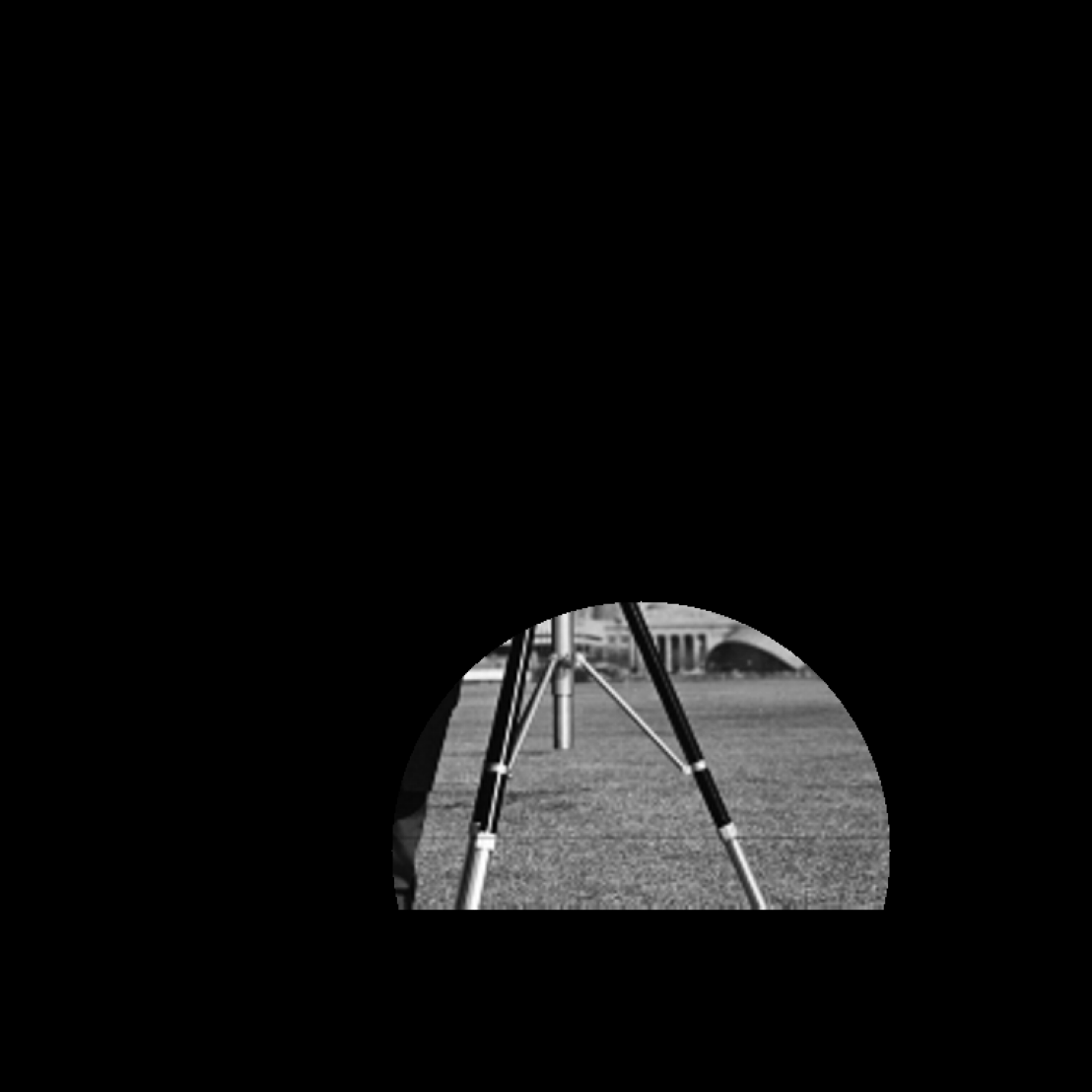}
        &
        \includegraphics[width=0.12\textwidth]{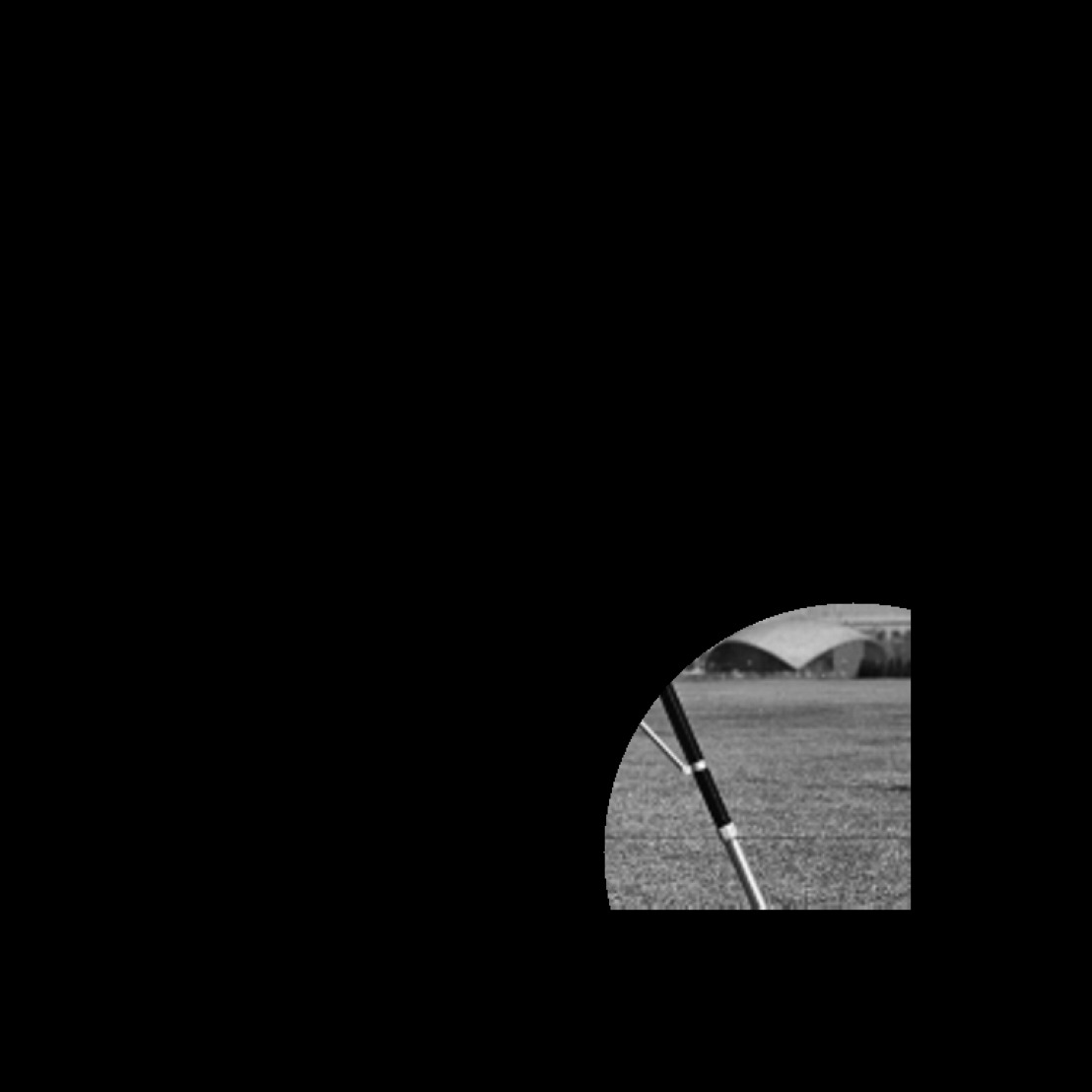}
    \end{tabular}
    \caption{\small{Scanning positions of the illumination window described in Section~\ref{subsubsec: blind_experimental_setup}}. Here, each probe illuminates a circle of radius $175$ pixels and shifts $150$ pixels at a time, leading to roughly $50\%$ overlap between consecutive probes.}
    \label{fig: circular_probes}
\end{figure}
\subsubsection{Experimental Setup}
\label{subsubsec: blind_experimental_setup}
For the blind experiments, we follow a very similar setup to that in~\cite{pham2019semi}. 
We use the same images of the baboon and cameraman where the baboon image is set to be the magnitude and the cameraman is set to be the phase of the ground truth. We also set the range of the true phase to be $[0, \frac{\pi}{2}]$ as before.

\begin{figure}[t]
    \centering
    \begin{tabular}{ccc}
        ePIE Phase & PFT Warmup Phase & Hybrid ePIE Phase
        \\
        \includegraphics[width=0.31\textwidth]{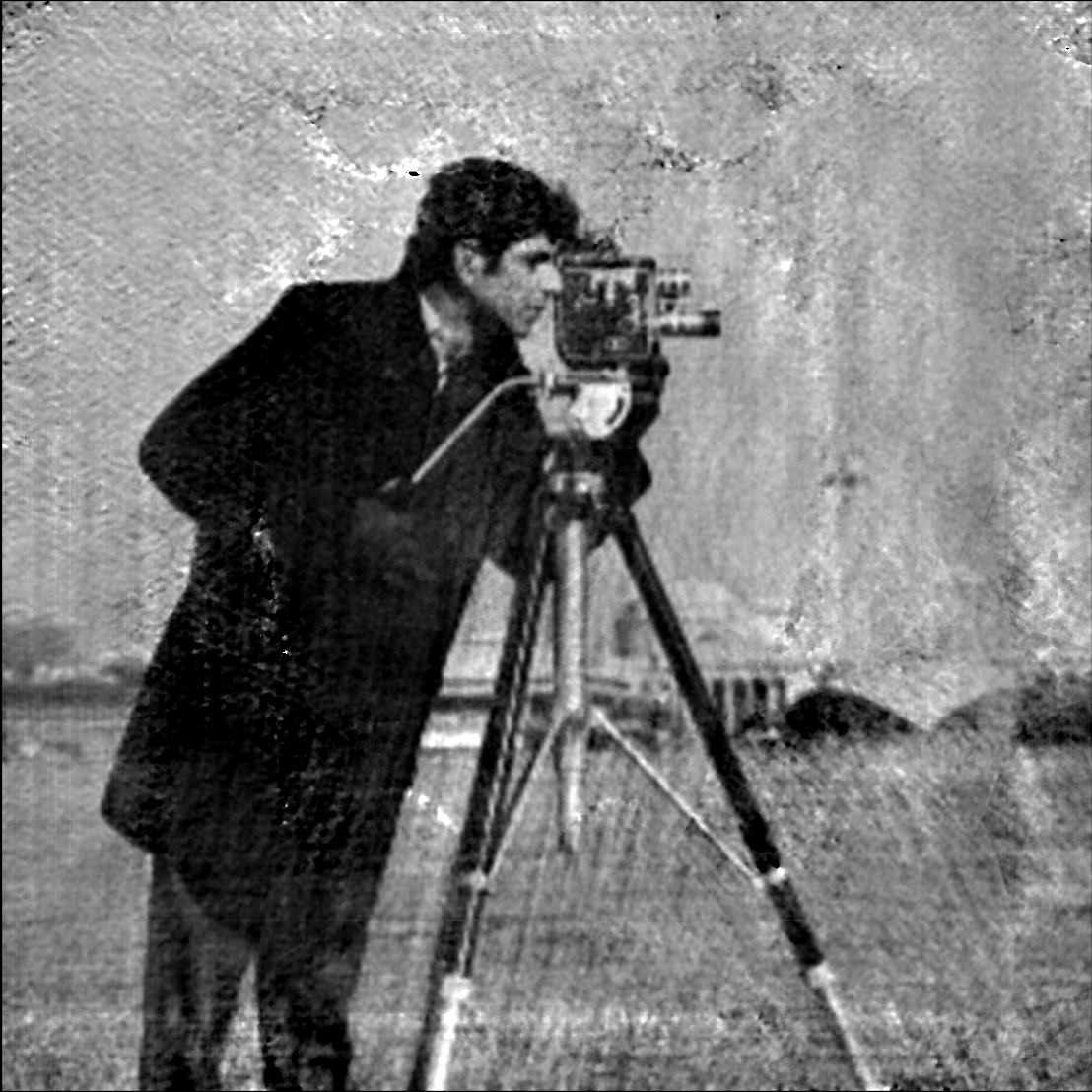}
        &
        \includegraphics[width=0.31\textwidth]{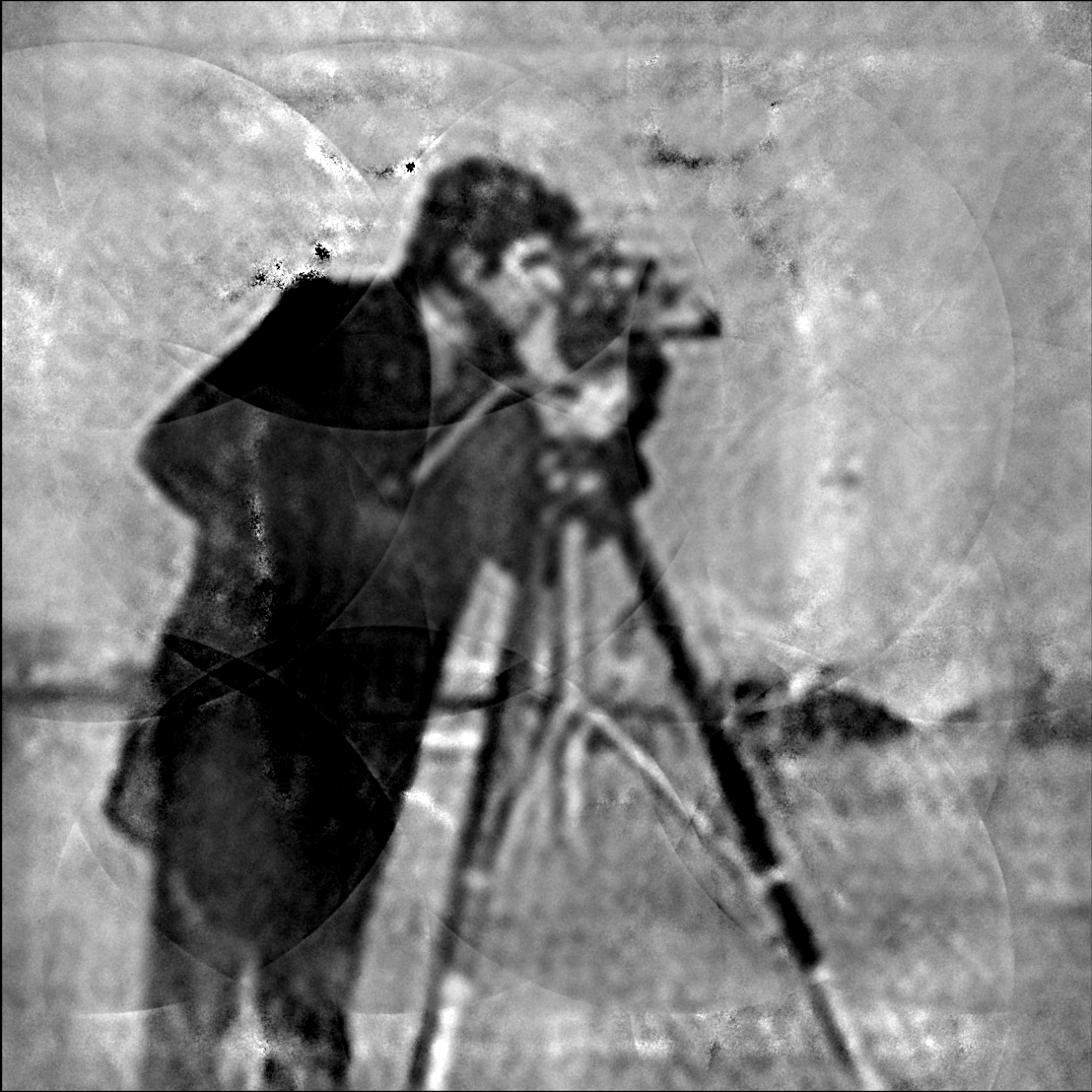}
        &
        \includegraphics[width=0.31\textwidth]{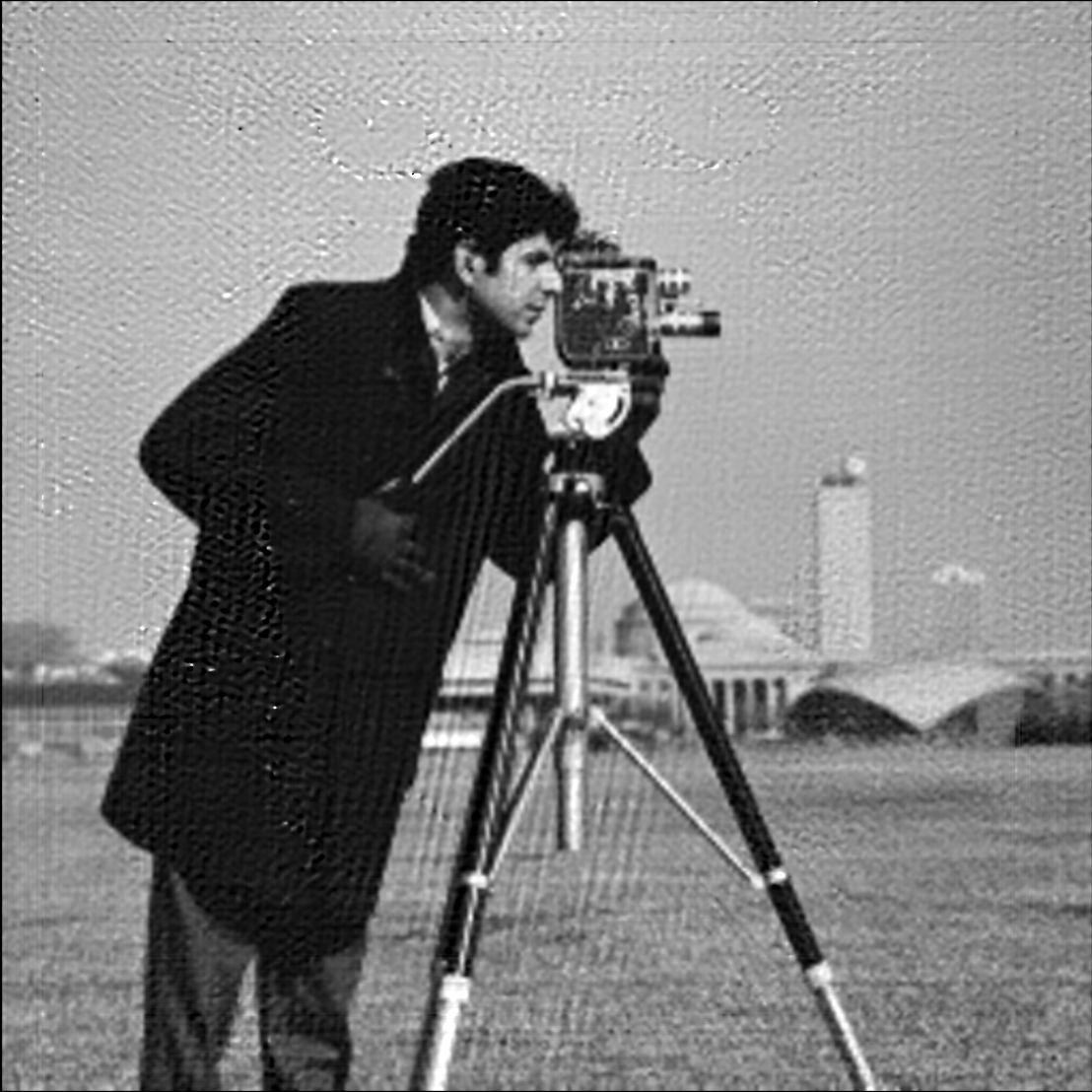}
        \\
        ePIE Magnitude & PFT Warmup Magnitude & Hybrid ePIE Magnitude
        \\
        \includegraphics[width=0.31\textwidth]{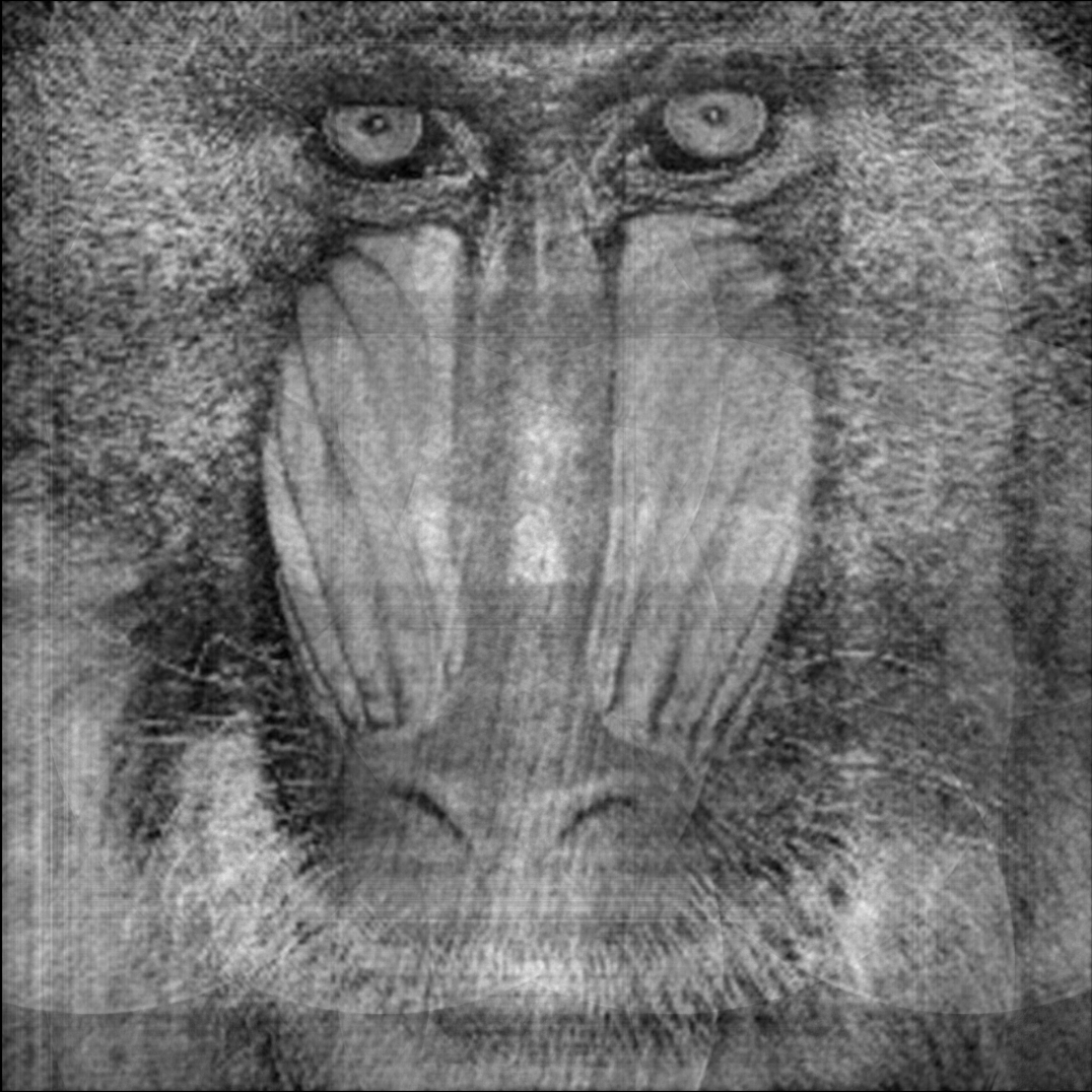}
        &
        \includegraphics[width=0.31\textwidth]{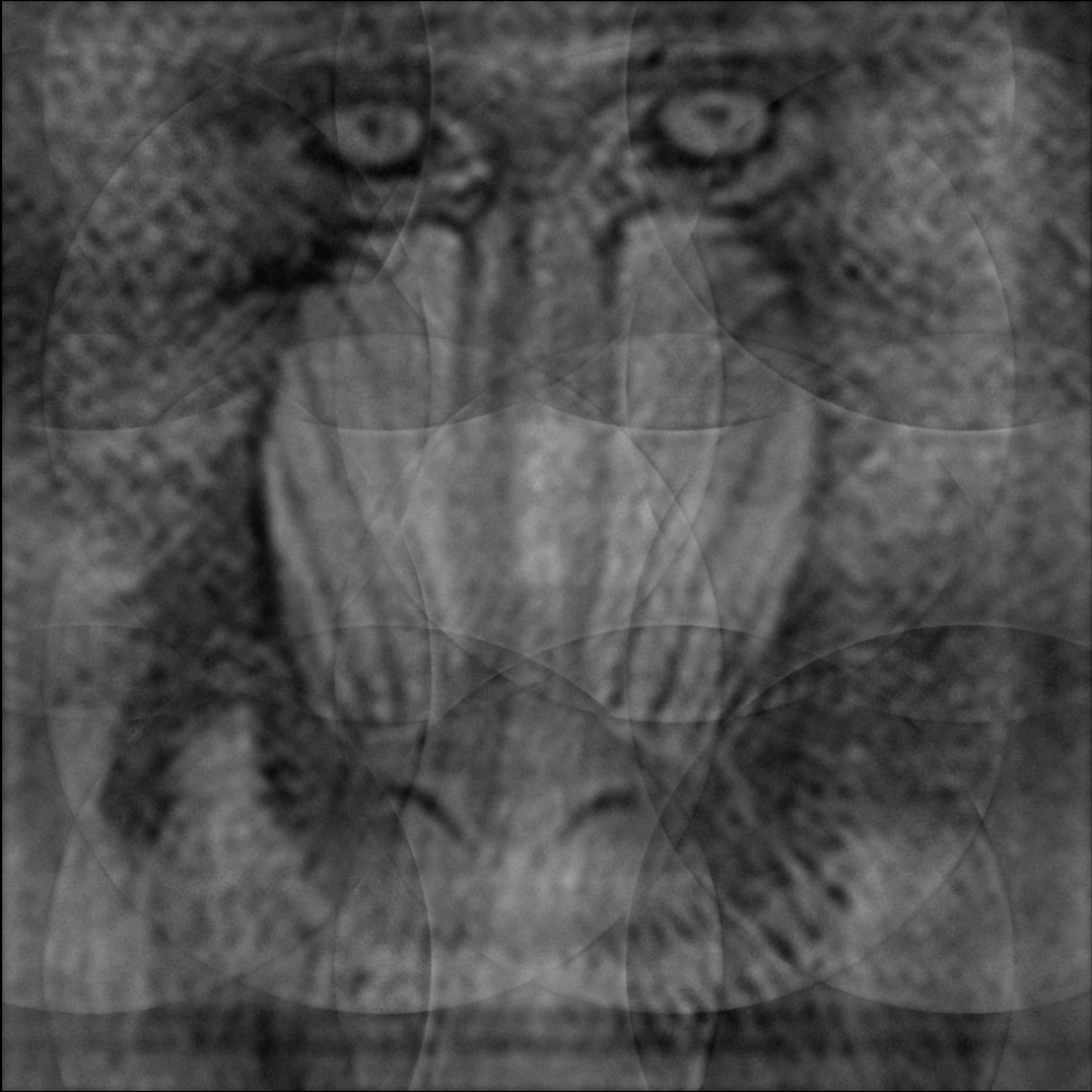}
        &
        \includegraphics[width=0.31\textwidth]{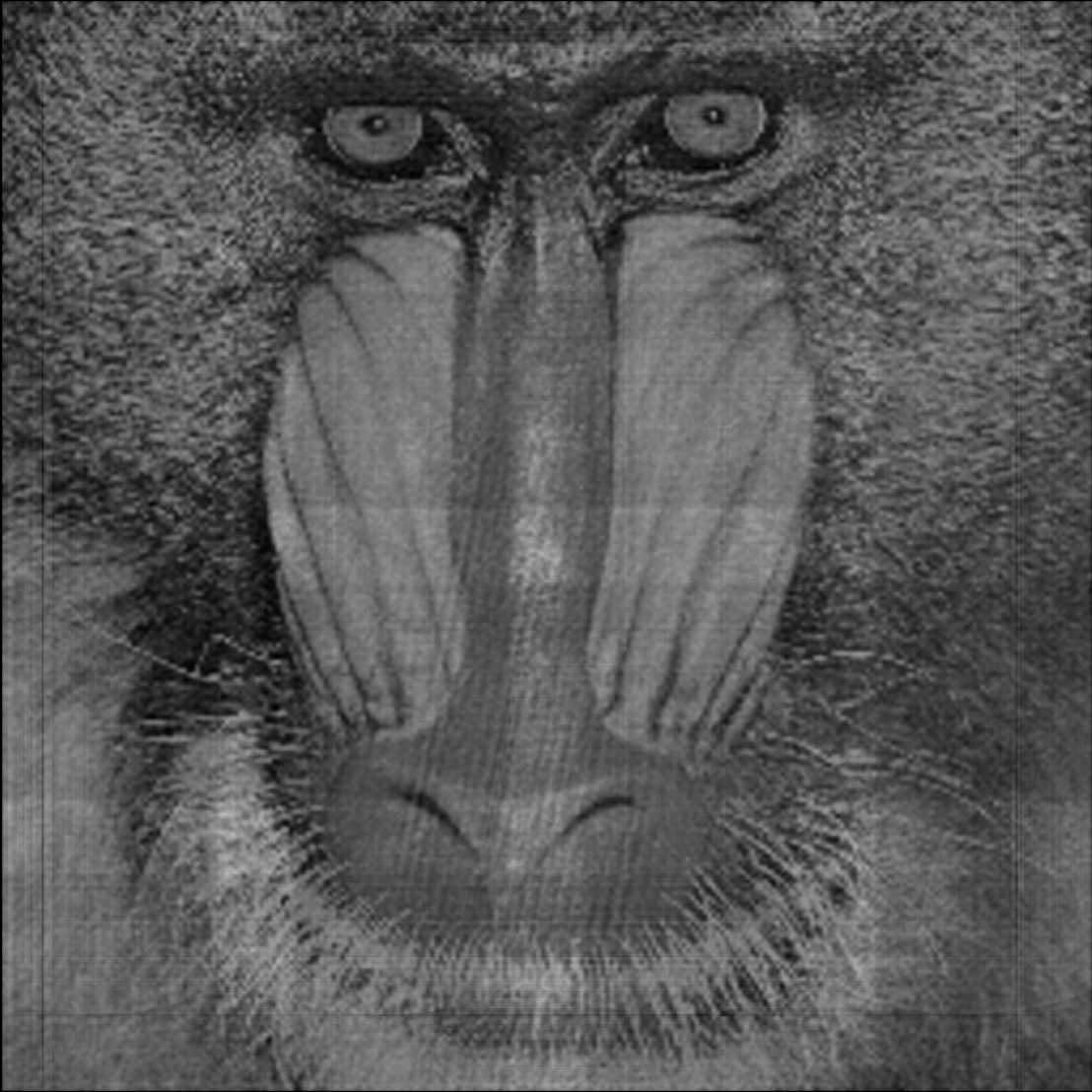}
    \end{tabular}
    \label{fig: blind_reconstructions}
    \caption{\small{Large-scale phase and magnitude reconstructions of the FFT-based ePIE (first column), PFT-based ePIE warmup (second column), and Hybrid ePIE (third column. These reconstructions are for the large-scale ptychography problem where $n_1 = n_2 = 8200$.}}
    \label{fig: blind_reconstruction_images}
\end{figure}

Following~\cite{pham2019semi}, we use circular probes for the blind ptychography problem and the probe, $\omega$ is also to be reconstructed. Note this is a very different experimental setup to that the nonblind setting.  
For a specified radius of pixels $r$, each illumination window $Q_{k} \in \RR^{n \times n}$ for $k = 1, \dots, N$ is a matrix that satisfies
\begin{equation}
    [Q_{k}]_{i, j} = \begin{cases}
    1 & \text{ if pixel } (i, j) \text{ lies within circle of radius $r$} \\
    0 & \text{ otherwise}.
\end{cases}
\end{equation}
Moreover, the true probe function in this experiment is given by 
\begin{equation}
    \omega(x,y) = \exp\left(-\frac{x^2 + y^2}{2\sigma^2}\right).
\end{equation}
where $\sigma = 10^{6}$. 
The probe is shifted $\frac{150n}{768}$ pixels at a time starting from the top-left corner to the right until it reaches the top-right corner. 
The probe is then shifted downwards and then to the left until it reaches the left edge again. This scanning procedure is continued until we have covered the entire image, and results in a total of $16$ probes with roughly $50\%$ overlap between adjacent probes. 
An illustration of the illumination window is shown in Figure~\ref{fig: circular_probes} for the small-scale experiment (Section~\ref{subsubsec: blind_distribution_relerrs}) and in Figure~\ref{fig: circular_probes_large} for the large-scale experiment (Section~\ref{subsubsec: blind_large_scale}). We remark, unlike in the non-blind case, there are nuances in the experimental setup with circular illumination windows that do not allow for the exact same settings when moving up to the large-scale setting. Moreover, following~\cite{pham2019semi} due to the circular nature of the probes, images need to be padded so that the entire image can be covered by the probes. This leads to reconstructions that include the padding (see Figure~\ref{fig: blind_reconstruction_images}).

Finally, we remark that in this implementation of ePIE~\eqref{eq: pie_update}, \eqref{eq: epie_probe_update}, we follow the techniques proposed in~\cite{pham2019semi} where the probes, i.e., the value of $j$, are chosen uniformly at random with no replacement; that is, rather than choose $j$ sequentially in each iteration (as was done in the non-blind case), we choose $j$ randomly until all N probes have been utilized. 
This random selection is chosen following, where it is observed that this randomness helps avoid local minima. 

\begin{figure}[t]
    \centering
    \begin{tabular}{ccc}
        a) $z$ & b) Magnitude of $z$ & c) Phase of $z$
        \\
        \includegraphics[width=0.3\textwidth]{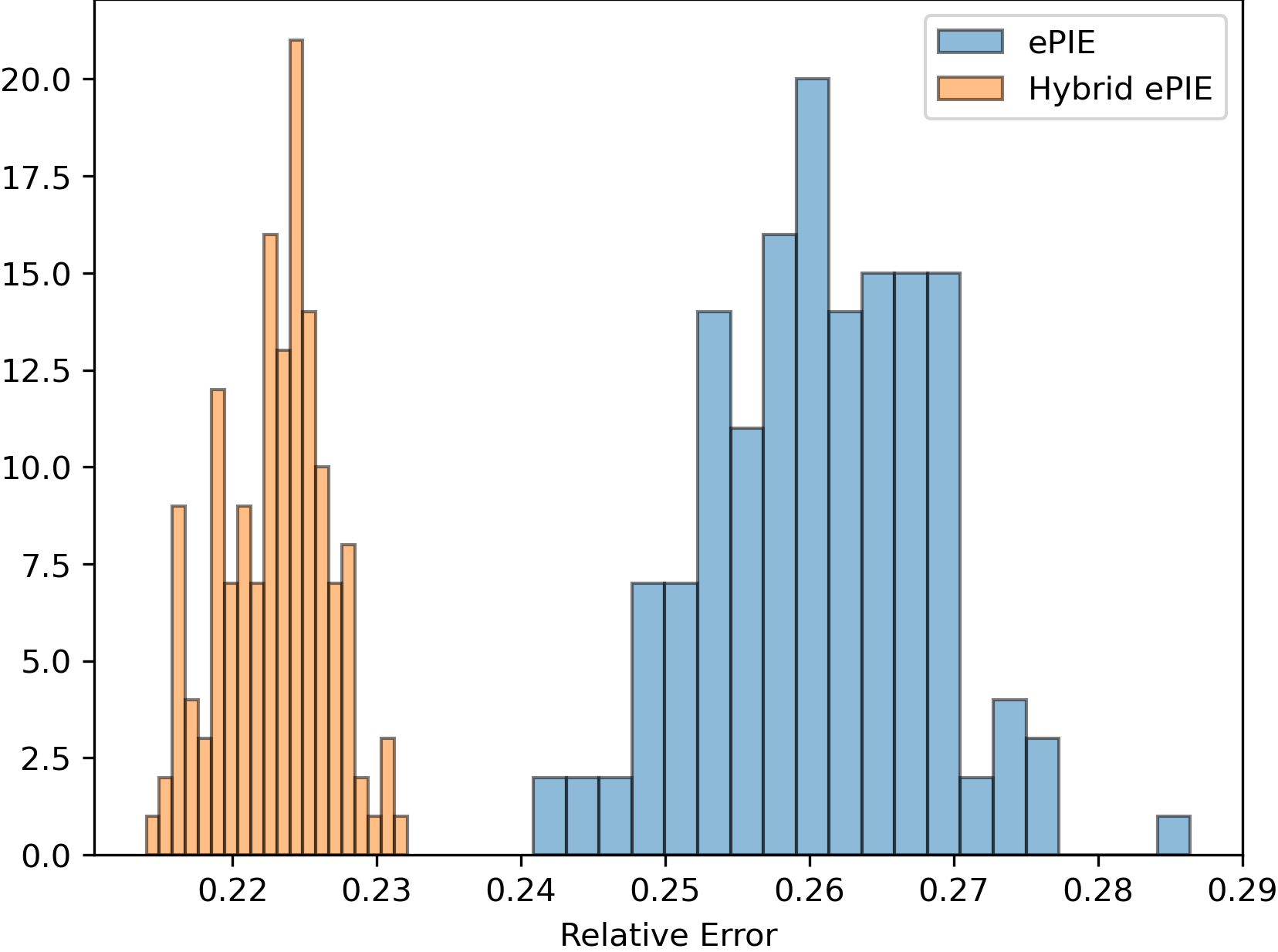}
        &
        \includegraphics[width=0.3\textwidth]{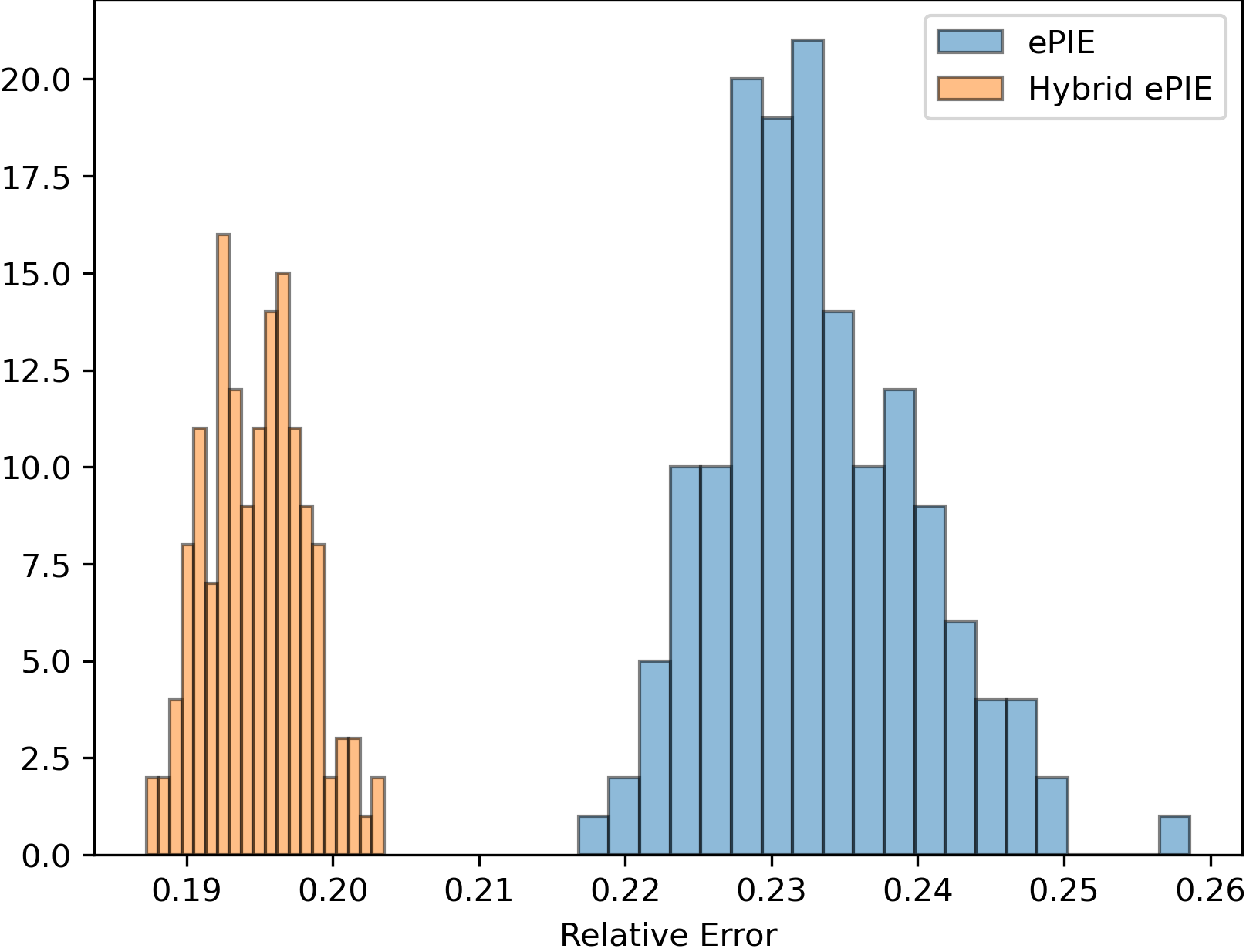}
        &
        \includegraphics[width=0.3\textwidth]{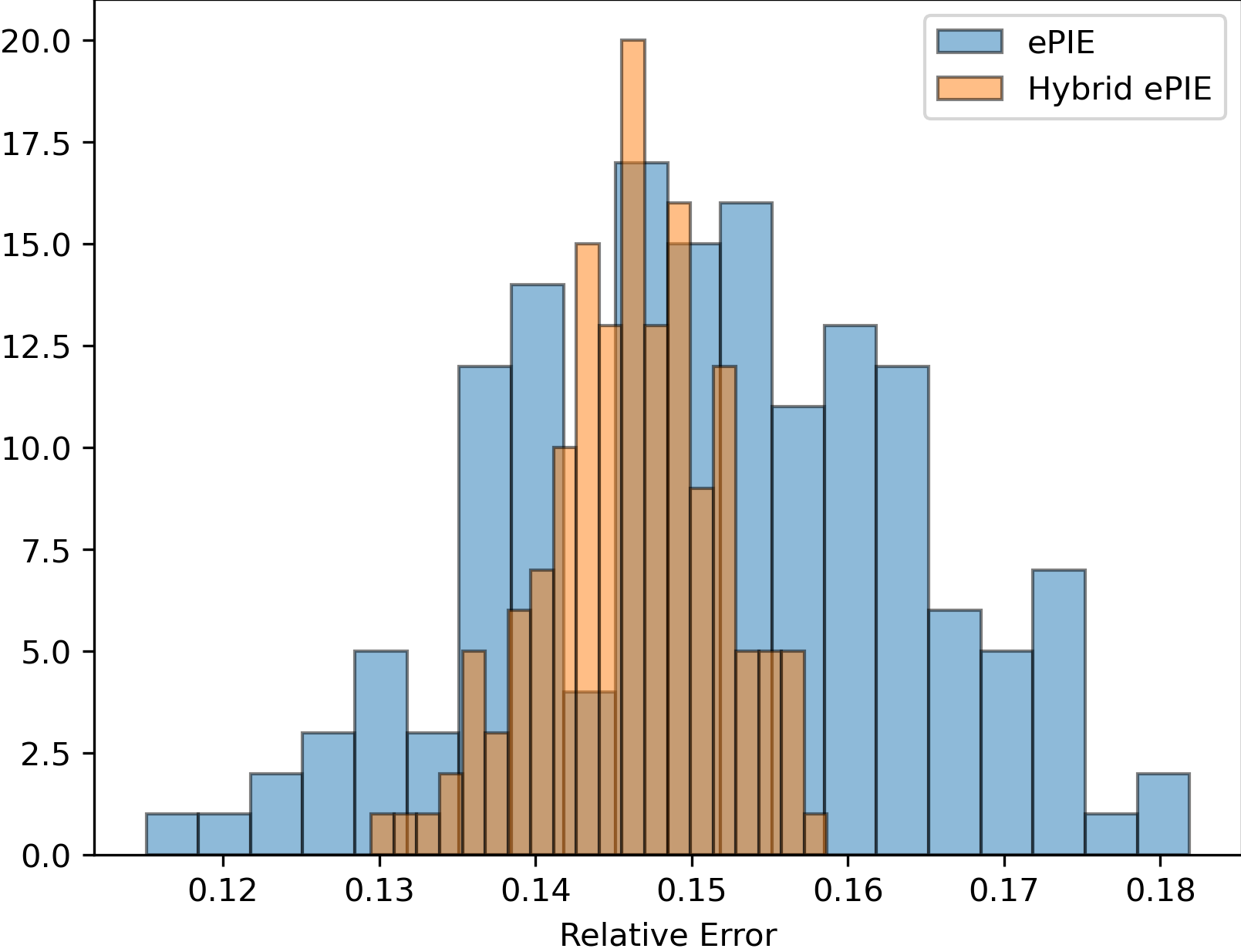}
    \end{tabular}
    \caption{\small{Histogram of the final reconstruction relative errors. The blue histogram shows the relative error frequency for ePIE and the orange histogram shows the relative error frequency for hybrid ePIE. a) shows the relative error of the reconstructed object, b) shows relative errors of only the magnitude of the object, and c) shows the relative errors of the phase of the object.}}
    \label{fig: epie_rel_err}
\end{figure}

\begin{figure}[t]
    \centering
    \begin{tabular}{cccc}
        Magnitude SSIM & Phase SSIM & Magnitude PSNR & Phase PSNR
        \\
        \includegraphics[width=0.23\textwidth]{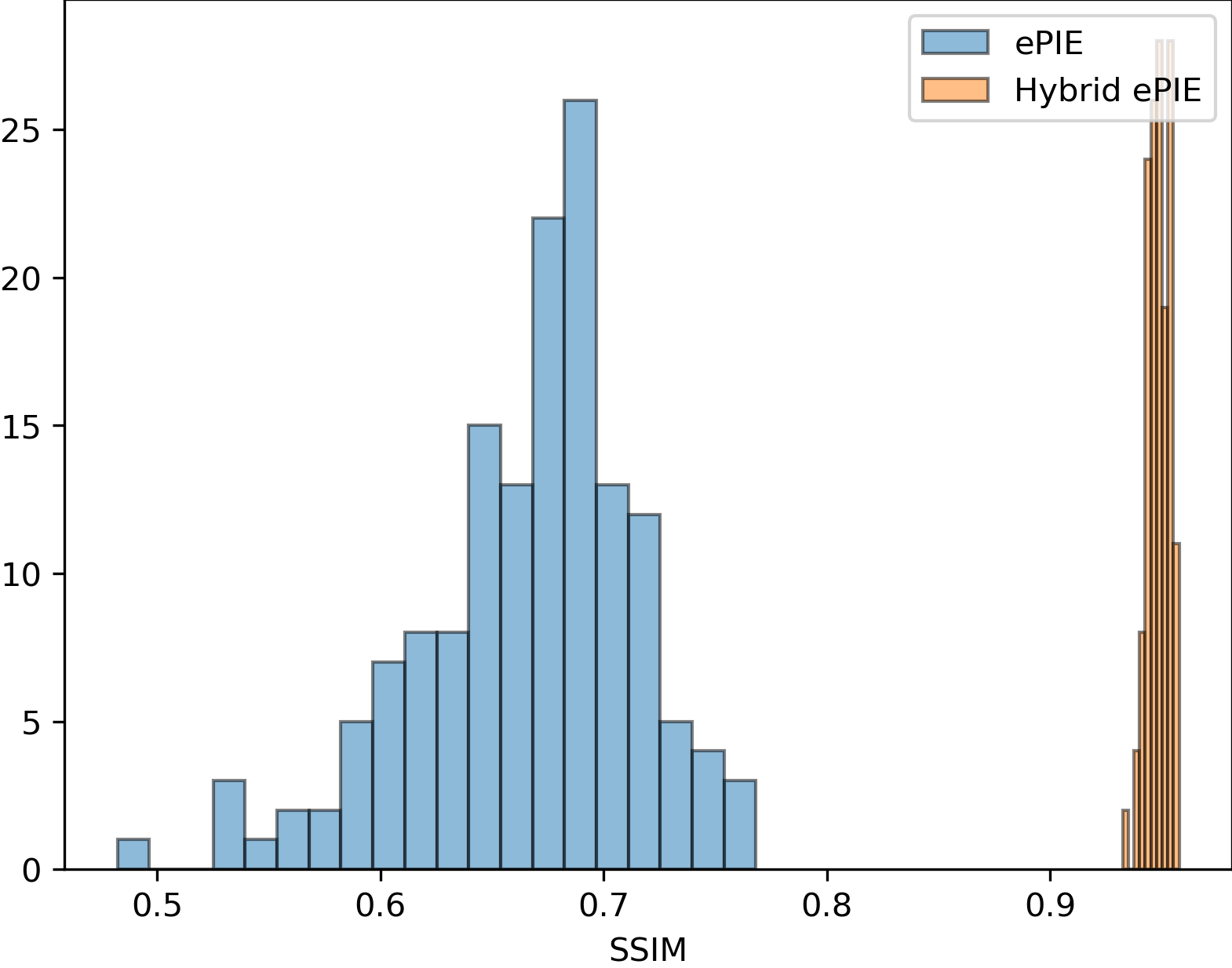}
        &
        \includegraphics[width=0.23\textwidth]{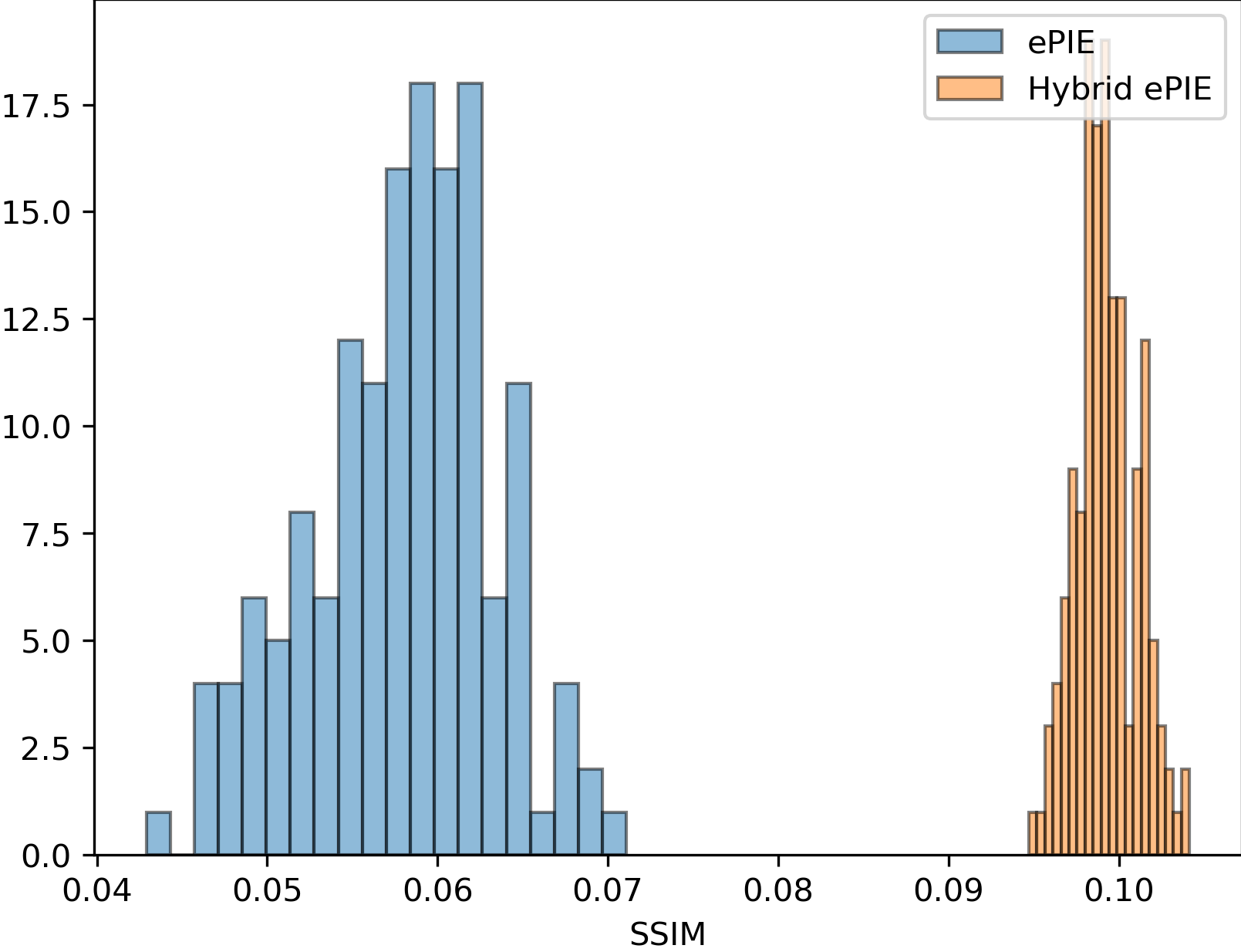}
        &
        \includegraphics[width=0.23\textwidth]{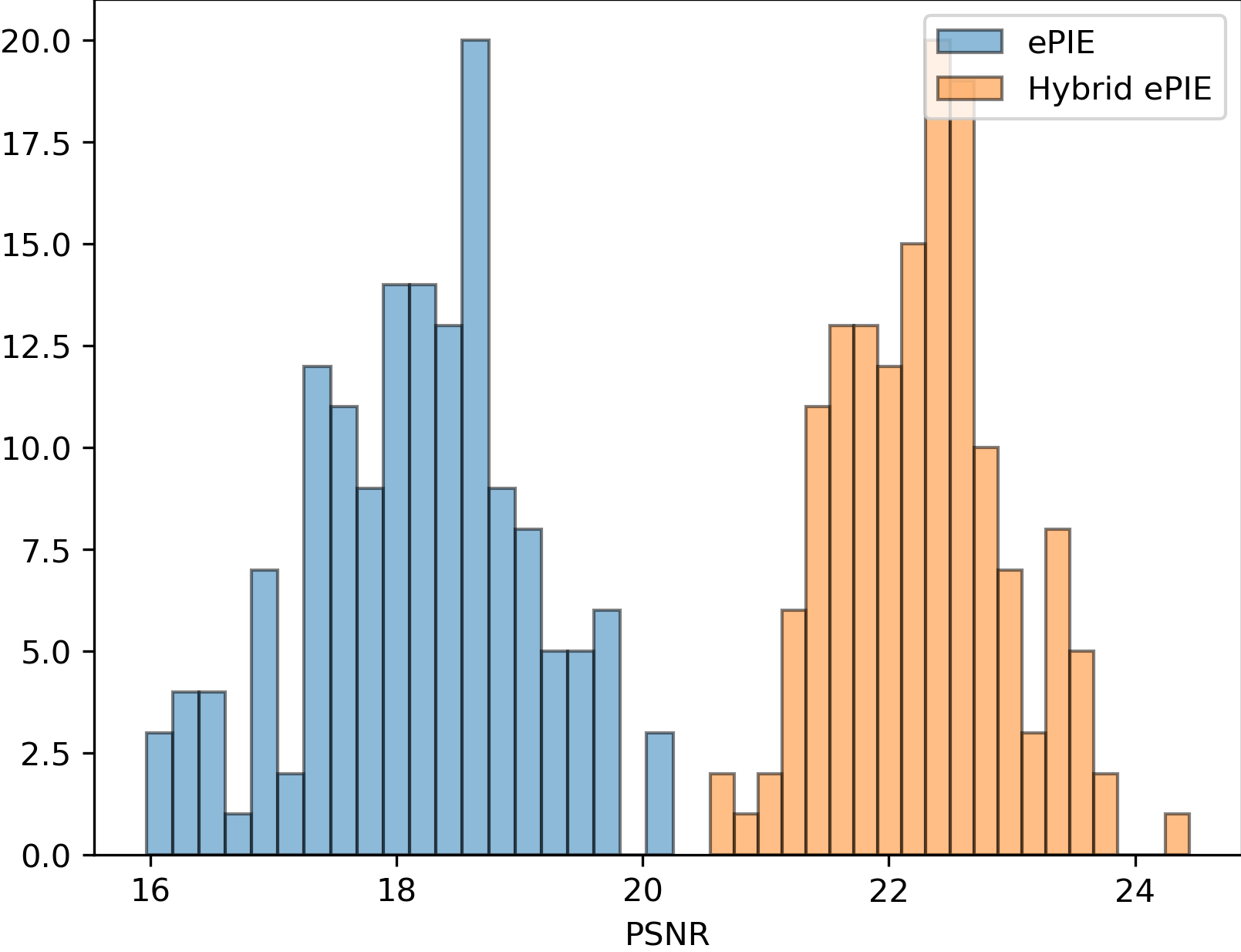}
        &
        \includegraphics[width=0.23\textwidth]{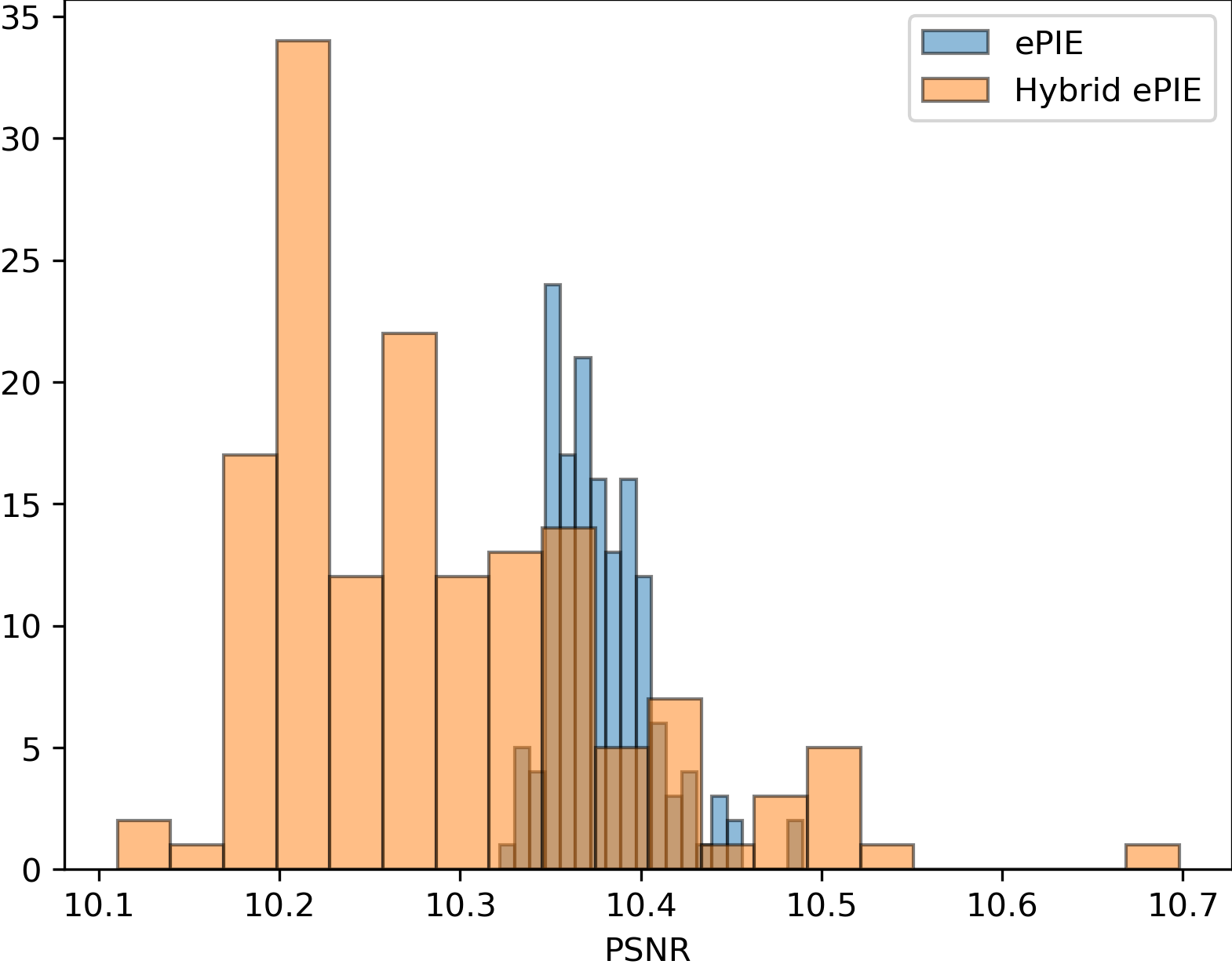}
    \end{tabular}
    \caption{\small{Histogram of the final SSIM and PSNR values for the reconstructed images. The blue histogram shows the SSIM/PSNR for ePIE and the orange histogram shows the SSIM/PSNR for hybrid ePIE.}}
    \label{fig: epie_ssim_psnr}
\end{figure}

\subsubsection{Distribution of Reconstruction Quality}
\label{subsubsec: blind_distribution_relerrs}
We now consider the distribution of the relative errors for the reconstructed images from the ePIE and warm-started ePIE algorithms on small-scale images of size $512 \times 512$. As with the non-blind ptychography experiments, we generate a sample of $150$ guess images where the guess magnitude is initialized with uniform random values in the range $[0, 1]$ and $[0, \frac{\pi}{2}]$ for the guess phase. We then compute the relative error of the final reconstructions using~\eqref{eq: rel_error}, as well as the PSNR, and SSIM of the magnitude and phase, respectively, to further analyze the quality of the reconstructions.

Figure~\ref{fig: epie_rel_err} shows the relative errors between a) the full object, b) the magnitude of the object, and c) the phase of the object. The results show that hybrid ePIE tends to find a better distribution of reconstructions compared to the distribution of vanilla ePIE. In particular, we see that the relative errors of the magnitude are lower in the hybrid approach. Moreover, we observe that the relative errors of the phase are similar in value (but with lower variance) between hybrid ePIE and traditional ePIE.
We also compare the distribution of the final SSIM and PSNR values of the reconstructed image from ePIE and hybrid ePIE shown in Figure~\ref{fig: epie_ssim_psnr}. We observe that the distribution of the SSIM values are overall higher in the hybrid approach. However, the hybrid ePIE magnitude reconstructions achieve higher PSNR values and overall slightly lower phase reconstruction PSNR values than vanilla ePIE. In other words, the PSNR histograms support the result of the relative errors, demonstrating the reconstruction of the magnitude makes up for the reconstruction of the phase in the hybrid algorithm. 

\begin{figure}[t]
    \small
    \centering
    \begin{tabular}{cccc}
        \includegraphics[width=0.23\textwidth]{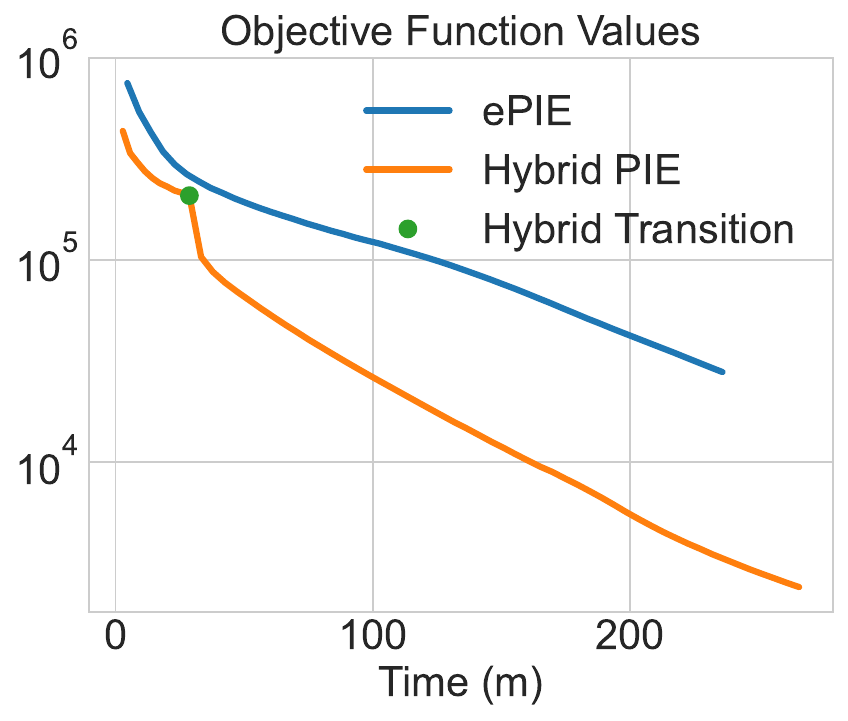}
        &
        \includegraphics[width=0.23\textwidth]{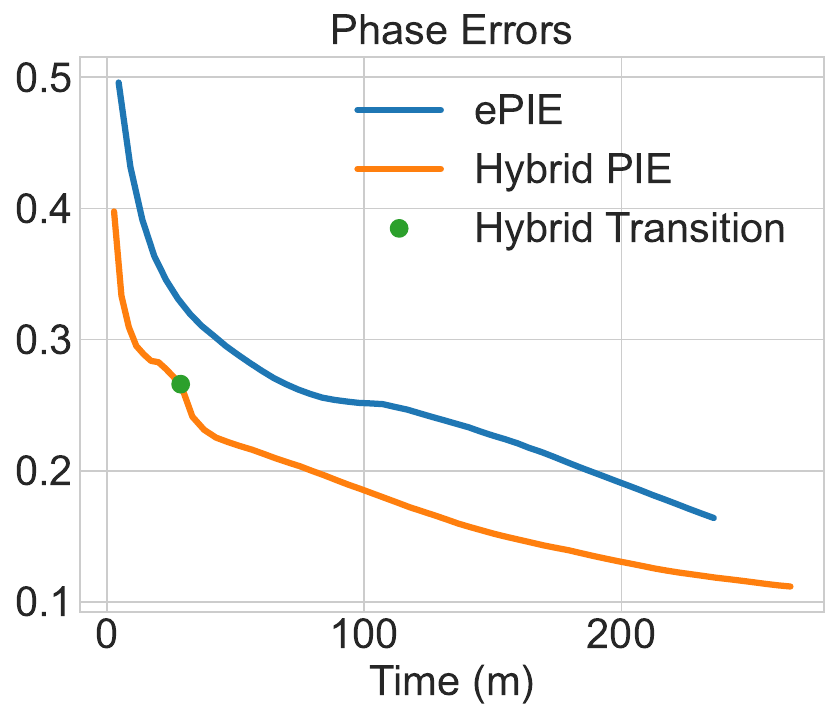}
        & 
        \includegraphics[width=0.23\textwidth]{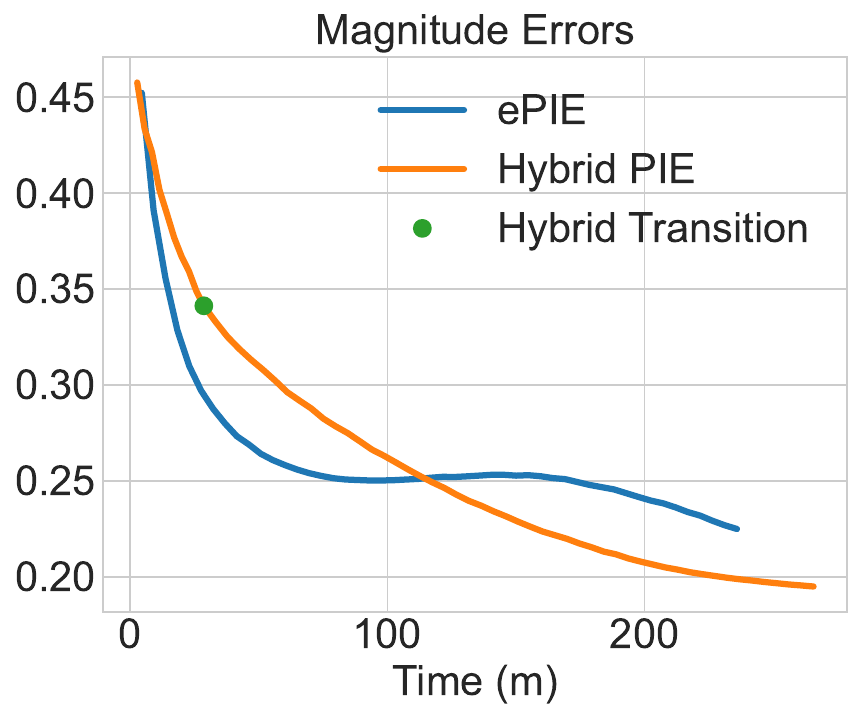}
        &
        \includegraphics[width=0.23\textwidth]{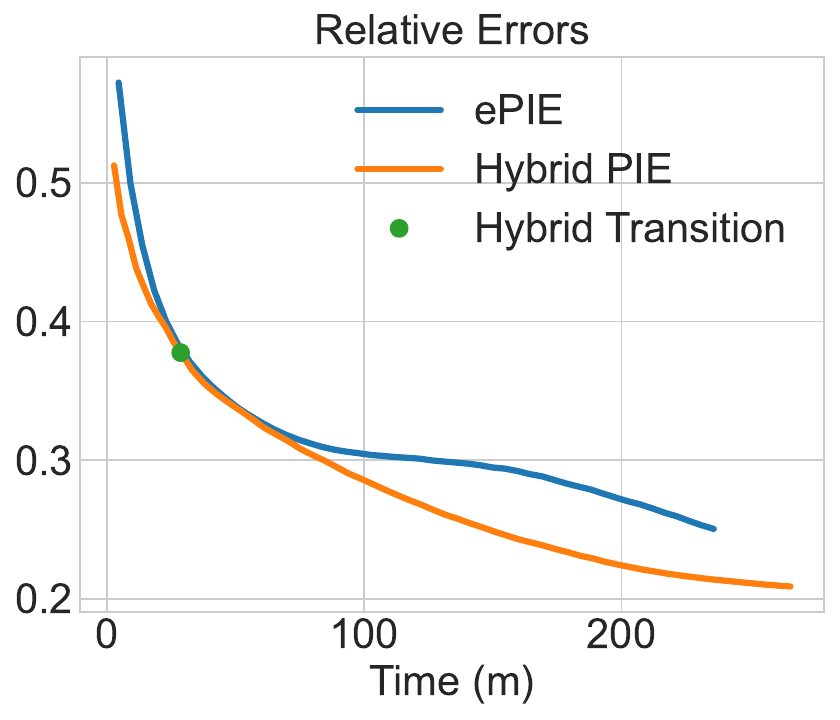}
    \end{tabular}
    \caption{\small{Objective values on full reconstruction (first column), phase relative errors (second column), magnitude relative errors (third column), and relative errors (fourth column) for a large-scale ptychographic reconstruction explained in Section~\ref{subsubsec: blind_experimental_setup}. The x-axis represents time in minutes.
    Green dot represents the transition from the PFT-based ePIE to the FFT-based ePIE in the hybrid algorithm.
    }
    }
    \label{fig: blind_rel_errors}
\end{figure}

\begin{figure}[t]
    \small
    \centering
    \begin{tabular}{cccc}
        \includegraphics[width=0.23\textwidth]{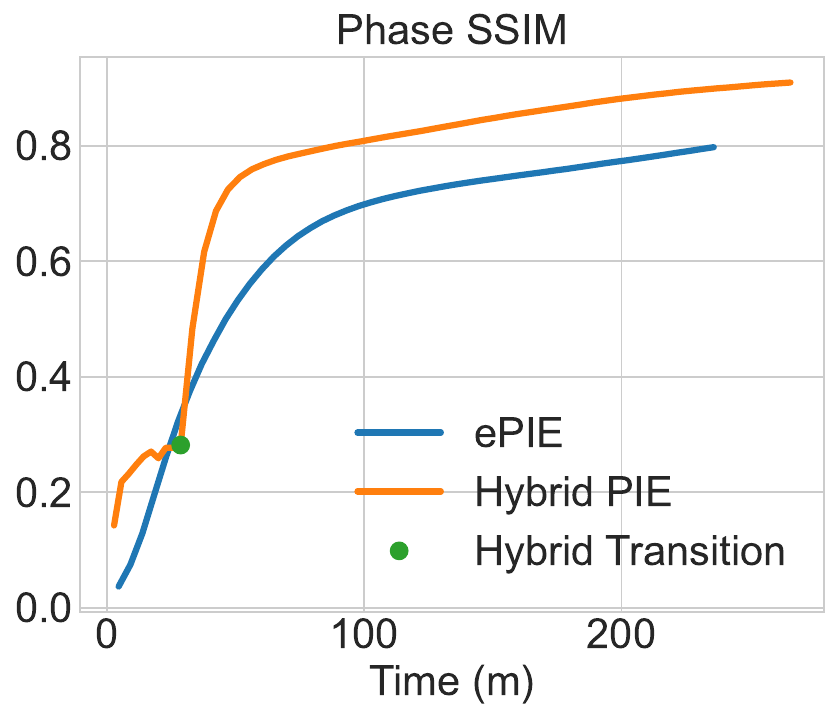}
        &
        \includegraphics[width=0.23\textwidth]{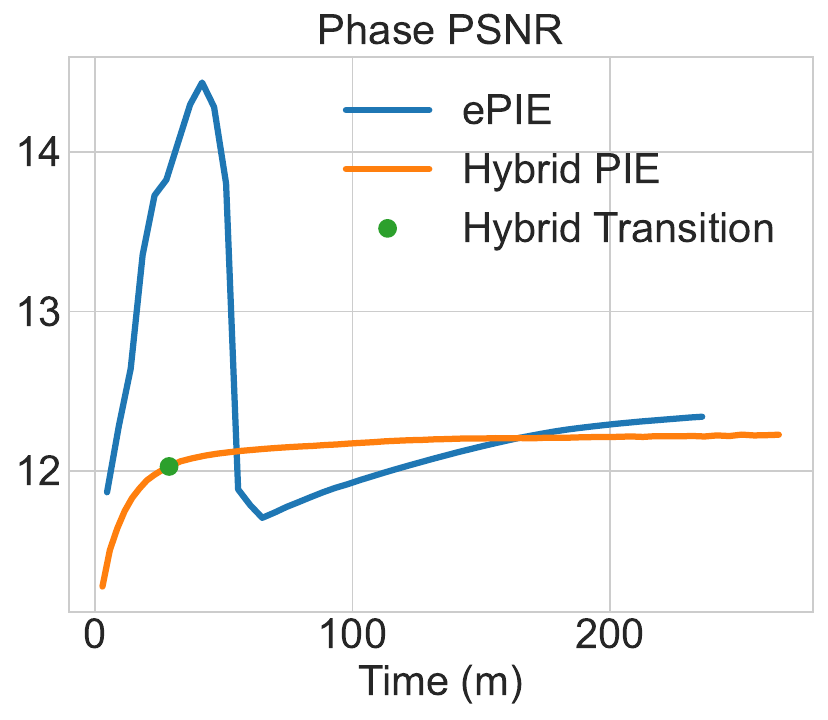}
        &
        \includegraphics[width=0.23\textwidth]{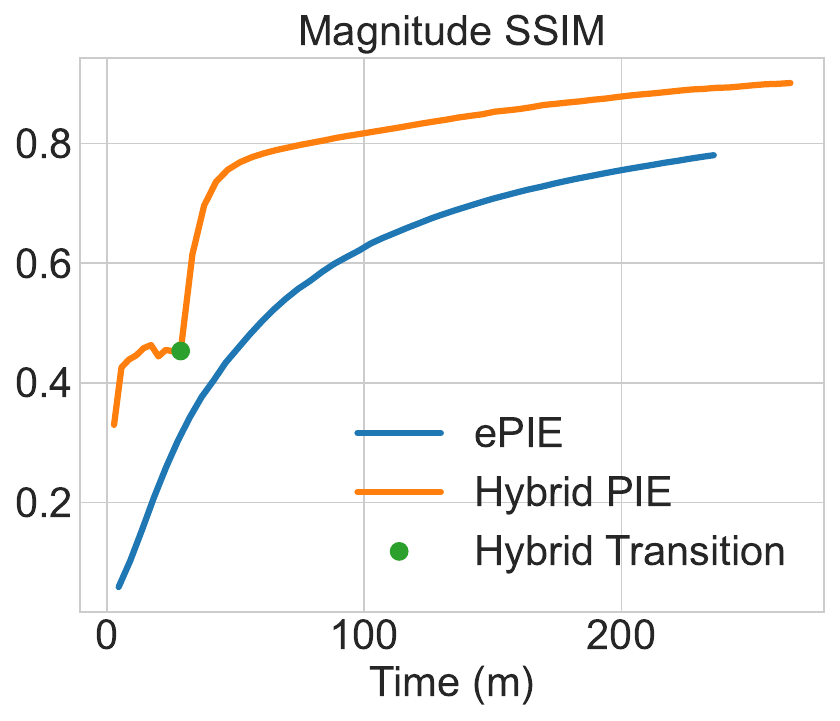}
        &
        \includegraphics[width=0.23\textwidth]{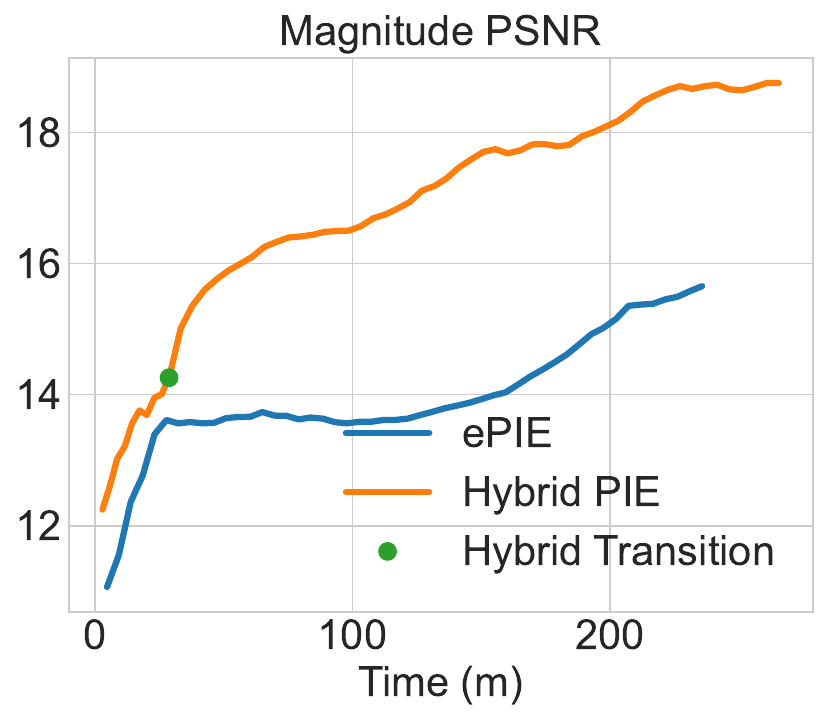}
    \end{tabular}
    \caption{SSIM and PSNR for ePIE (blue) and Hybrid ePIE (orange) for the phase and magnitude over time (in minutes). The green dot represents the transition from the PFT-based ePIE to the FFT-based ePIE in the hybrid algorithm.}
    \label{fig: blind_ssim_psnr}
\end{figure}

\subsubsection{Large-Scale Reconstruction}
\label{subsubsec: blind_large_scale}
As in the nonblind case, we demonstrate the computational benefits of using a PFT-based ePIE in the proposed algorithm. Here, $n_1 = n_2 = m_1 = m_2 = 8200$, and $\widetilde{m}_1 = \widetilde{m}_2 = p_1 = p_2 = 64$, leading to a cropped image of size $128 \times 128$. The degree of the approximating polynomial in this setting is given by degrees $r_1 = r_2 = 13$, which is obtained from the precomputed values $\xi(\varepsilon, r)$ defined in Section~\ref{subsec: pft_offline_derivation}.
Here, we run the PFT-based ePIE algorithm until a tolerance of $\epsilon_{\text{PFT}} = 10^{-2}$ or a maximum number of $10$ iterations is reached. Afterwards, we run the standard FFT-based ePIE until a tolerance of $\epsilon = 5 \times 10^{-4}$ or a maximum number of $50$ iterations is reached. Here, we choose step sizes $\gamma = \beta = 10$ for the FFT-based ePIE and $\gamma_{PFT} = \beta_{PFT} = 2 \times 10^{-3}$ for the PFT-based ePIE. These were chosen tuned using a logarithmic gridsearch over the set $\{10^{-6}, 10^{-5}, \ldots, 10^2, 10^{3}\}$, and afterwards using a standard gridsearch starting from the currently chosen parameter, e.g., $10^{-2}$, to the next order of magnitude, e.g., $10^{-1}$.  
The intuition for choosing a higher tolerance for the PFT-based ePIE is based on the fact that the PFT-based ePIE will primarily capture large features, whereas the FFT-based ePIE needs more iterations to capture fine details. For all reconstructions, we use total variation regularization~\cite{rudin1992nonlinear} with parameters $\lambda = 10^{-6}$ for the FFT-based ePIE and $\lambda = 10^3$ for the PFT-based ePIE. These are the parameters that led to the lowest relative error and were found over a logarithmic gridsearch over the set~$\{10^{-6}, 10^{-5}, \ldots, 10^2, 10^{3}\}$.

In Figure~\ref{fig: blind_rel_errors}, we show the objective values of the full reconstruction (first column), the phase (second column), the magnitude (third column), and the relative errors (fourth column). Here, the x-axis represents time in minutes.
Similar results are shown in Figure~\ref{fig: blind_ssim_psnr}, where we only show the phase and magnitude structural similarity index (SSIM)~\cite{wang2004image} and peak signal to noise ratio (PSNR)~\cite{poobathy2014edge}. 
We observe that while using the PFT-based ePIE as warm up leads to overall improved reconstructions, there are some trade-offs, this can be seen in, e.g., the hybrid ePIE magnitude reconstruction. However, for this particular large-scale run, we obtain a very good quality phase reconstruction as seen in Figure~\ref{fig: blind_reconstruction_images}. This is appealing as the primary object of interest in ptychography (also known as ptychographic phase retrieval). 

As a final note, there is a nontrivial padding element for this ptychographic experiment, where the image recovered consists of the reconstruction and corresponding noisy padding. To address this, we use the built-in function \texttt{match\_template} algorithm from scikit-image~\cite{lewis1995fast} to generate Figure~\ref{fig: blind_reconstruction_images} .

\section{Discussion}
We present a new hybrid algorithm for ptychographic phase retrieval based on the fast partial Fourier transform (PFT), which only computes the coefficients corresponding to low frequencies. The PFT is utilized within the ptychographic iterative engine (PIE). The hybrid algorithm consists of using PFT-based PIE in early iterations as a warm up followed by the standard FFT-based PIE algorithm. The core idea is to let the PFT-based PIE iterations capture large features corresponding to low frequencies whereas the FFT-based PIE iterations capture fine details in the reconstruction. Our numerical results demonstrate that the proposed hybrid PIE algorithm accelerates convergence by reducing the time-to-solution. This work also provides a PyTorch implementation of the PFT with automatic differentiation capabilities, enabling ease of use of the PFT within, e.g., deep learning architectures and other optimization algorithms. 
Future works include exploring a differentiable optimization methodology~\cite{agrawal2019differentiable, mckenzie2024differentiating, heaton2023explainable} of the proposed hybrid algorithm as well as distributed methods such as alternating direction method of multipliers~\cite{bui2024stochastic, boyd2011distributed, fung2019uncertainty}. 
Our code can be found in \url{https://github.com/mines-opt-ml/pft-for-ptycho}.

\section*{Acknowledgments}
\noindent Samy Wu Fung was partially funded by National Science Foundation award DMS-2110745. Stanley Osher was partially funded by Air Force Office of Scientific Research (AFOSR) MURI
FA9550-18-502, Office of Naval Research (ONR)
N00014-20-1-2787, and STROBE: a National Science Foundation Science and Technology Center under Grant No. DMR-1548924.
We thank Yong-chan Park for his help on setting up the PyTorch-based PFT code.

\bibliographystyle{abbrv}
\bibliography{references}

\newpage
\appendix
\section{2D PFT Algorithm}
\label{app: pft_2d_description}
We present the two-dimensional (2D) extension of the PFT algorithm. Analogous to the 1D case, we begin by recalling the 2D DFT given by:
\begin{equation}
    \hat{z}_{t_1, t_2} = \sum_{(k_1, k_2) \in [n_1]\times[n_2]} z_{k_1, k_2}e^{-2\pi it_1k_1/n_1}e^{-2\pi it_2k_2/n_2}
\end{equation}
where now $z \in \CC^{n_1 \times n_2}$ is a complex-valued matrix of size $n_1 \times n_2$. We assume that $n_1 = p_1q_1$ and $n_2 = p_2q_2$ are composite integers where $p_1, p_2, q_1, q_2 > 1$. Then rearranging the above expression, we get:
\begin{equation}
\label{eq: 2d_modified_summation}
    \hat{z}_{t_1, t_2} = \sum_{k_1, k_2, j_1, j_2} z_{q_1k_1 + j_1, q_2k_2 + j_2} \prod_{\nu} e^{-2\pi it_{\nu}(j_\nu - q_\nu/2)/n_\nu} \cdot e^{-2\pi it_\nu k_\nu/p_\nu} \cdot e^{-\pi i t_\nu/p_\nu}
\end{equation}
where $k_1 \in [p_1], k_2 \in [p_2], j_1 \in [q_1], j_2 \in [q_2]$, and $\nu = 1, 2$. As before, we want to use polynomial approximations for the exponential $e^{\pi ix}$. Afterwards, by re-scaling the polynomials and using exponent laws, one can get an approximation of the twiddle factors in the collection $\left\{e^{-2\pi it_\nu(j_\nu - q_\nu/2)/n_\nu}\right\}_{j_\nu = 0}^{q_\nu - 1}$. Using the same definitions and notation as in Section~\ref{subsec: pft_offline_derivation}, choose the minimum $r_\nu$ satisfying $\xi(\varepsilon, r_\nu) \geq \widetilde{m}_\nu/p_\nu$ to get the re-scaled polynomial approximations $\left\{\cP_{r_\nu - 1, \xi(\varepsilon, r_\nu)}(-2t_\nu(j_\nu - q_\nu/2)/n_\nu)\right\}_{j_\nu = 0}^{q_\nu - 1}$. Algorithm~\ref{alg: 2d_pft_offline} shows how to build the polynomial approximation. The polynomial coefficients $w_{\varepsilon, r_\nu - 1, j}$ are precomputed, and we obtain them from the code database of~\cite{park2021fast}. 
\begin{algorithm}[H]
    \caption{: Configuration (Offline) Phase of 2D PFT}\label{alg: 2d_pft_offline}
    \hspace*{\algorithmicindent} \textbf{Input}: Input size $(n_1, n_2) \in \NN^{2}$, crop size $(\widetilde{m}_1, \widetilde{m}_2) \in \NN^{2}$, divisors $(p_1, p_2) \in \NN^{2}$, and tolerance $\varepsilon$ \\
    \hspace*{\algorithmicindent} \textbf{Output}: Matrices $B_1 \in \CC^{q_1 \times r_1}, B_2 \in \CC^{q_2 \times r_2}$, tensor $W \in \mathbb{C}^{(2\widetilde{m}_1 + 1) \times (2\widetilde{m}_2 + 1) \times r_1 \times r_2}$, configuration results $p_1, p_2, q_1, q_2, r_1, r_2$
    \begin{algorithmic} [1]
    \For{$\nu = 1, 2$}
        \State $q_\nu = n_\nu/p_\nu$
        \State $r_\nu = \min\{r_\nu \in \NN : \xi(\varepsilon, r_\nu) \geq \widetilde{m}/p\}$
        \Comment{degree of polynomial $\cP$ approximating $e^{\pi ix}$ within tolerance $\varepsilon$}
        \For{$l \in [q_\nu], j \in [r_\nu]$}
            \State $x = (1 - 2l/q_\nu)$
            \State $B_{\nu} = w_{\varepsilon, r_\nu - 1, j} \cdot x^{j}$ \Comment{Using precomputed $w_{\varepsilon, r_{\nu}-1, j}$}
        \EndFor
    \EndFor
    \For{$(t_1, t_2) = \left([-\widetilde{m}_1, \ldots,  \widetilde{m}_1] \times [-\widetilde{m}_2, \ldots, \widetilde{m}_2]\right)$}
        \State $W[t_1, t_2, j_1, j_2] = (t_1/p_1)^{j_1} e^{-\pi it_1/p_1} \cdot (t_2/p_2)^{j_2} e^{-\pi it_2/p_2}$
        \Comment{Precompute remaining terms}
    \EndFor
    \end{algorithmic}
\end{algorithm}

Substituting the approximating polynomials for the twiddle factors in~\eqref{eq: 2d_modified_summation} and performing some algebraic manipulations (in particular, swapping and rewriting summations), we represent the summations as matrix-matrix multiplications between $C^{(k_1, k_2)} = B_{1}^{\top} \times Z^{(k_1, k_2)} \times B_{2}$ where $Z^{(k_1, k_2)} \in \CC^{q_1 \times q_2}$ are slices of the data matrix $z$, the matrices $B_{1} \in \CC^{q_1 \times r_1}$ and $B_{2} \in \CC^{q_2 \times r_2}$ are defined in line $6$ of Algorithm~\ref{alg: 2d_pft_offline} followed by a series of 2D FFT computations on $C$ for each $(k_1, k_2) \in [p_1] \times [p_2]$. As in the 1D PFT, we note that the primary cost of Algorithm~\ref{alg: 2d_pft_online} lies in the matrix multiplications between $Z^{(k_1, k_2)}$, $B_1$, and $B_2$ and the 2D FFT computations on each of the matrices $C^{(k_1, k_2)} \in \CC^{r_1 \times r_2}$. This gives us time complexity $\cO(n + \widetilde{m}\log\widetilde{m})$, where $n = n_1n_2$ and $\widetilde{m} = \widetilde{m}_1\widetilde{m}_2$. For thorough details, we refer the reader to~\cite[Appendix B]{park2021fast}.
\begin{algorithm}[H]
    \caption{: Computation (Online) Phase of 2D PFT}\label{alg: 2d_pft_online}
    \hspace*{\algorithmicindent} \textbf{Input}: 2D array $Z$ of size $n_1 \times n_2$, crop size $(\widetilde{m}_1, \widetilde{m}_2)$, tensor $W \in \mathbb{C}^{(2\widetilde{m}_1 + 1) \times (2\widetilde{m}_2 + 1) \times r_1 \times r_2}$, configuration results $B_1 \in \CC^{q_1 \times r_1}, B_2 \in \CC^{q_2 \times r_2}, p_1, p_2, q_1, q_2, r_1, r_2$ \\
    \hspace*{\algorithmicindent} \textbf{Output}: 2D array $\hat{Z}_{PFT}$ of estimated Fourier coefficients of $Z$
    \begin{algorithmic} [1]
        \State $Z = Z\text{.reshape}(p_1, q_1, p_2, q_2)$
        \State $Z = Z\text{.permute}(0, 2, 1, 3)\text{.contiguous()}$
        \State $Z = Z\text{.reshape}(p_1, p_2, q_1, q_2)$ \Comment{reshape $z$ into $p_1 \times p_2 \times q_1 \times q_2$ tensor}
        \For{$(k_1, k_2) \in [p_1] \times [p_2]$}
            \State $C[k_1, k_2, :, :] = B_{1}^{\top} \times Z[k_1, k_2, :, :] \times B_{2}$
            \Comment{matrix multiply $B_{1}^{\top}$ by $Z[k_1, k_2, :, :]$ by $B_{2}$}
        \EndFor
        \For{$(j_1, j_2) \in [r_1] \times [r_2]$}
            \State $\hat{C}[:, :, j_1, j_2] = $ FFT2$(C[:, :, j_1, j_2])$
            \Comment{apply 2D FFT to matrices $C[:, :, j_1, j_2]$}
        \EndFor
        \For{$(t_1, t_2) = \left([-\widetilde{m}_1, \ldots,  \widetilde{m}_1] \times [-\widetilde{m}_2, \ldots, \widetilde{m}_2]\right)$}
            \State $\hat{Z}_{PFT}[t_1, t_2] = \sum_{j_1 \in [r_1], j_2 \in [r_2]} \hat{C}[t_1\%p_1, t_2\%p_2, j_1, j_2] \cdot W[t_1, t_2, j_1, j_2]$
            \Comment{Hadamard product}
        \EndFor
    \end{algorithmic}
\end{algorithm}

\section{Illumination Windows For Large-Scale Blind Ptychography}
We provide an illustration of the illumination window for the large-scale experiment. The size of the image is $8200 \times 8200$, each circular probe has a radius of $2553$ pixels, and shifts $2188$ pixels at a time.
\begin{figure}[H]
    \setlength\tabcolsep{1 pt}
    \centering
    \begin{tabular}{cccccccc}
        $Q_1$ & $Q_2$ & $Q_3$ & $Q_4$ & $Q_5$ & $Q_6$ & $Q_7$ & $Q_8$
        \\
        \includegraphics[width=0.12\textwidth]{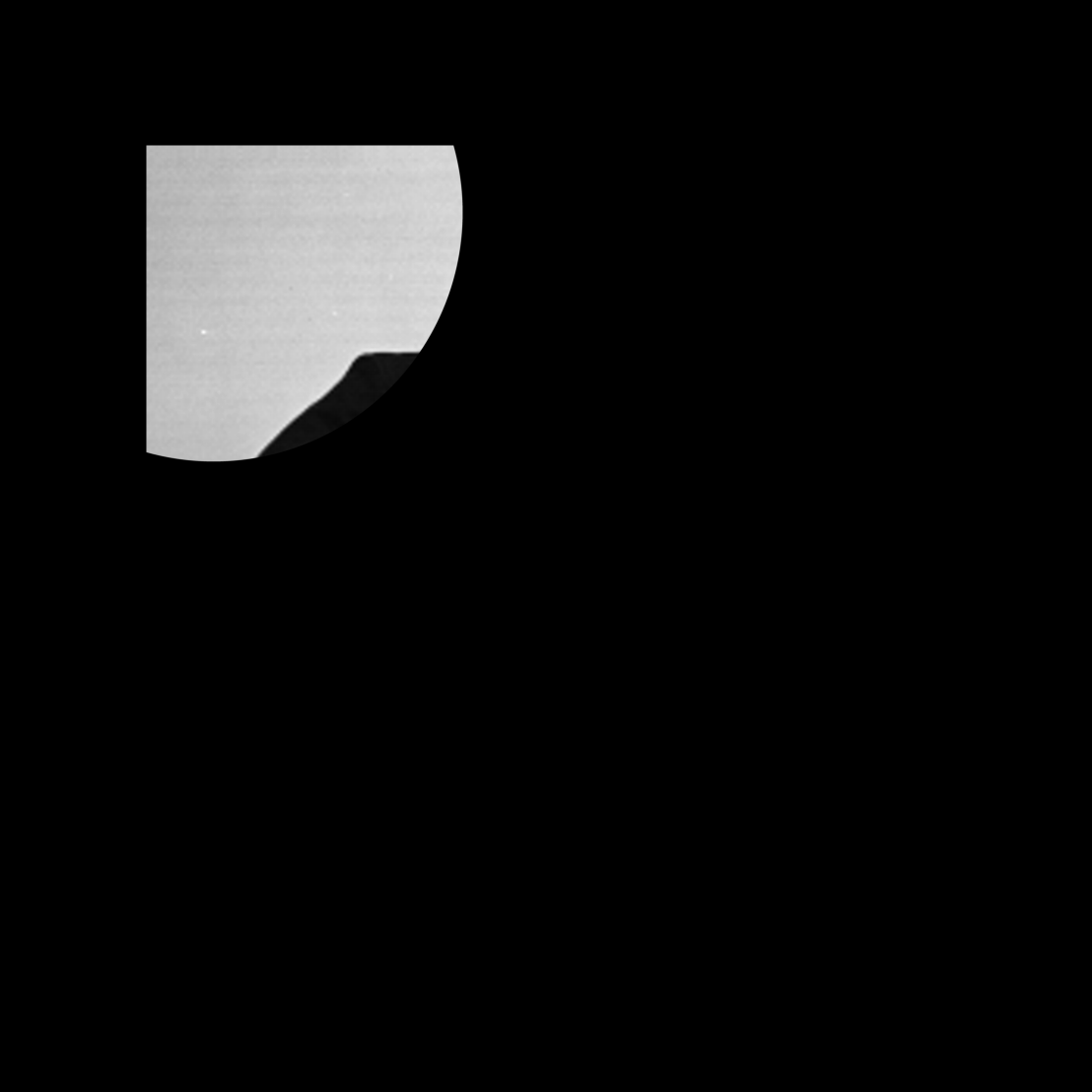}
        &
        \includegraphics[width=0.12\textwidth]{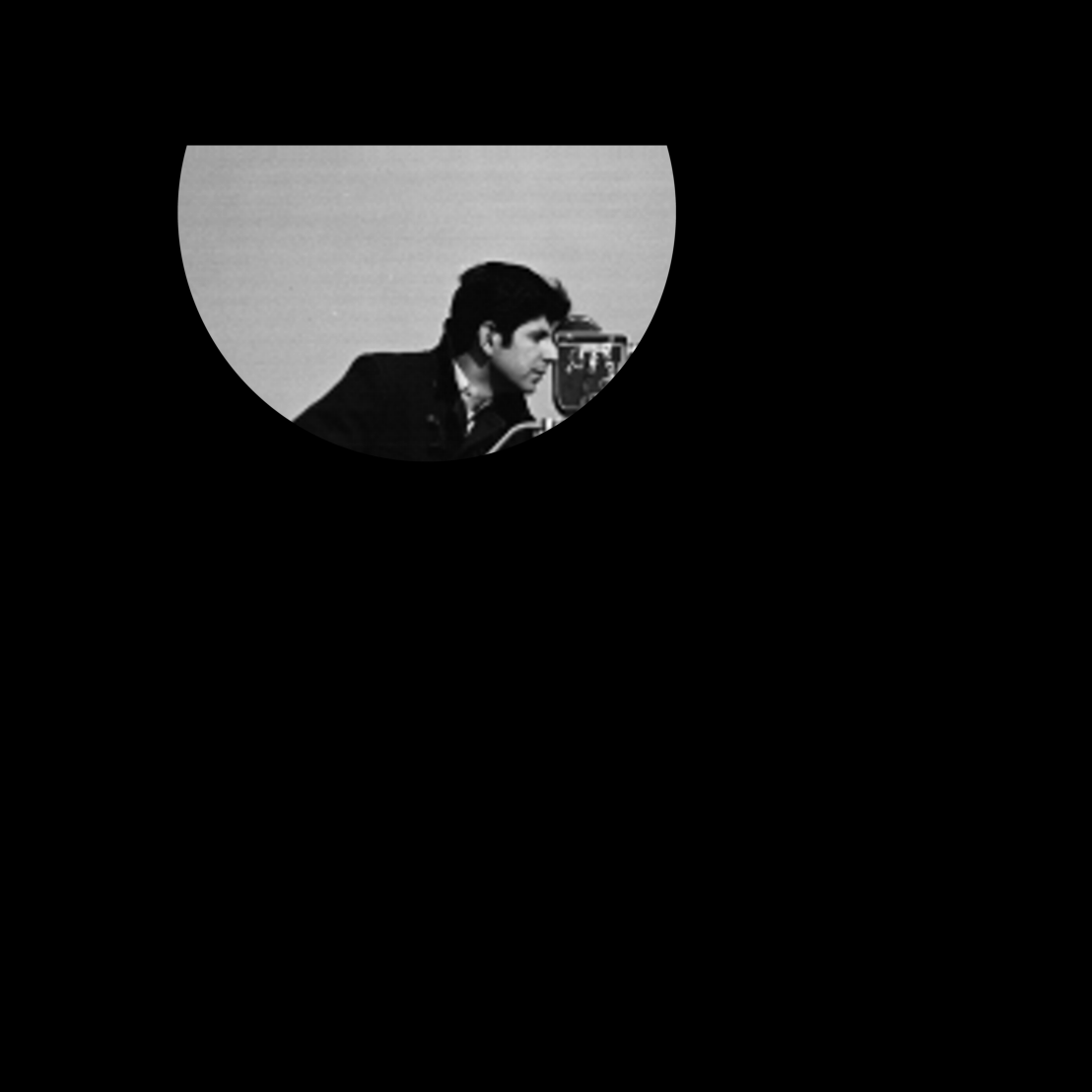}
        &
        \includegraphics[width=0.12\textwidth]{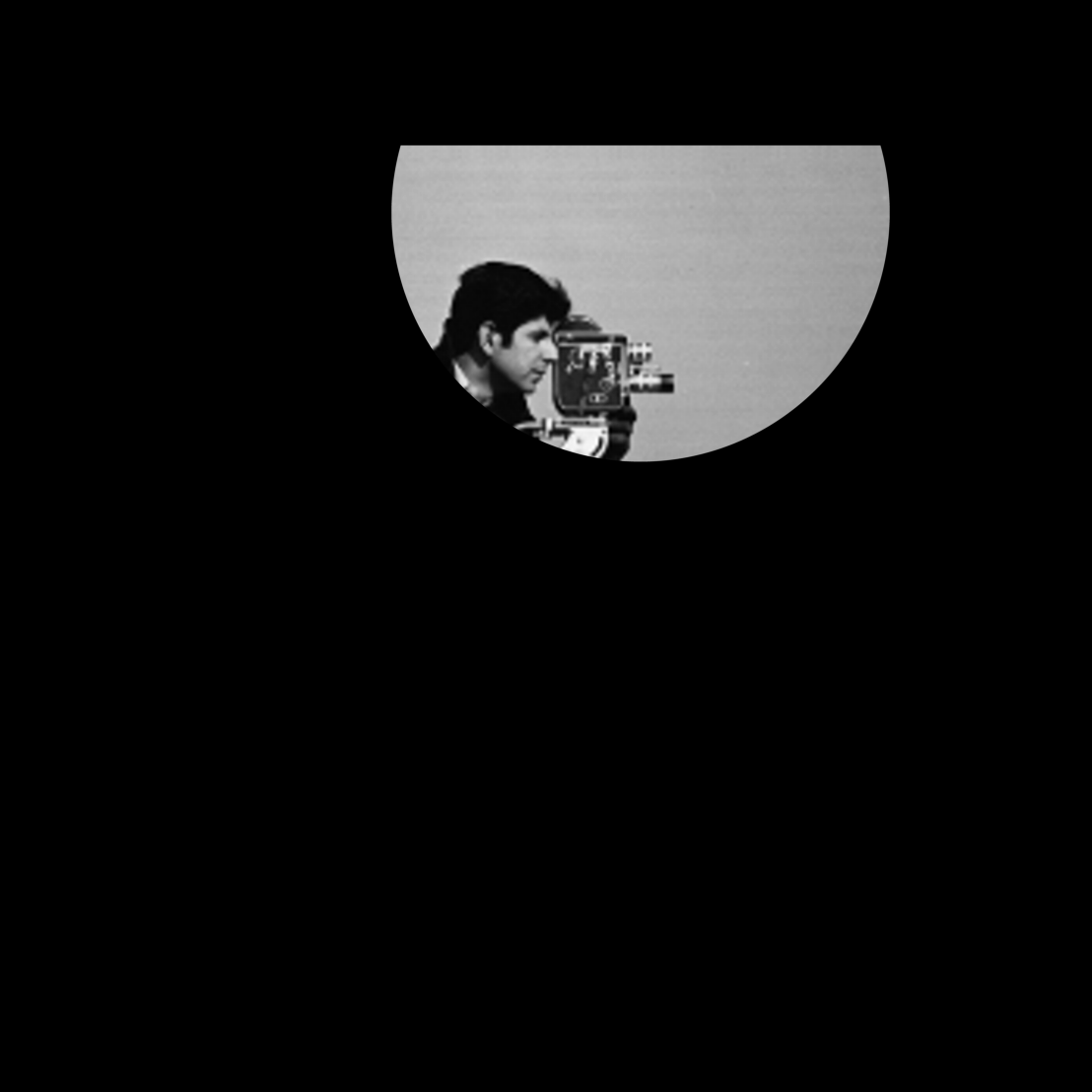}
        &
        \includegraphics[width=0.12\textwidth]{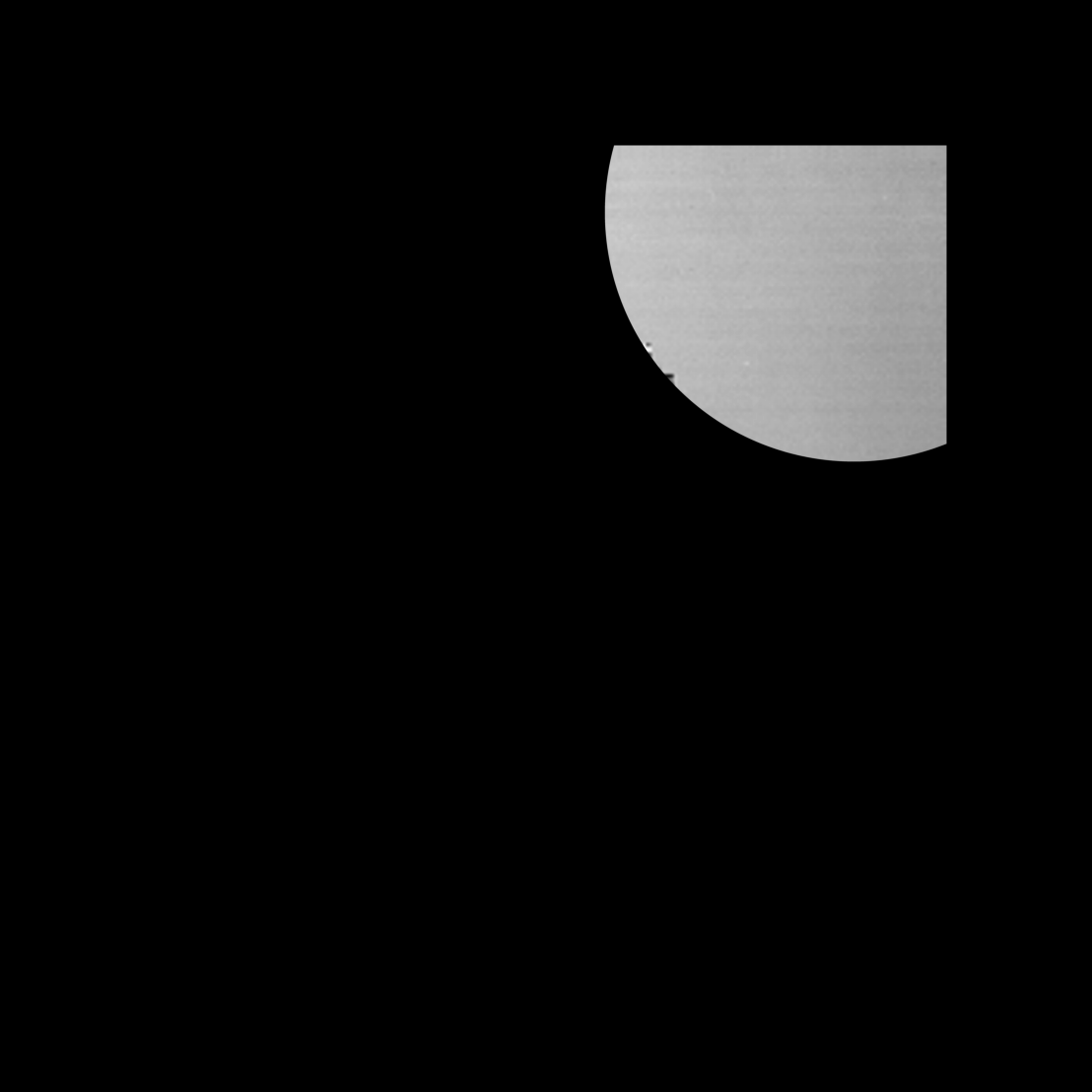} 
        &
        \includegraphics[width=0.12\textwidth]{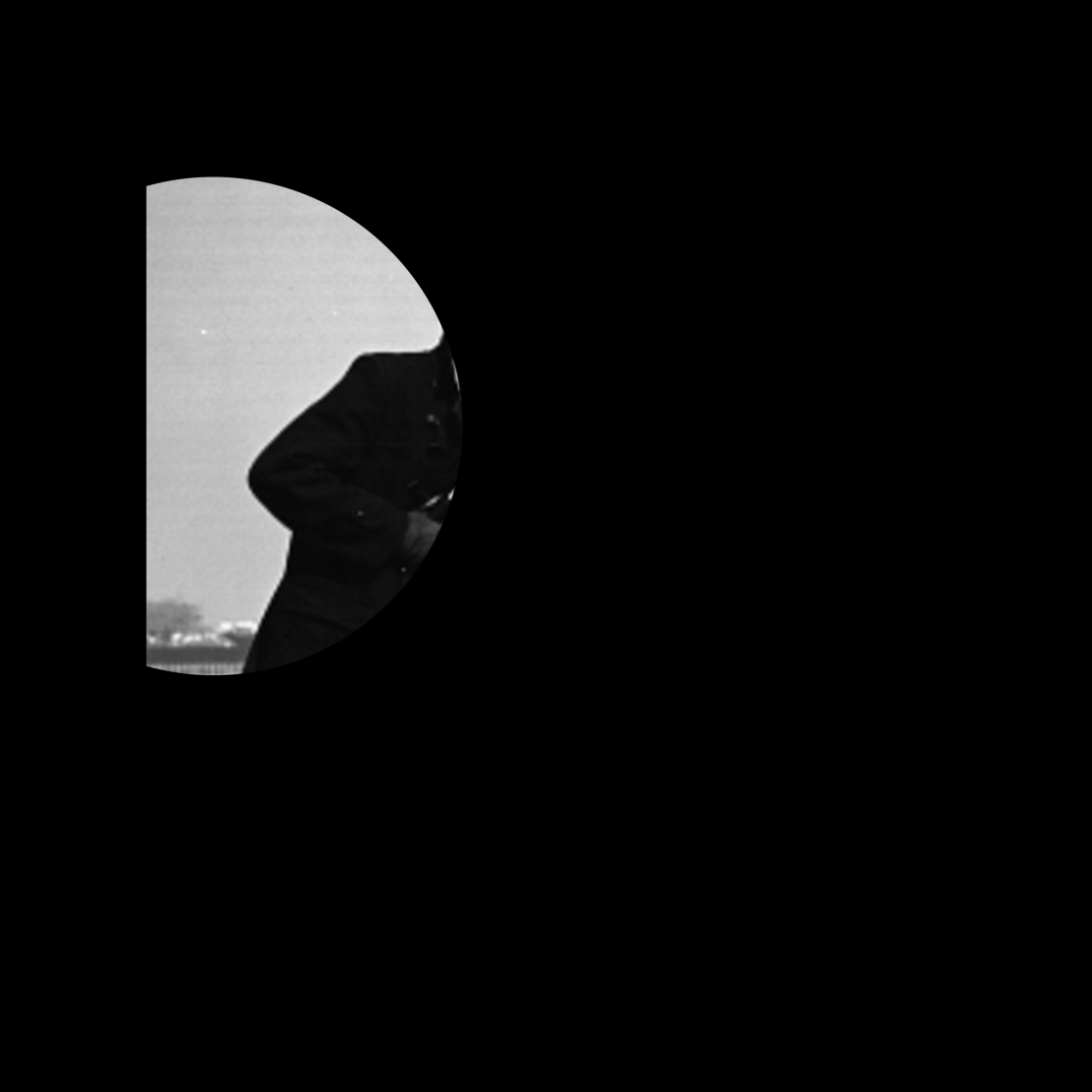}
        &
        \includegraphics[width=0.12\textwidth]{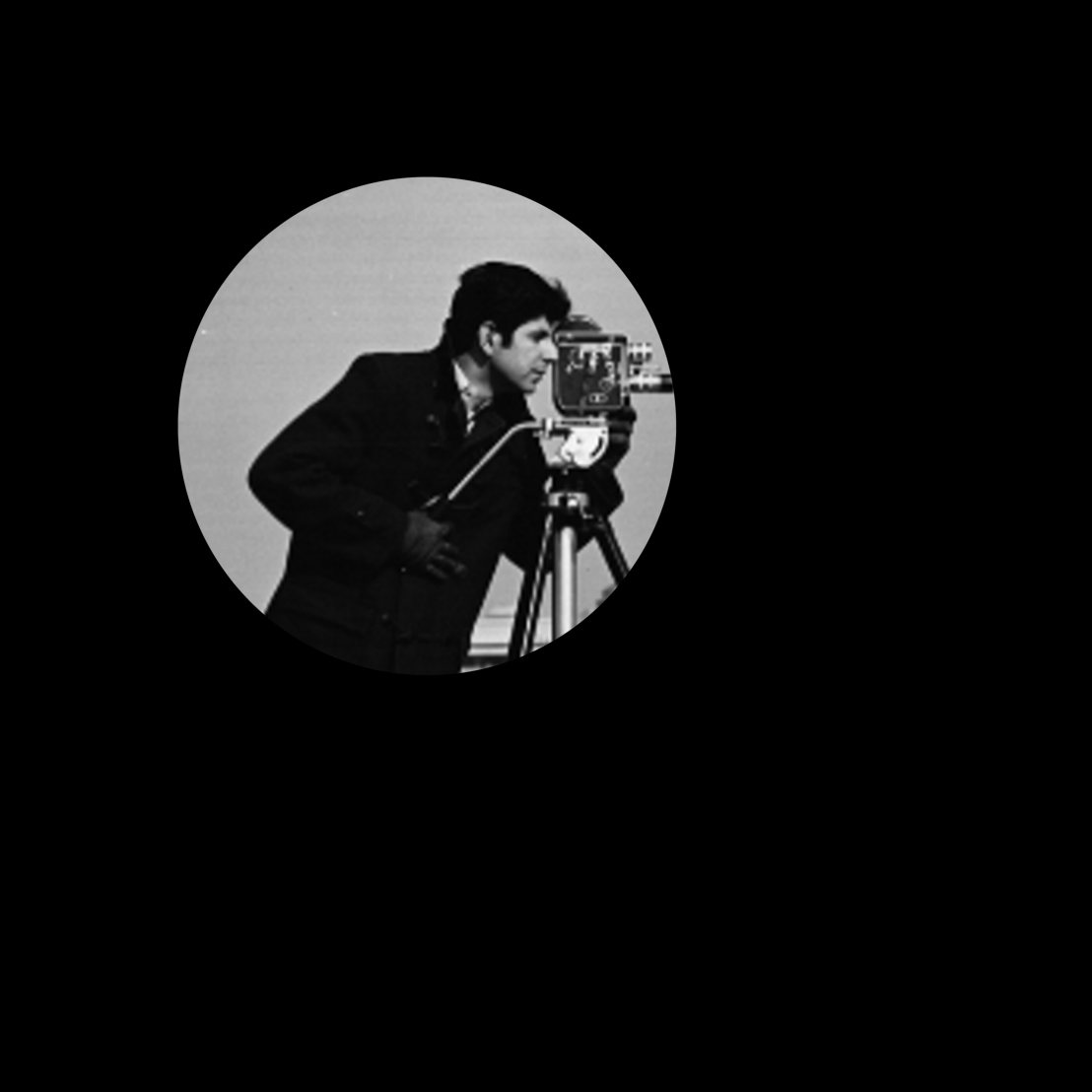}
        &
        \includegraphics[width=0.12\textwidth]{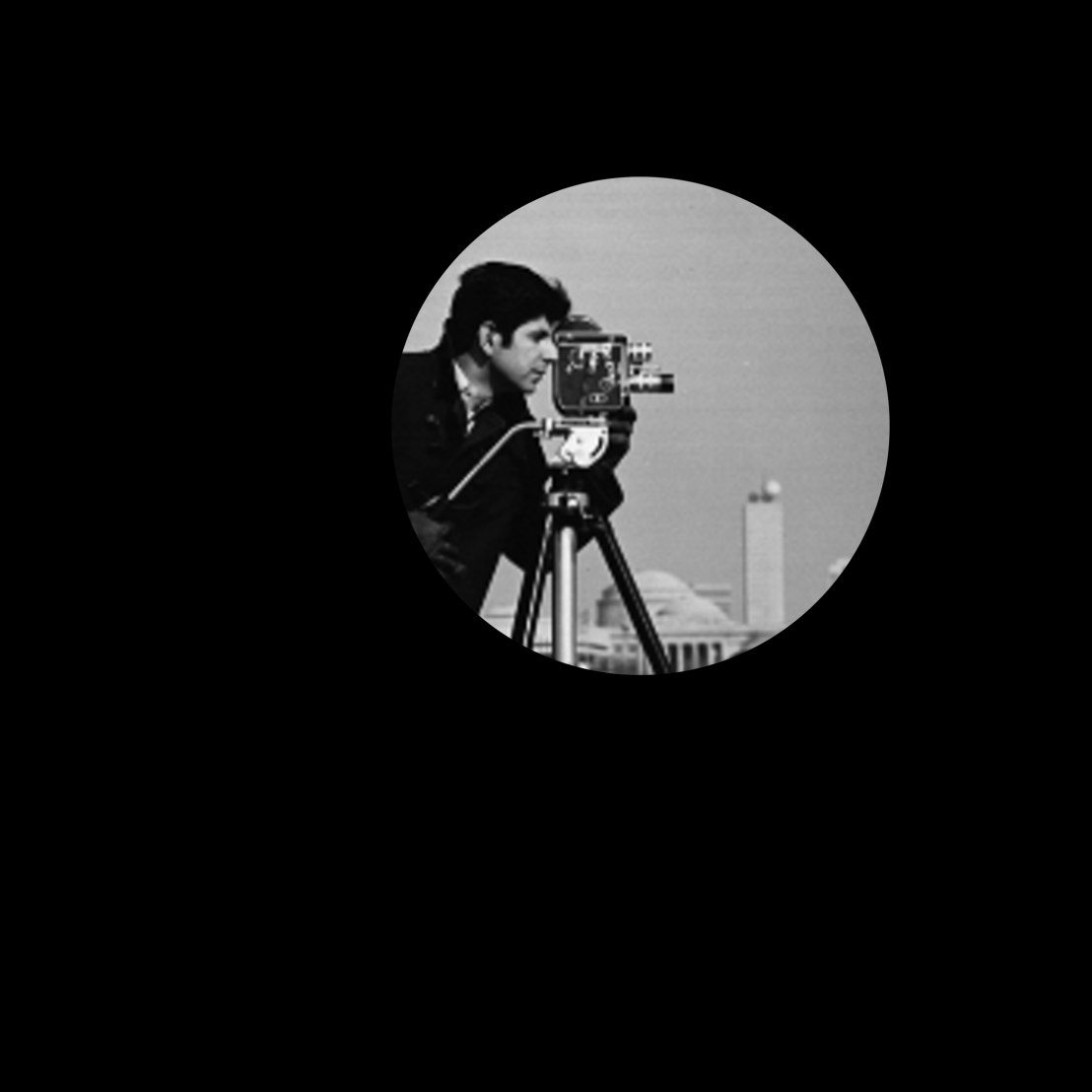}
        &
        \includegraphics[width=0.12\textwidth]{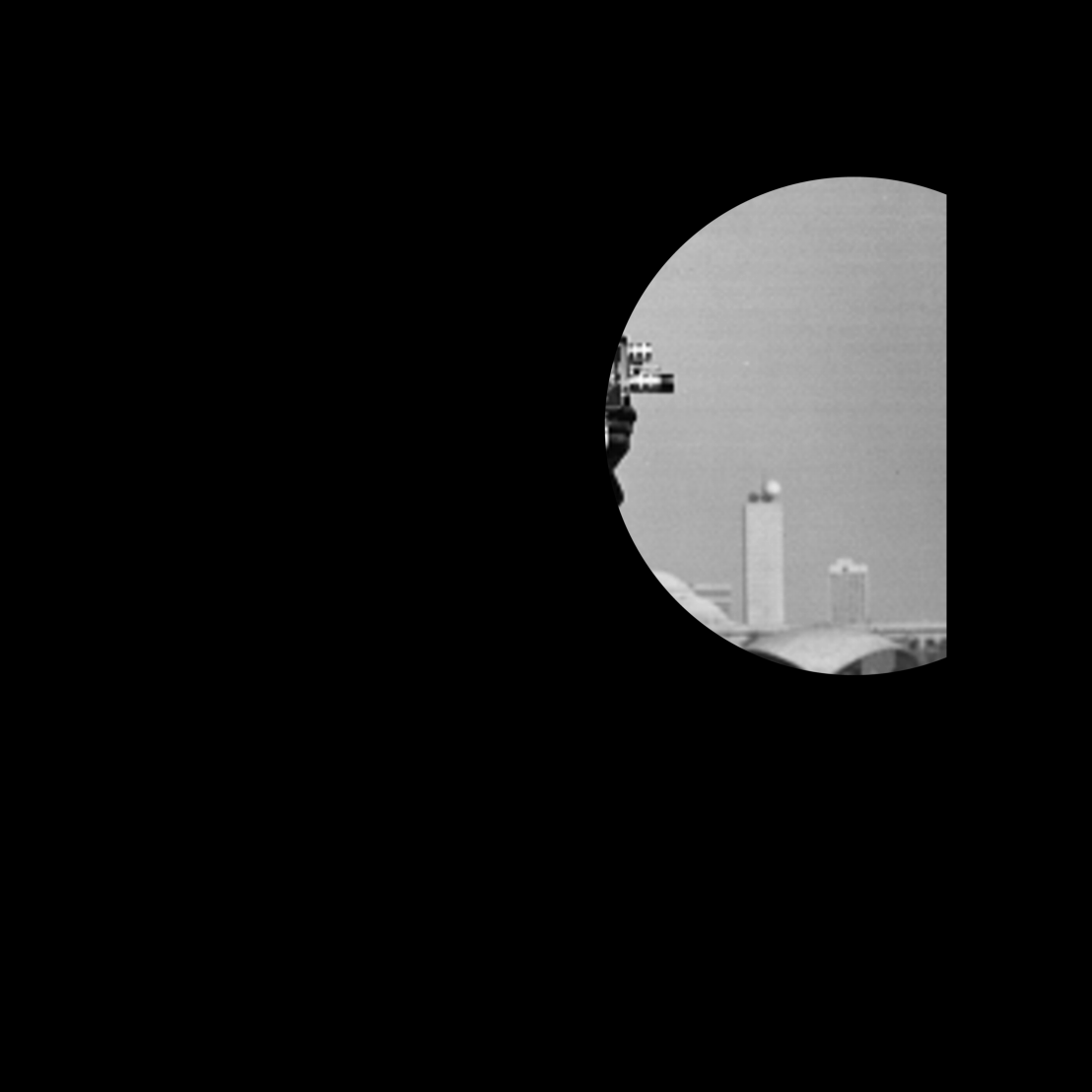} \\
        $Q_9$ & $Q_{10}$ & $Q_{11}$ & $Q_{12}$ & $Q_{13}$ & $Q_{14}$ & $Q_{15}$ & $Q_{16}$
        \\
        \includegraphics[width=0.12\textwidth]{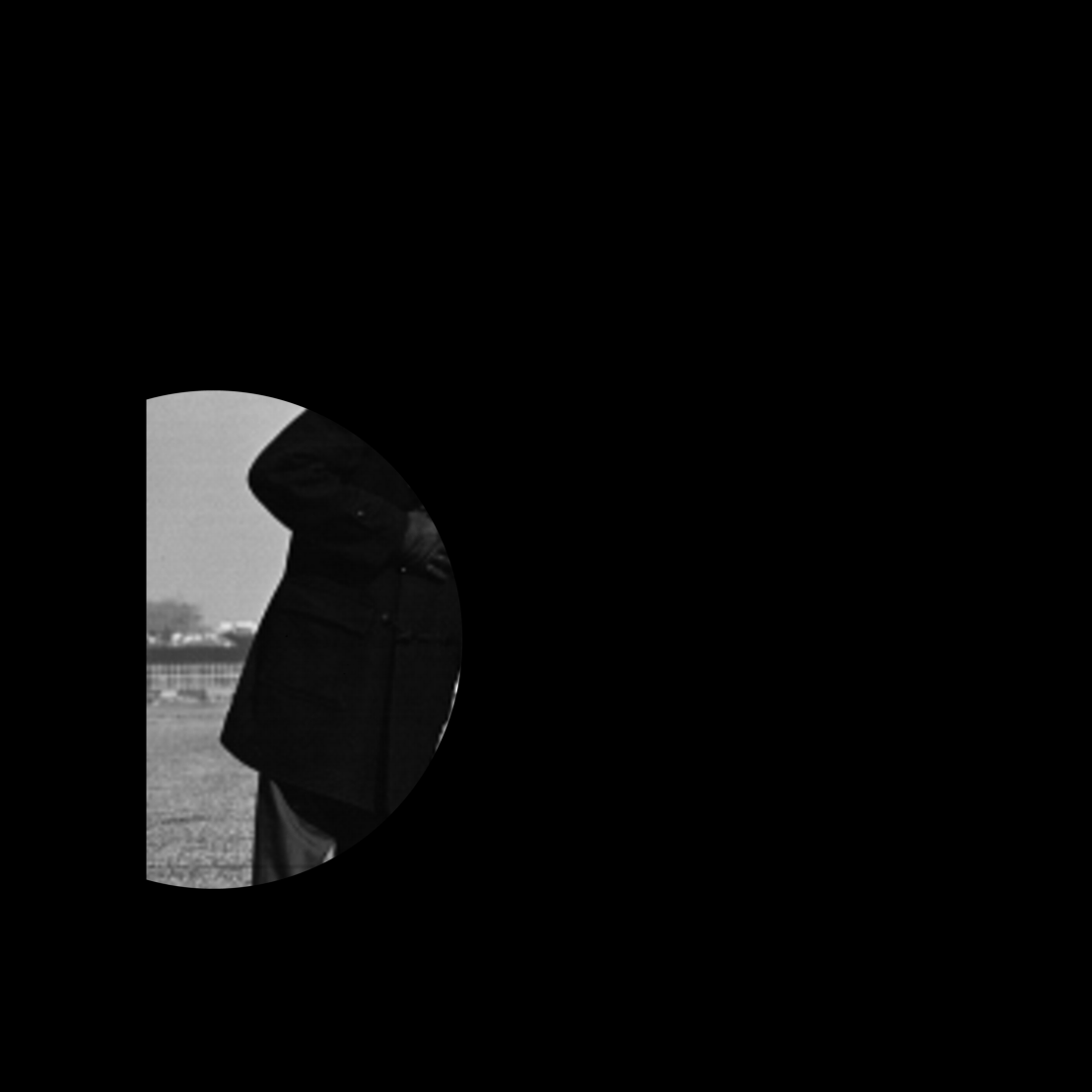}
        &
        \includegraphics[width=0.12\textwidth]{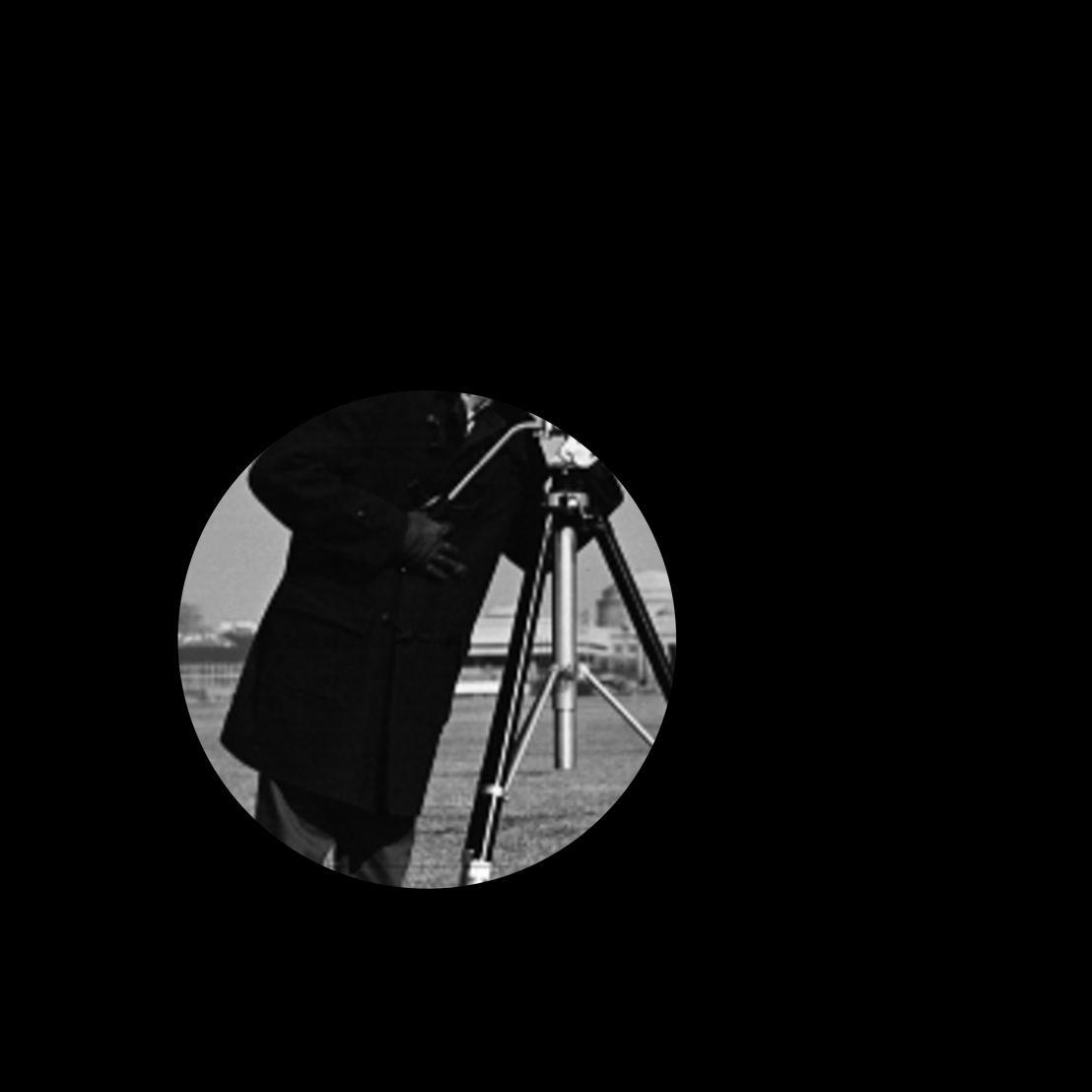}
        &
        \includegraphics[width=0.12\textwidth]{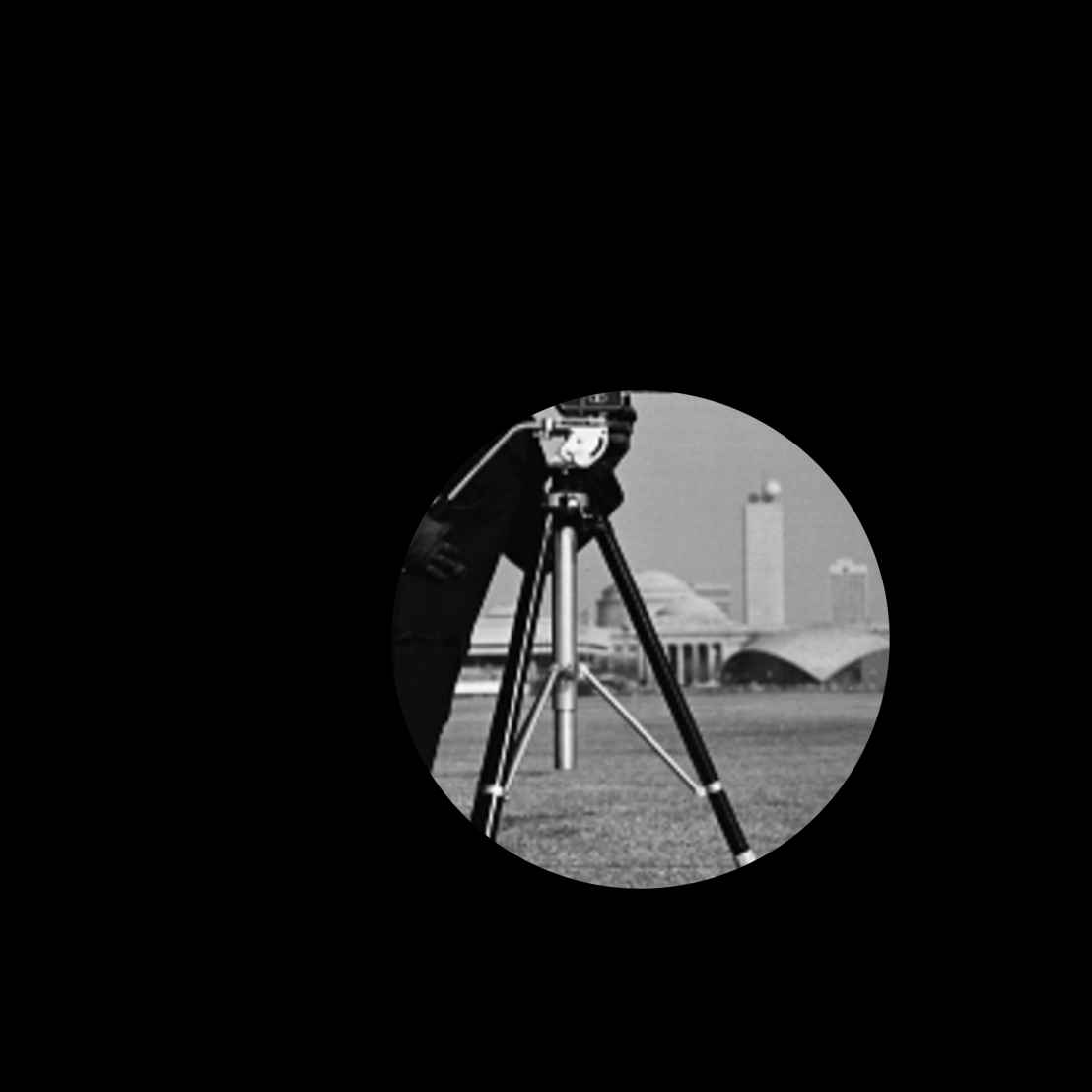}
        &
        \includegraphics[width=0.12\textwidth]{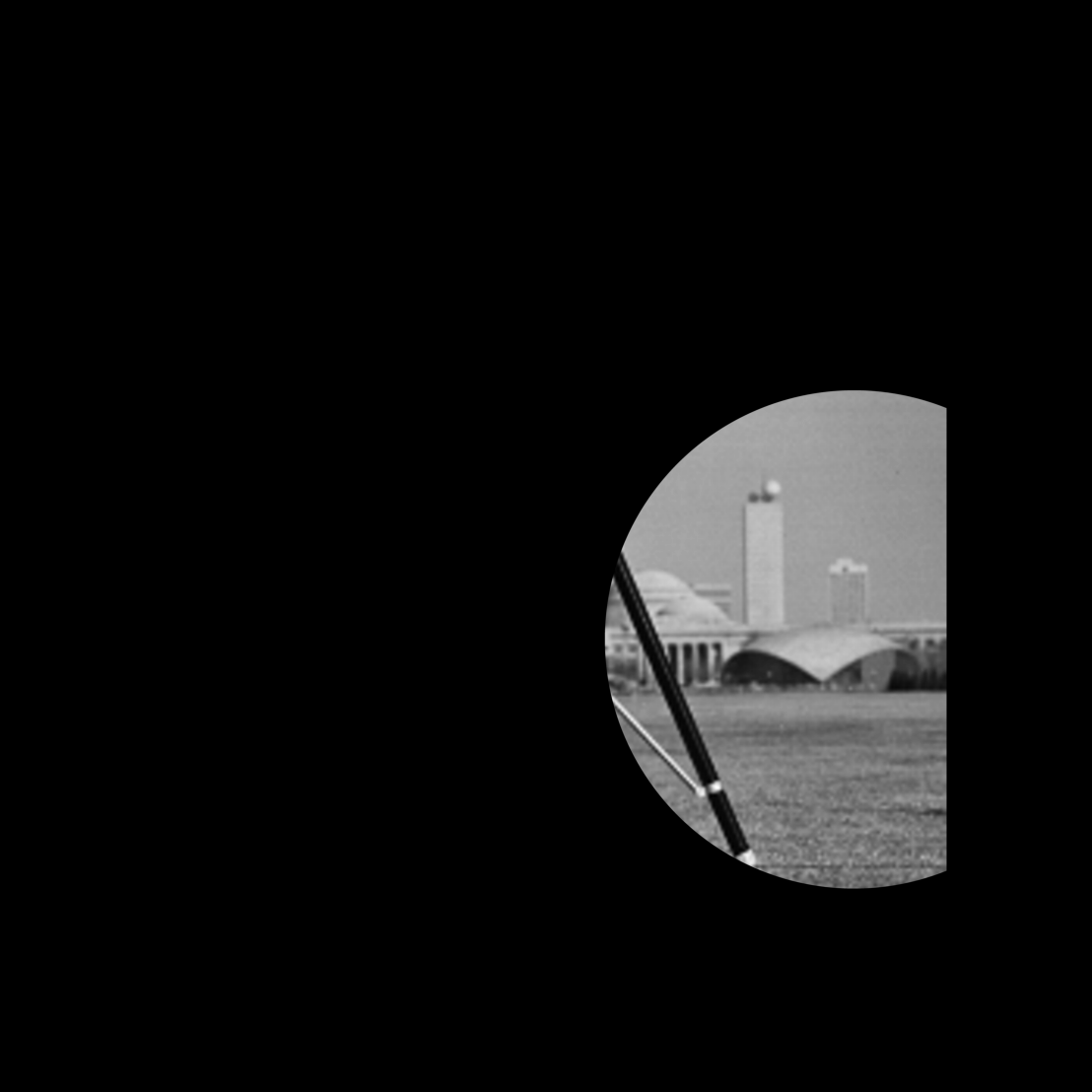}
        &
        \includegraphics[width=0.12\textwidth]{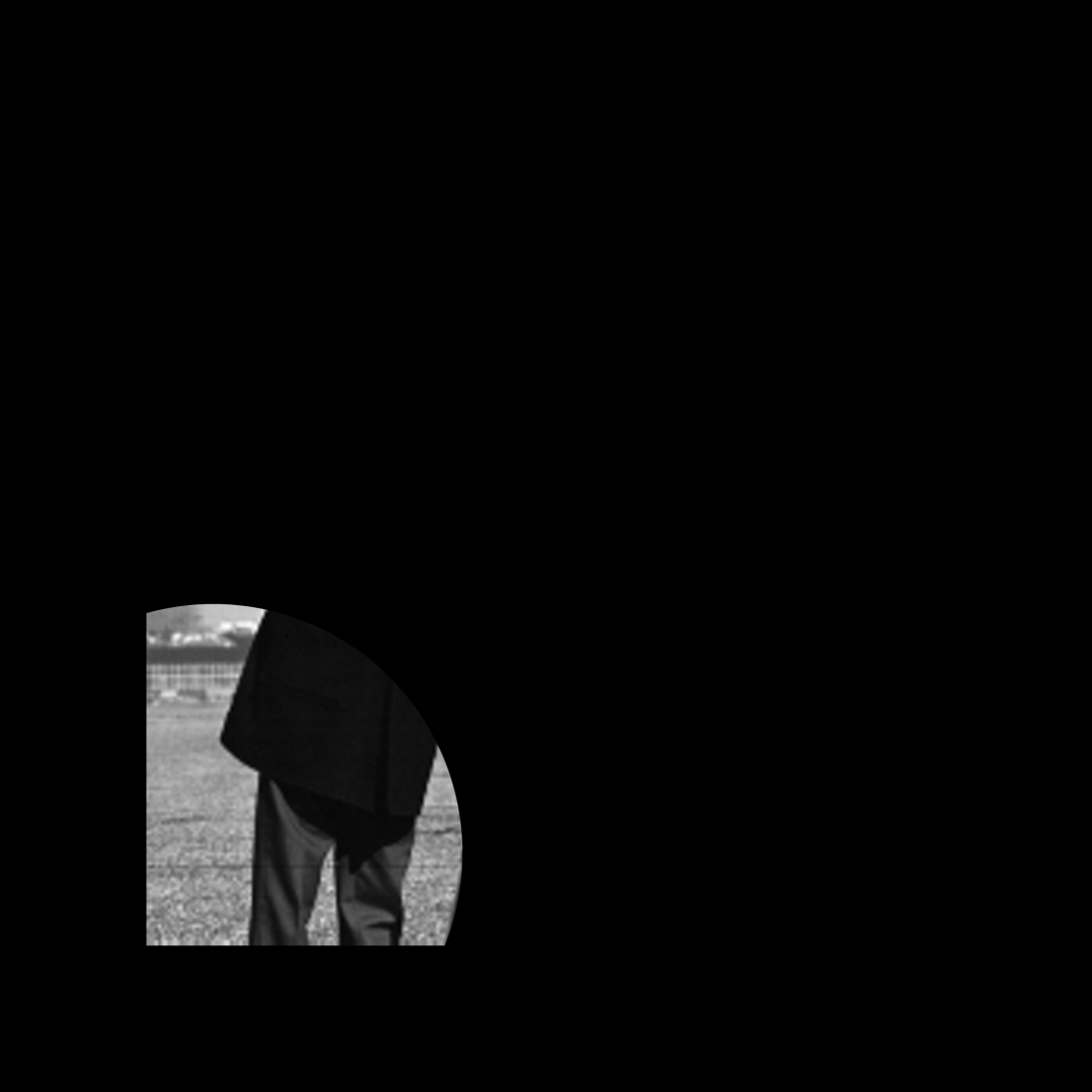}
        &
        \includegraphics[width=0.12\textwidth]{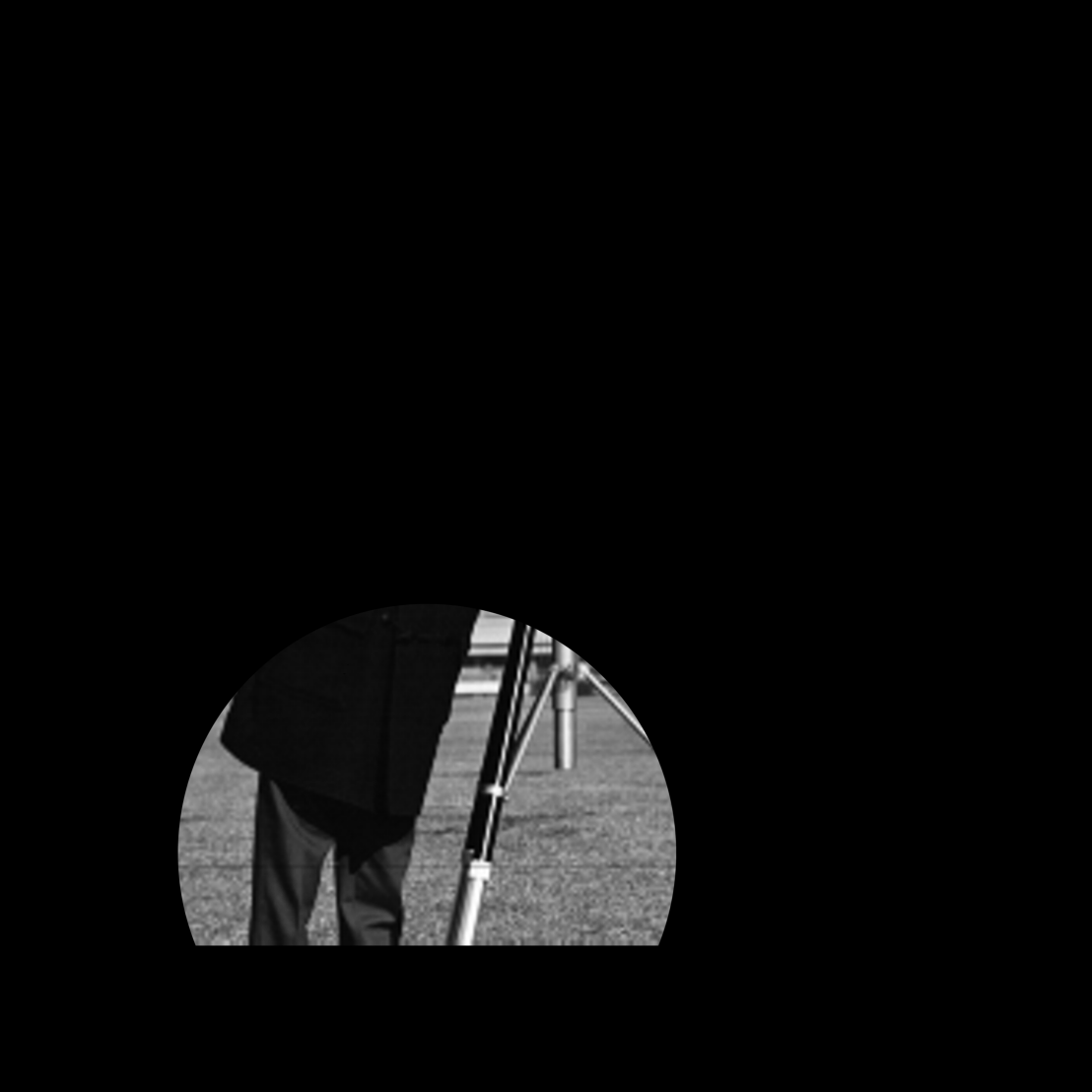}
        &
        \includegraphics[width=0.12\textwidth]{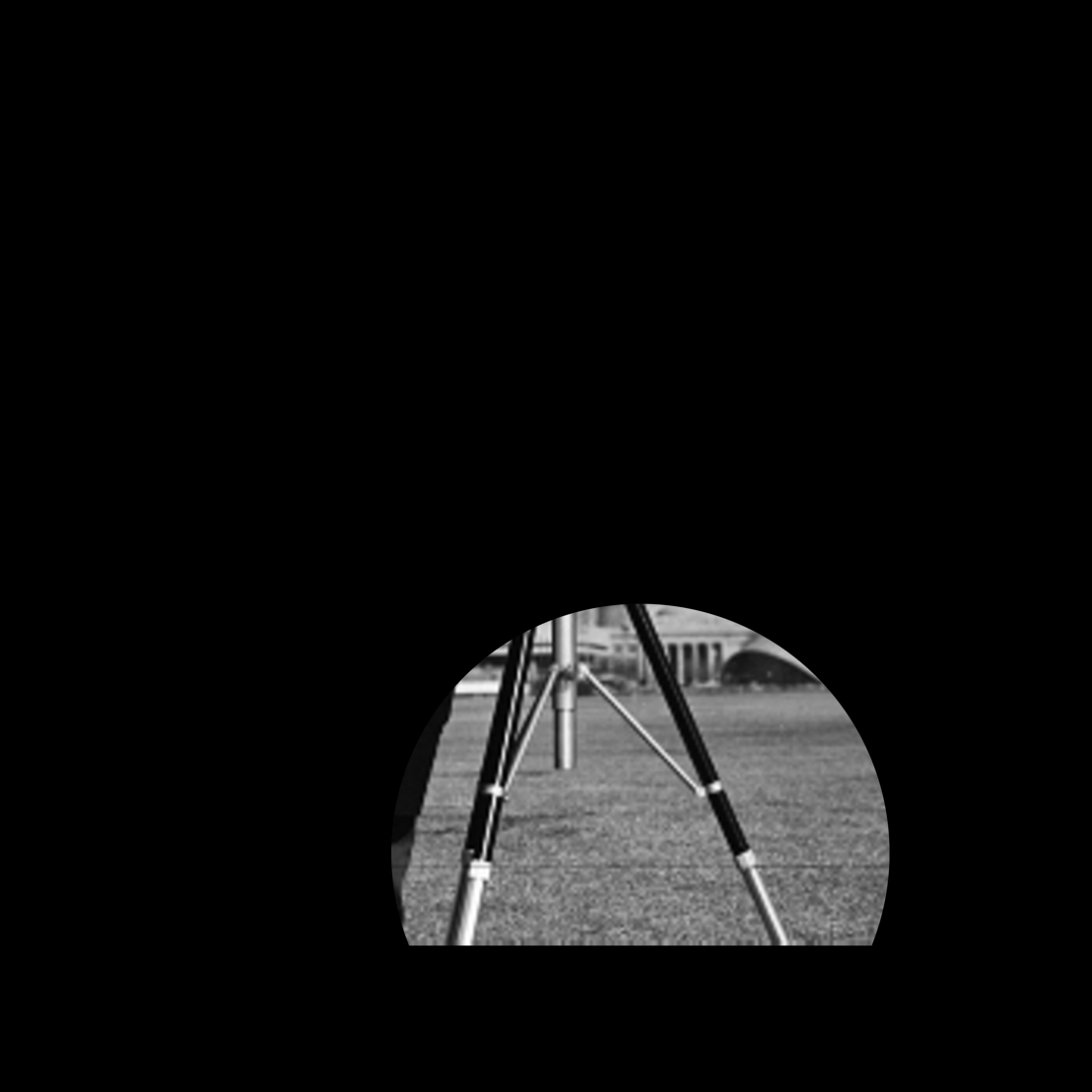}
        &
        \includegraphics[width=0.12\textwidth]{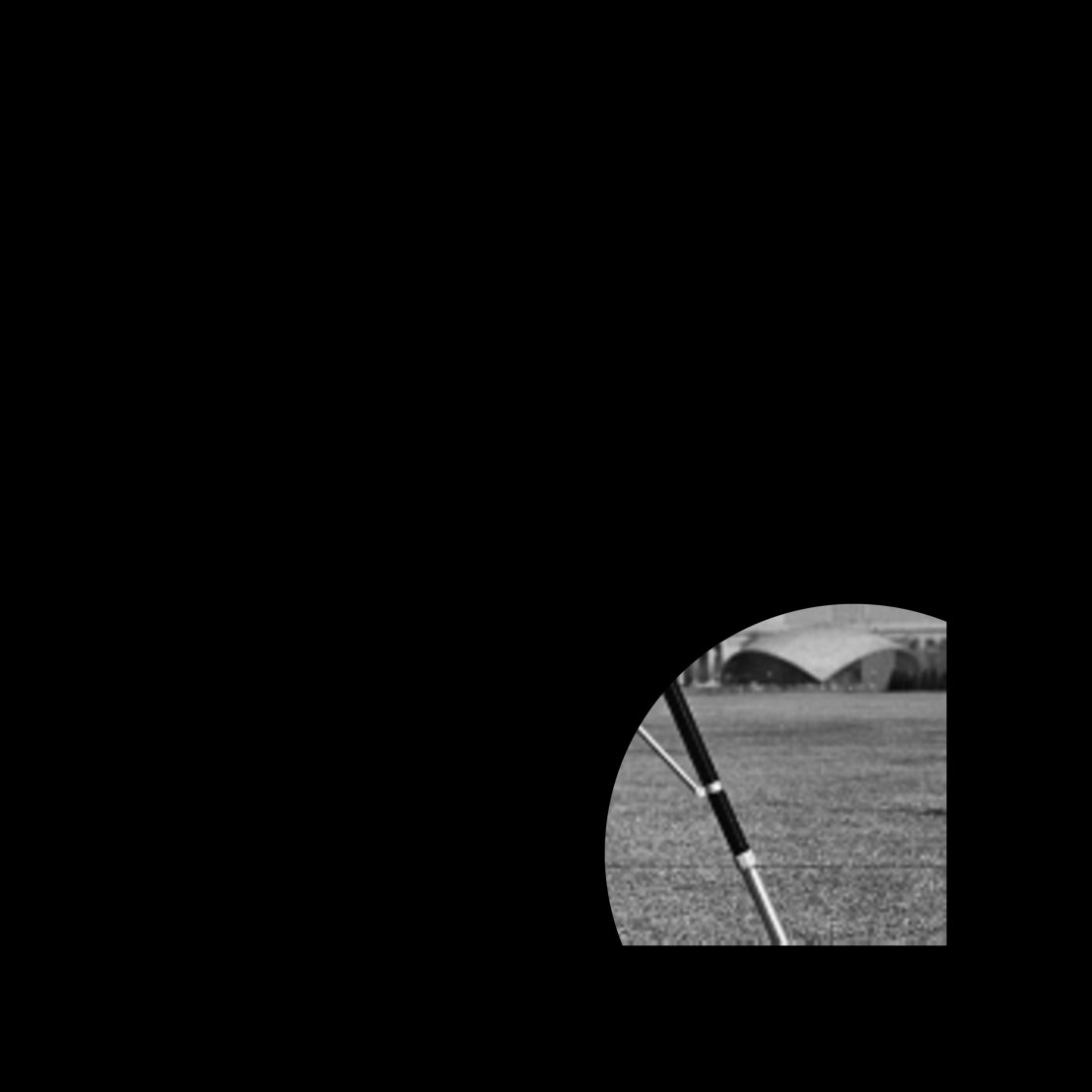}
    \end{tabular}
    \caption{\small{Scanning positions of the illumination window described in Section~\ref{subsubsec: blind_experimental_setup}}. Here, each probe illuminates a circle of radius $2553$ pixels and shifts $2188$ pixels at a time, leading to roughly $50\%$ overlap between consecutive probes.}
    \label{fig: circular_probes_large}
\end{figure}
\end{document}